\documentclass[11pt,reqno]{amsart}
\usepackage{amssymb}
\usepackage{bm}
\usepackage{xcolor}
\usepackage{amsmath}
\usepackage{enumitem}
\usepackage{mathrsfs}
\usepackage{hyperref}
\usepackage{titlesec}
\usepackage{tikz-cd}
\usepackage{tikz}
\usetikzlibrary{calc}
\usetikzlibrary{shapes.multipart}
\usetikzlibrary{patterns}

\titleformat{\section}  
{\fontsize{14}{16}\bfseries} 
{\thesection.} 
{1em} 
{} 
[] 
\titlespacing*{\section}{0pt}{1.5em}{0.8em}

\titleformat{\subsection}  
{\fontsize{13}{15}\bfseries} 
{\thesubsection.} 
{1em} 
{} 
[] 
\titlespacing*{\subsection}{0pt}{1.2em}{0.6em}

\newtheorem{theorem}{Theorem}[section]
\newtheorem{lemma}[theorem]{Lemma}

\newtheorem{proposition}[theorem]{Proposition}
\newtheorem{corollary}[theorem]{Corollary}
\newtheorem{claim}[theorem]{Claim}

\theoremstyle{definition}
\newtheorem{definition}[theorem]{Definition}
\newtheorem{notation}[theorem]{Notation}
\newtheorem{remark}[theorem]{Remark}

\theoremstyle{definition}
\newtheorem{example}[theorem]{Example}

\newenvironment{claimproof}[1]{\par\noindent\textit{Proof of claim.}\space#1}{\leavevmode\unskip\penalty9999 \hbox{}\nobreak\hfill\quad\hbox{$\square$}\vspace{0.3cm}}

\let \restr = \upharpoonright

\let \vect = \overrightarrow

\let \la = \langle
\let \ra = \rangle

\DeclareMathOperator{\dom}{dom}

\DeclareMathOperator{\cf}{cf}
\DeclareMathOperator{\ot}{ot}

\DeclareMathOperator{\rank}{rank}

\title{Forcing with symmetric systems of models of two types}

\author[Curial Gallart]{Curial Gallart}

\address{Curial Gallart, Institute of Mathematics, University of Vienna, Kolingasse 14-16, 1090 Vienna, Austria}

\email{curial.gallart.rodriguez@univie.ac.at}

\date{}

\begin{document}
	
	
	
	\thanks{This work was partially supported by the EPSRC project reference 2440142 and the Austrian and Polish Science Foundations NCN-FWF through the project PIN1355423. Part of this work is based on the author's PhD thesis \cite{Gallart2024thesis}. The author would like to express his most sincere gratitude to David Asperó for his encouragement and for the very fruitful discussions about this topic. The author would also like to thank Monroe Eskew for encouraging him to expand this work, which resulted in Sections \ref{section-decorations} and \ref{section-morass}. The author is also very grateful to the anonymous referees for their careful reading of the manuscript and their very insightful suggestions, which resulted in a much better presentation of the paper.}
	
	\begin{abstract}
		The purpose of this paper is to present a general method for forcing on $\omega_2$ and $\omega_3$ with finite conditions, while preserving all cardinals and some fragments of $\mathrm{GCH}$. This method is based on the technique of forcing with finite symmetric systems of elementary submodels, and improves earlier versions of this forcing by including models of two types. We will present several applications of the pure side condition forcing and variants thereof, by adding a Kurepa tree on $\omega_2$, a club subset of $\omega_2$ that avoids infinite sets from the ground model, a function bounding every canonical function below $\omega_3$ on a club, and a simplified $(\omega_2,1)$-morass.
	\end{abstract}
	
	\maketitle
	\pagestyle{myheadings}\markright{Forcing with symmetric systems of models of two types}

	\section{Introduction}
		
	The idea of adding finite sets of elementary submodels in the conditions of a forcing notion, known as \emph{forcing with side conditions}, dates back to the work of Todor\v cevi\'c in \cite{Todorcevic1985directedsets} and \cite{Todorcevic1984AnoteonthePFA}. This technique has been exploited by many set theorists to build nice notions of forcing adding generic objects of size $\aleph_1$ with finite conditions (see \cite{Todorcevic2014notesonforcingaxioms} for many applications). The addition of finite systems of countable elementary submodels in the conditions ensures that the forcing is proper, and hence preserves the first uncountable cardinal.
	
	Typically, a notion of forcing with side conditions has two components: the \emph{working part}, which is a finite approximation of the object we are ultimately interested in adding in the generic extension, and the \emph{side condition}, which consists of a finite system of elementary submodels that interact in some suitable way with the working part. 
	
	These systems of models have a certain structure that determines the preservation properties of the forcing notion. In Todor\v cevi\'c's inital applications, the side conditions consisted of finite $\in$-chains of countable elementary submodels of some $H(\kappa)$, which ensured properness of the resulting poset, and hence preservation of $\aleph_1$. Mitchell \cite{Mitchell2005Addingclubsofomega2withfiniteforcing} and Friedman \cite{Friedman2006Forcingwithfiniteconditions} independently developed forcing notions with side conditions that add club subsets of $\omega_2$ with finite conditions, while preserving $\aleph_1$ and $\aleph_2$. Both posets used finite sets of countable elementary submodels, but the requirement of forming an $\in$-chain was replaced with more complicated agreement and coherence conditions between the models. Mitchell-style side conditions were also used by Mitchell himself \cite{Mitchell2009I[omega2]canbethenonstationaryidealoncofomega1} to prove that it is consistent relative to large cardinals that there is no stationary subset of $S^{\omega_2}_{\omega_1}$ in the approachability ideal $I[\omega_2]$, and by Dolinar and Dzamonja \cite{DolinarDzamonja2013} to force a square sequence on $\omega_2$ with finite conditions.
	
	Neeman \cite{Neeman2014Forcingwithsequencesofmodelsoftwotypes}, building on the works of Mitchell and Friedman, developed a common framework for forcing at the level of $\omega_2$ with finite conditions. His poset captured the main features of the posets of Mitchell and Friedman, while retaining much of the simplicity of Todor\v cevi\'c's side conditions. The breakthrough in Neeman's poset was to use models of two different types, countable and transitive, to ensure the preservation of two cardinals (typically $\aleph_1$ and $\aleph_2$), while returning to the much simpler situation where the side conditions form an $\in$-chain of elementary submodels. Applications of this method include very elegant proofs of new and old consistency results (see \cite{Neeman2017Twoapplicationsoffinitesideconditionsatomega2}, \cite{Velickovic:Venturi2011properforcingremastered} and \cite{Mohammadpour2023SpecialisingtreeswithsmallapproximationsI}), but most importantly, Neeman found a new and revolutionary way to iterate proper forcings with finite support (although not in the classical way), which he used to obtain an alternative proof of the consistency of the Proper Forcing Axiom ($\mathrm{PFA}$).  
	
	Some of the technical obstacles that prevented set theorists from tackling the problem of finding consistent high forcing axioms disappear when countable support iterations are replaced by finite support ones. This prompted set theorists to generalize Neeman's new iteration theory for proper forcings to other classes of posets, while developing new systems of side conditions in the process. Gitik and Magidor \cite{GitikMagidor2016SPFAbyfiniteconditions} and Veli\v ckovi\'c \cite{Velikovic2014iterationofsemiproper}, independently, found alternative finite support proofs of the consistency of the Semiproper Forcing Axiom ($\mathrm{SPFA}$), and Neeman himself announced in \cite{NeemanSlides1} and \cite{NeemanSlides2} that his iteration technique could be used to force a high analog of $\mathrm{PFA}$.
	
	Independently of Neeman, but around the same time, Asperó and Mota in \cite{Aspero:Mota2015a} and \cite{Aspero:Mota2015b} developed a new method for building finite support forcing iterations with finite symmetric systems of countable elementary submodels as side conditions. These iterations were used to prove the consistency of certain fragments of $\mathrm{PFA}$\footnote{These include the restriction of $\mathrm{PFA}$ to the classes of finitely proper forcings and forcings with the $\aleph_{1.5}$-chain condition.} together with arbitrarily large values of the continuum, thus showing that many known consequences of $\mathrm{PFA}$, in particular very strong failures of Shelah's Club Guessing principle, are independent of the size of the continuum. 
	
	The Asperó-Mota side conditions are no longer finite $\in$-chains of countable elementary submodels. Rather, they are finite matrices of countable elementary submodels, exhibiting some form of symmetry\footnote{Typically, models of the same rank are isomorphic and the system is closed under these isomorphisms.}. These side conditions had already been hinted in the early work of Todor\v cevi\'c on forcing with side conditions \cite{Todorcevic1984AnoteonthePFA}, and explicitly used in \cite{Todorcevic1985directedsets} to show that, consistently, there are only five cofinal types of directed sets of size $\aleph_1$. While Todor\v cevi\'c's original poset consisting of $\in$-chains of countable elementary submodels collapses cardinals, the matrix version ensures that the poset preserves all cardinals. In fact, Todor\v cevi\'c's matrices enjoy a strong form of the $\aleph_2$-chain condition, known as the $\aleph_2$-properness isomorphism condition, which is preserved under countable support iterations of length $<\omega_2$ and implies $\mathrm{CH}$ (see \cite{Abraham2010Properforcing}, \cite{Burke1998Forcingaxioms}, \cite{Shelah2017Properandimproperforcing} and \cite{Todorcevic1985directedsets} for more on this property). These properties have been a key ingredient in the proof of many consequences of $\mathrm{PFA}$ without the need to require large cardinals.
	
	These side conditions have proven to be as interesting as their classical counterparts and have found many applications in very different contexts. Besides the ones that have already been listed above, let us mention that, using the poset of matrices of countable elementary submodels, Kuzeljevi\'c and Todor\v cevi\'c in \cite{KuzeljevicTodorcevic2017Forcingwithmatrices} forced a Kurepa tree, an almost Souslin Kurepa tree and a diamond sequence; Abraham, Cummings and Smyth in \cite{AbrahamCummings2012moreresultsinpolychromatic} and \cite{AbrahamCummingsSmyth2007SomeresultsinpolychromaticsRamseytheory} forced certain colorings of $[\omega_2]^2$, which give negative polychromatic partition relations; Zapletal in \cite{Zapletal1997Stronglyalmostdisjointfunctions} forced strong almost disjoint families of functions from $\omega_1$ to $\omega$ of arbitrary size; Miyamoto in a series of papers (\cite{Miyamoto2012AnoteonageneralizedMAfortheNorwichposets}, \cite{Miyamoto2019AnoteonPffcandWTCGS12}, \cite{Miyamoto2023Astronglysigmaclosedposetthatforcesasimplifiedmorass}, \cite{Miyamoto2014Matricesofisomorphicmodelsandmorasslikestructures}, \cite{Miyamoto2015Squaresbymatriceswithcoherentsequences}, \cite{Miyamoto2022Forcingaclubbyageneralizedfastfunction}), some of them unpublished, forced different square principles on $\omega_2$, multiple variants of $(\omega_1,1)$-morasses, and $\omega_2$-Souslin trees. 
	
	In this paper, we will present a higher analog of the matrix-like side conditions used by Todor\v cevi\'c and Asperó-Mota, by integrating the technology of Neeman's two-type side conditions forcing. Our \emph{symmetric systems of models of two types} will have strong symmetry requirements to ensure that the preservation properties of the different forcing notions that we will introduce, strengthen those of Neeman's poset. In particular, under suitable extra assumptions, all the forcings that we will define in this paper will be strongly proper and strongly $\aleph_1$-proper (at least in a stationary set, which is enough to ensure the preservation of $\aleph_1$ and $\aleph_2$), will have the $\aleph_3$-Knaster condition, and will preserve $2^{\aleph_1}=\aleph_2$. 
	
	The paper is organized as follows. Section \ref{section-prelim} establishes the notation and provides all the necessary background in strongly proper forcing and elementary submodels. Section \ref{section-pure} introduces symmetric systems of models of two types, proves standard amalgamation lemmas needed for the preservation theorems, and ends with two applications: the pure side condition forcing adds an $\omega_1$-club subset of $\omega_2$ and a Kurepa tree on $\omega_2$. Section \ref{section-decorations} introduces the decorated version of the forcing with two-type symmetric systems, proves analogous amalgamation lemmas, which ensure the preservation of all cardinals and $2^{\aleph_1}=\aleph_2$, and ends with two applications: the forcing adds a club subset of $\omega_2$ which does not contain any infinite set from the ground model and a function on $\omega_2$ that bounds every canonical function below $\omega_3$ on a club. Section \ref{section-morass} introduces a variant of the decorated forcing from Section \ref{section-decorations}, which has the same preservation properties and adds a simplified $(\omega_2,1)$-morass.
	
	

	\section{Preliminaries}\label{section-prelim}
	
	Our notation will be standard and follow \cite{Kunen2011Settheory} and \cite{Jech2003Settheory}. We refer the reader to these two sources for any undefined notions. Unless otherwise specified, lowercase Greek letters $\alpha,\beta,\gamma,\delta,\varepsilon,\xi,\eta$ will be used to denote ordinals, while $\kappa,\lambda,\mu,\nu,\theta$ will be used to denote infinite cardinals. We will denote by $\mathrm{OR}$ the class of all ordinals. Let $X$ be any set. If $X$ is well-ordered, we will denote its order-type by $\ot(X)$. We will denote by $\mathcal{P}(X)$ the power set of $X$. If $\mu$ is a cardinal, we will denote by $[X]^\mu$ the set of all subsets of $X$ of size $\mu$. The sets $[X]^{<\mu}$ and $[X]^{\leq\mu}$ are defined in the obvious way. If $f$ is a function and $X$ is a set, the point-wise image of $X$ by $f$, denoted $f"(X)$, is the set $\{y:\exists x\in X(y=f(x))\}$. If $\mu<\lambda$ are two infinite regular cardinals, we will denote the set $\{\alpha<\lambda:\operatorname{cf}(\alpha)=\mu\}$ by $S_\mu^\lambda$. If $\mathbb{P}$ is a forcing notion, $G$ is a $\mathbb{P}$-generic filter over $V$ and $\tau$ is a $\mathbb{P}$-name, we will denote the interpretation of $\tau$ by $\tau^G$. If $x$ is an element of $V$, we will denote by $\check{x}$ its canonical $\mathbb{P}$-name.

	For the rest of the paper, fix a regular cardinal $\kappa>\omega_2$ and a relational structure $(H(\kappa);\in,\lhd)$, where $\lhd$ is a fixed well-order of $H(\kappa)$. We will usually refer to the structure $(H(\kappa);\in,\lhd)$ simply by $H(\kappa)$.
	
	The following are some basic facts about elementary submodels, which can be found in any standard reference for set-theoretic notions such as \cite{Kunen2011Settheory} and \cite{Jech2003Settheory}. We include them without proof, and we will use them throughout the paper, sometimes without mention. Additional information on elementary submodels and their applications in set theory can be found in \cite{Barwise1975Admissiblesetsandstructures}, \cite{Dow1998Anintroductiontoapplicationsofelementarysubmodels}, \cite{JustWeese1997DiscoveringmodernsetthoryII}, and \cite{Kunen2009Thefoudnationsofmathematics}.
	
	\begin{theorem}[Tarski-Vaught test]
		Let $M$ be a model and let $A\subseteq M$. Then, $A$ is the domain of an elementary submodel $N\preceq M$ iff for every formula $\varphi(y,\overline{x})$ and every tuple $\overline{a}$ of $A$ such that $M\models\exists y\varphi(y,\overline{a})$, there is $b\in A$ such that $M\models\varphi(b,\overline{a})$.
	\end{theorem}
	
	\begin{proposition}\label{prelim-elem-sub}
		If $Q_0,Q_1\preceq H(\kappa)$ are such that $Q_0\subseteq Q_1$, then $Q_0\preceq Q_1$.
	\end{proposition}
	
	\begin{proposition}\label{prelim-trans}
		Let $Q\preceq H(\kappa)$, and let $\mu<\kappa$ be a cardinal such that $\mu\subseteq Q$. For every $A\in Q$, if $H(\kappa)\models|A|=\mu$, then $A\subseteq Q$.
	\end{proposition}
	
	\begin{proposition}
		Let $Q\preceq H(\kappa)$. If $A\in H(\kappa)$ is definable over $H(\kappa)$ with parameters in $Q$, then $A\in Q$.
	\end{proposition}
	
	\begin{proposition}
		Let $\theta$ be a large enough regular cardinal such that $\theta>\kappa$. If $Q^*$ is an elementary submodel of $H(\theta)$ such that $H(\kappa),\lhd\in Q^*$, then $Q^*\cap H(\kappa)\preceq (H(\kappa);\in,\lhd)$.
	\end{proposition}
	
	\begin{proposition}\label{prelim-ordinal-height}
		Let $Q$ be an elementary submodel of $H(\kappa)$ such that $|Q|=\mu<\mu^+<\kappa$ and $\mu\subseteq Q$. Then $Q\cap\mu^+\in\mu^+$ is a limit ordinal.
	\end{proposition}
	
	\begin{proposition}
		Let $Q_0,Q_1\preceq H(\kappa)$ such that $|Q_0|=|Q_1|=\mu<\mu^+<\kappa$, $\mu\in Q_0\cap Q_1$ and $\mu\subseteq Q_0\cap Q_1$. If $\Psi$ is an isomorphism between the structures $(Q_0;\in,\lhd)$ and $(Q_1;\in,\lhd)$, then $\Psi$ is the identity on $Q_0\cap\mu^+$. In particular, $Q_0\cap\mu^+=Q_1\cap\mu^+$.
	\end{proposition}
	
	Let $\mathcal{S}$ denote the set of countable elementary submodels of $(H(\kappa);\in,\lhd)$. It is well-known that $\mathcal{S}\subseteq[H(\kappa)]^\omega$ is a club in the sense of Jech \cite{Jech2003Settheory}. We will tend to use capital letters $M,N$ to refer to the elements of $\mathcal{S}$, which we will call \emph{countable} or \emph{small models}. 
	
	Let $\mathcal{L}$ denote an arbitrary set of $\aleph_1$-sized elementary substructures $X$ of $(H(\kappa);\in,\lhd)$ with $\omega_1\subseteq X$. We will use the capital letters $X,Y,Z$ to refer to the elements of $\mathcal{L}$, which we will call \emph{uncountable} or \emph{large models}. 
	
	To refer to models of arbitrary size, we will use $Q$, as well as other capital letters further down the alphabet. If $Q$ is an elementary submodel of $H(\kappa)$, we will usually refer to the structure $(Q;\in,\lhd\cap (Q\times Q))$ by $(Q;\in,\lhd)$. Moreover, we will indistinctly use $Q$ to refer to the structure $(Q;\in,\lhd)$ and its universe. Throughout the paper, it will be clear from the context to which one we are referring to.
	
	\begin{definition}\label{def-unbd-stat}
		Let $\mathcal{T}$ be a subset of $\mathcal{P}(H(\kappa))$.
		\begin{itemize}
			\item $\mathcal{T}$ is \emph{unbounded in $H(\kappa)$} if for any $x\in H(\kappa)$, there is $T\in\mathcal{T}$ such that $x\in T$ (\cite{Neeman2014Forcingwithsequencesofmodelsoftwotypes}).
			
			\item $\mathcal{T}$ is \emph{stationary in $H(\kappa)$} if for any $f:[H(\kappa)]^{<\omega}\to H(\kappa)$, there is $T\in\mathcal{T}$ such that $f"([T]^{<\omega})\subseteq T$.
		\end{itemize}
	\end{definition}
	
	The following lemma is a folklore result that can be found in many of the standard set-theoretic references listed at the beginning of this section. 
		
	\begin{lemma}\label{prelim-folklore}
		Let $\theta$ be a regular cardinal. Let $Q$ be an elementary submodel of $H(\theta)$.	Let $\mathbb{P}$ be a forcing poset such that $\mathbb{P}\in Q$, and let $G$ be a $\mathbb{P}$-generic filter over $V$. Then, in $V[G]$,
		\[
		Q[G]\preceq H(\theta)[G]=H(\theta)^{V[G]}.
		\]
	\end{lemma}


	\subsection{Strongly proper forcing and cardinal preservation}
			
	This subsection gives the basic necessary background on strongly proper forcing. This notion, which strengthens the usual notion of properness, was isolated by Mitchell in \cite{Mitchell2005Addingclubsofomega2withfiniteforcing}. 
	
	\begin{definition}
		Let $\mathbb{P}$ be a poset. For any set $M$ and any $p\in\mathbb{P}$, we say that $p$ is \emph{strongly $(M,\mathbb{P})$-generic} if any of the following three equivalent conditions are satisfied:
		\begin{itemize}
			\item For any $D\subseteq\mathbb{P}\cap M$ dense in $\mathbb{P}\cap M$, $D$ is predense below $p$ in $\mathbb{P}$.
			
			\item For any $P\subseteq\mathbb{P}\cap M$ predense in $\mathbb{P}\cap M$, $P$ is predense below $p$ in $\mathbb{P}$.
			
			\item For any $A\subseteq\mathbb{P}\cap M$ maximal antichain in $\mathbb{P}\cap M,$ $A$ is predense below $p$ in $\mathbb{P}$.
		\end{itemize}
	\end{definition}
	
	\begin{lemma}\label{strong-prop-lemma}
		Let $\mathbb{P}$ be a poset. For any set $M$ and any $p\in\mathbb{P}$, the following two are equivalent:
		\begin{itemize}
			\item $p$ is strongly $(M,\mathbb{P})$-generic.
			
			\item For every $q\leq p$ in $\mathbb{P}$, there is some $q_M\in\mathbb{P}\cap M$ such that every condition in $\mathbb{P}\cap M$ extending $q_M$ is compatible with $q$ in $\mathbb{P}$.
		\end{itemize}
		
		Such a condition $q_M$ is typically called a \emph{residue of $q$ in $M$}.
	\end{lemma}
	
	\begin{definition}
		If $\mathcal{K}\subseteq\mathcal{S}\cup\mathcal{L}$, a poset $\mathbb{P}\subseteq H(\kappa)$ is called \emph{strongly $\mathcal{K}$-proper}, if for any $Q\in\mathcal{K}$ and any $p\in\mathbb{P}\cap Q$, there is a strongly $(Q,\mathbb{P})$-generic condition $q\in\mathbb{P}$ such that $q\leq p$. 
	\end{definition}
	
	\begin{lemma}\label{preservacio-proper}
		If $\mathcal{L}$ is a stationary set of large models in $H(\kappa)$ and $\mathbb{P}\subseteq H(\kappa)$ is a strongly $\mathcal{S}$-proper and strongly $\mathcal{L}$-proper poset, then forcing with $\mathbb{P}$ does not collapse $\aleph_1$ nor $\aleph_2$.
	\end{lemma}
	\begin{proof}
		We only prove the preservation of $\aleph_1$, since the preservation of $\aleph_2$ is proved in the exact same way.
		
		For the preservation of $\aleph_1$, let $\dot{f}$ be a $\mathbb{P}$-name for a function from $\omega^V$ to $\omega_1^V$. We prove that every $p\in \mathbb{P}$ can be extended to a condition forcing that the function $\dot{f}$ is bounded below $\omega_1^V$. Fix some condition $p\in\mathbb{P}$. Find, for a large enough regular cardinal $\theta>\kappa$, an elementary submodel $M^*\preceq H(\theta)$ such that $H(\kappa),\mathbb{P},p,\dot{f}\in M^*$ and $M:=M^*\cap H(\kappa)\in\mathcal{S}$. Since $\mathbb{P}\subseteq H(\kappa)$, the condition $p$ is a member of $M$, and thus, by strong $\mathcal{S}$-properness, there is a strongly $(M,\mathbb{P})$-generic condition $q\leq p$. We claim that $q$ forces $\mathrm{ran}(\dot{f})\subseteq M$. Since $\omega^V<M\cap\omega_1^V<\omega_1^V$ by Proposition \ref{prelim-ordinal-height}, this implies that $q$ forces $\dot{f}$ to be bounded below $\omega_1^V$.
		
		Suppose on the contrary that there is some $q'\leq q$ forcing that $\dot{f}(n)=\check{\alpha}$, for some $n<\omega$ and some $\alpha<\omega_1^V$ such that $\alpha\geq M\cap\omega_1^V$. Since $q'$ is also strongly $(M,\mathbb{P})$-generic, there is a residue $q'_M$ of $q'$ in $M$. By elementarity, there is a condition $r\in \mathbb{P}\cap M^*$ and an ordinal $\beta\in M^*\cap\omega_1^V$ such that $r\leq q'_M$ and $r\Vdash\dot{f}(n)=\check{\beta}$. Since $\mathbb{P}\cap M^*=\mathbb{P}\cap M$ and $M^*\cap\omega_1^V=M\cap\omega_1^V$, we have $r\in\mathbb{P}\cap M$ and $\beta\in M\cap\omega_1^V$, and we reach a contradiction because $r$ and $q'$ are compatible by Lemma \ref{strong-prop-lemma}, and $\beta<M\cap\omega_1^V\leq\alpha$.
	\end{proof}

	We finish this subsection with a well-known preservation theorem, which combined with the preservation theorems for strongly $\mathcal{S}$-proper and strongly $\mathcal{L}$-proper forcings, will ensure that all the forcing notions presented in this paper preserve all cardinals.
		
	\begin{definition}
		Let $\lambda$ be an uncountable cardinal. A poset $\mathbb{P}$ has the \emph{$\lambda$-Knaster condition} if for every $A\subseteq\mathbb{P}$ of size $\lambda$ there is a subset $B\subseteq A$ of the same size consisting of pairwise compatible conditions.
	\end{definition}
	
	\begin{remark}
		Every poset $\mathbb{P}$ with the $\lambda$-Knaster condition has the $\lambda$-chain condition (denoted $\lambda$-c.c.).
	\end{remark}
	
	\begin{lemma}\label{preservation-cc}
		If $\lambda$ is an uncountable cardinal and $\mathbb{P}$ is a forcing notion with the $\lambda$-c.c., then every cardinal $\mu\geq\lambda$ is preserved after forcing with $\mathbb{P}$.
	\end{lemma}

	

	\subsection{Appropriate sets of elementary submodels}
	
	Given a model $Q\in\mathcal{S}\cup\mathcal{L}$, we will denote $Q\cap\omega_1$ by $\delta_Q$, and $\sup(Q\cap\omega_2)$ by $\varepsilon_Q$. We will call $\delta_Q$ the \textit{$\omega_1$-height of $Q$}, and $\varepsilon_Q$ the \textit{$\omega_2$-height of $Q$}. 
	
		
	\begin{definition}
		Given a model $Q$, we let
		\[
		Q[\omega_1]:=\{f(\alpha):f\in Q,f\text{ a function with 	}\operatorname{dom}(f)=\omega_1, \alpha\in\omega_1\},
		\]
		and we call it the \emph{$\omega_1$-hull of $Q$}, or simply the \emph{hull of $Q$}.
	\end{definition}
	
	\begin{remark}
		If $Q\in\mathcal{S}\cup\mathcal{L}$ and $M\in\mathcal{S}\cap Q$, then $M[\omega_1]$ is definable in $Q$, and thus $M[\omega_1]\in Q$. Consequently, if $M\in\mathcal{S}$, then $M[\omega_1]$ cannot be countably-closed, otherwise $M$ would be a member of $M[\omega_1]$, and hence $M[\omega_1]\in M[\omega_1]$, which is impossible by the axiom of foundation.
	\end{remark}
	

	\begin{proposition}\label{prop7}
		Let $Q\in\mathcal{S}\cup\mathcal{L}$. Then, $Q[\omega_1]$ is the smallest elementary submodel of $H(\kappa)$ that contains $Q\cup\omega_1$ as a subset. 
	\end{proposition}
	\begin{proof}
		First note that since $\omega_1\in Q$, the identity function $id:\omega_1\to\omega_1$ is definable in $Q$, and thus, $\omega_1\subseteq Q[\omega_1]$. Moreover, for each $a\in Q$, the constant function sending all $\alpha\in\omega_1$ to $a$ is definable in $Q$, and so $Q\subseteq Q[\omega_1]$. Therefore, $Q\cup\omega_1\subseteq Q[\omega_1]$.
		
		We use the Tarski-Vaught test to show that $Q[\omega_1]$ is an elementary submodel of $H(\kappa)$. Let $\varphi(y,x_0,\dots,x_n)$ be a first-order formula in the language of set theory. Let $a_0,\dots,a_n\in Q[\omega_1]$ be such that $H(\kappa)\models\exists y\varphi(y,a_0,\dots,a_n)$. We have to find $b\in Q[\omega_1]$ such that $H(\kappa)\models\varphi(b,a_0,\dots,a_n)$. Let $f_i\in Q$ and $\alpha_i\in\omega_1$ be such that $a_i=f_i(\alpha_i)$, for all $i\leq n$. Fix a bijection $F:\omega_1^{<\omega}\to\omega_1$. There exists a function $g$ with $\dom(g)=\omega_1$ in $H(\kappa)$ such that for any $y\in H(\kappa)$ and $\beta_0,\dots,\beta_n<\omega_1$,
		\[
		H(\kappa)\models \varphi\big(y,f_0(\beta_0),\dots,f_n(\beta_n)\big)\to \varphi\Big(g\big(F(\beta_0,\dots,\beta_n)\big),f_0(\beta_0),\dots,f_n(\beta_n)\Big).
		\]
		Note that $g\circ F$ is defined with $f_0,\dots,f_n$ and $\omega_1$ as parameters. Hence, we can assume that $g\circ F\in Q$. Now, let $b=g\big(F(\alpha_0,\dots,\alpha_n)\big)$, and note that $b\in Q[\omega_1]$ and $H(\kappa)\models \varphi(b,a_0,\dots,a_n)$. Therefore, $Q[\omega_1]\preceq H(\kappa)$.
		
		The minimality of $Q[\omega_1]$ follows from the fact that if $P$ is an elementary submodel of $H(\kappa)$ such that $Q\cup\omega_1\subseteq P$, then $Q[\omega_1]\subseteq P$.
	\end{proof}
	
	\begin{corollary}
		If $X$ is a large model, then $X[\omega_1]=X$.
	\end{corollary}
	
		\begin{proposition}\label{heighthull}
		Let $M\in\mathcal{S}$, and let $\alpha\in M$ be an ordinal such that $cf(\alpha)>\omega_1$. Then, $\sup(M[\omega_1]\cap\alpha)=\sup(M\cap\alpha)$. 
	\end{proposition}
	\begin{proof}
		The inequality $\geq$ is clear. For the other direction, let $\xi\in M[\omega_1]\cap\alpha$. By definition, there are a function $f\in M$ and an ordinal $\beta\in\omega_1$ such that $f(\beta)=\xi$. Define a function $g$ on $\omega_1$ by letting $g(\gamma)=f(\gamma)$ if $f(\gamma)\in\alpha$ and $g(\gamma)=0$ otherwise, for each $\gamma\in\omega_1$. Note that since $f,\omega_1,\alpha\in M$, the function $g$ is definable in $M$. Therefore, as $\operatorname{dom}(g)=\omega_1$ and $cf(\alpha)>\omega_1$, we have $\sup(\operatorname{ran}(g))\in M\cap\alpha$. Hence
		\[
		\xi=f(\beta)=g(\beta)<\sup(M\cap\alpha).
		\] 
	\end{proof}
	
	\begin{corollary}\label{prelim-height-hull}
		If $Q\in\mathcal{S}\cup\mathcal{L}$, then $\varepsilon_Q=\varepsilon_{Q[\omega_1]}=Q[\omega_1]\cap\omega_2$. In particular, if $Q\in\mathcal{S}$, then $\varepsilon_{Q[\omega_1]}$ is an ordinal of countable cofinality. 
	\end{corollary}

	\begin{definition}
		Let $Q_0,Q_1\in\mathcal{S}\cup\mathcal{L}$. A map $\Psi$ between $Q_0[\omega_1]$ and $Q_1[\omega_1]$ is called an \emph{$\omega_1$-isomorphism} if $\Psi$ is the unique isomorphism between the structures $(Q_0[\omega_1];\in,Q_0,\lhd)$ and $(Q_1[\omega_1];\in,Q_1,\lhd)$, and $\Psi$ is the identity on $Q_0[\omega_1]\cap Q_1[\omega_1]$.
		
		If $\Psi$ is an $\omega_1$-isomorphim between $Q_0[\omega_1]$ and $Q_1[\omega_1]$, we will say that $Q_0$ and $Q_1$ are \emph{$\omega_1$-isomorphic} and we will denote this by $Q_0[\omega_1]\cong Q_1[\omega_1]$. Moreover, we will denote the (unique) $\omega_1$-isomorphism $\Psi$ by $\Psi_{Q_0[\omega_1],Q_1[\omega_1]}$.
		
		If $X_0,X_1\in\mathcal{L}$ and $X_0=X_0[\omega_1]\cong X_1[\omega_1]=X_1$, we will simply write $X_0\cong X_1$ and $\Psi_{X_0,X_1}$.
	\end{definition}

	\begin{proposition}\label{prop5}
		Let $P,Q_0,Q_1\in\mathcal{S}\cup\mathcal{L}$. Suppose that $P\in Q_0[\omega_1]$ and $Q_0[\omega_1]\cong Q_1[\omega_1]$, witnessed by $\Psi$. Then, the following hold:
		\begin{enumerate}
			\item $(P;\in,\lhd)\preceq(Q_0[\omega_1];\in,\lhd)$.
			
			\item $(\Psi"(P);\in,\lhd)\preceq(H(\kappa);\in,\lhd)$.
			
			\item $\Psi(P)=\Psi"(P)$.
		\end{enumerate} 
	\end{proposition}
	\begin{proof}
		For (1), note that since $(Q_0[\omega_1];\in,\lhd)\preceq(H(\kappa);\in,\lhd)$, $\omega_1\subseteq Q_0[\omega_1]$, and $H(\kappa)\models|P|\leq\omega_1$, we have $P\subseteq Q_0[\omega_1]$ by Proposition \ref{prelim-trans}. Hence, the conclusion of (1) follows from Proposition \ref{prelim-elem-sub}.
		
		For (2), note first that it is straightforward to see that the restriction map $\Psi\restr P$ is an isomorphism between the structures $(P;\in,\lhd)$ and $(\Psi"(P);\in,\lhd)$. Let $\varphi(y,x_1,\dots,x_n)$ be a first-order formula in the language of set theory and let $p_1,\dots,p_n\in P$. Suppose that $H(\kappa)\models\varphi\exists y\varphi(y,\Psi(p_1),\dots,\Psi(p_n))$. Since $Q_1[\omega_1]\preceq H(\kappa)$ and $\Psi"(P)\subseteq Q_1[\omega_1]$, the model $Q_1[\omega_1]$ also satisfies the formula $\exists y\varphi(y,\Psi(p_1),\dots,\Psi(p_n))$, and thus $Q_0[\omega_1]\models\exists y\varphi(y,p_1,\dots,p_n)$. Since $p_1,\dots,p_n\in P$, we have $P\models\exists y\varphi(y,p_1,\dots,p_n)$ by (1). Let $q\in P$ such that $P\models\varphi(q,p_1,\dots,p_n)$. By chasing back the diagram, we can show that $H(\kappa)\models\varphi(q,p_1,\dots,p_n)$, and by the Tarski-Vaught test, we can conclude that $(\Psi"(P);\in,\lhd)\preceq(H(\kappa);\in,\lhd)$.

		For (3), note that since $(Q_0[\omega_1];\in,\lhd)\preceq(H(\kappa);\in,\lhd)$, the model $Q_0[\omega_1]$ knows that $|P|\leq\omega_1$. Hence, $Q_1[\omega_1]$ knows that $|\Psi(P)|\leq\Psi(\omega_1)=\omega_1$. By assumption $(Q_1[\omega_1];\in,\lhd)\preceq(H(\kappa);\in,\lhd)$, hence $H(\kappa)$ knows that $|\Psi(P)|\leq\omega_1$. Now, apply Proposition \ref{prelim-trans} to conclude $\Psi(P)\subseteq Q_1[\omega_1]=\Psi"(Q_0[\omega_1])$. Therefore, $\Psi(P)\subseteq\Psi"(P)$. However, $\Psi"(P)\subseteq\Psi(P)$. Hence, we have $\Psi(P)=\Psi"(P)$ as we wanted. 
	\end{proof}
	
	
	
	\begin{proposition}\label{prop-iso}
		Let $M_0,M_1\in\mathcal{S}$ and suppose that $M_0[\omega_1]\cong M_1[\omega_1]$. Then $\Psi_{M_0[\omega_1],M_1[\omega_1]}\upharpoonright M_0$ is the unique isomorphism between the structures $(M_0;\in,\lhd)$ and $(M_1;\in,\lhd)$. 
	\end{proposition}
	
	\begin{proposition}
		Let $M_0,M_1\in\mathcal{S}$ be such that $M_0[\omega_1]\cong M_1[\omega_1]$. Then, $\varepsilon_{M_0}=\varepsilon_{M_1}$.
	\end{proposition}
	\begin{proof}
		Follows from Proposition \ref{prop-iso} and Corollary \ref{prelim-height-hull}.
	\end{proof}
	
	\begin{definition}\label{def-appropriate}
		A set $\mathcal{L}$ of large models is said to be \emph{appropriate for $\mathcal{S}$} if the following hold:
		\begin{itemize}
			\item If $X\in\mathcal{L}$ and $M\in\mathcal{S}$ are such that $X\in M$, then $X\cap M\in X\cap \mathcal{S}$.
			
			\item If $Q_0,Q_1\in\mathcal{S}\cup\mathcal{L}$ are $\omega_1$-isomorphic and $X\in Q_0[\omega_1]\cap\mathcal{L}$, then $\Psi_{Q_0[\omega_1],Q_1[\omega_1]}(X)\in\mathcal{L}$. 
		\end{itemize}
	\end{definition}
	
	As we shall see in subsequent sections, all the forcing notions introduced in this paper will be proved to be strongly $\mathcal{S}$-proper and strongly $\mathcal{L}$-proper. Hence, we will be interested in sets $\mathcal{L}$ of large models appropriate for $\mathcal{S}$ and stationary in $H(\kappa)$, so that, in light of Lemma \ref{preservacio-proper}, both cardinals $\aleph_1$ and $\aleph_2$ are preserved. However, since we want to make this work as general as possible, throughout the paper, the stationarity of $\mathcal{L}$ will not be assumed a priori, and we will specify exactly where it is needed.
	
	\begin{example}
		Assuming $\mathrm{CH}$, if $\mathcal{L}$ denotes the set of elementary submodels $X\preceq (H(\kappa);\in,\lhd)$ of size $\aleph_1$ and such that $^\omega X\subseteq X$, then $\mathcal{L}$ is appropriate for $\mathcal{S}$ and stationary in $H(\kappa)$.
	\end{example}
	
	\begin{example}
		Recall that a large model $X\in\mathcal{L}$ is said to be \emph{internally club} if $X=\bigcup_{i<\omega_1}M_i$, where $\la M_i:i<\omega_1\ra$ is an increasing and continuous sequence such that for all $i<\omega_1$, $M_i\in X\cap [H(\kappa)]^\omega$. Let $\mathcal{L}$ denote the set of internally club large models. We claim that $\mathcal{L}$ is appropriate for $\mathcal{S}$. Let $X\in\mathcal{L}$ and suppose that $\la M_i:i<\omega_1\ra$ is the $\lhd$-least sequence witnessing that $X$ is internally club. Suppose first that $N\in\mathcal{S}$ is such that $X\in N$. Then $\la M_i:i<\omega_1\ra\in N$ and $X\cap N=M_{\delta_N}\in X\cap\mathcal{S}$. Suppose now that $Q_0,Q_1\in\mathcal{S}\cup\mathcal{L}$ are $\omega_1$-isomorphic and $X\in Q_0[\omega_1]\cap\mathcal{L}$. Since $\la M_i:i<\omega_1\ra\in Q_0[\omega_1]$ and $M_i\in Q_0[\omega_1]$ for every $i<\omega_1$, the isomorphic copy $\Psi_{Q_0[\omega_1],Q_1[\omega_1]}(X)$ of $X$ is internally club, and this is witnessed by the sequence $\la\Psi_{Q_0[\omega_1],Q_1[\omega_1]}(M_i):i<\omega_1\ra$.
	\end{example}
	
	\begin{proposition}
		Let $\mathcal{L}$ be an appropriate set of large models for $\mathcal{S}$. Then for any $X\in\mathcal{L}$ we have the following:
		\begin{itemize}
			\item $X$ is internally club.
			
			\item $\cf(\varepsilon_X)=\omega_1$.
		\end{itemize}
	\end{proposition}
	\begin{proof}
		Let $\la M_i:i<\omega_1\ra$ be a continuously increasing $\in$-chain of small models such that $X\in M_i$ for every $i<\omega_1$. Since $\mathcal{L}$ is appropriate for $\mathcal{S}$, we have $X\cap M_i\in X$ for all $i<\omega_1$. We claim that $X=\bigcup\{X\cap S_i:i<\omega_1\}$. Let $e:\omega_1\to X$ be a surjection. Since $X\in M_0$ and $M_0$ is an elementary submodel of $H(\kappa)$, we may assume that $e\in M_0$. Since $M_i\in M_j$ entails $M_i\subseteq M_j$, we have $e\in M_i$ for all $i<\omega_1$. Now by induction, we have $i\subseteq S_i$ for any $i<\omega_1$. Hence, $e"(i)\subseteq X\cap S_i$, and we can conclude that 
		\[
		X=e"(\omega_1)=\bigcup\{X\cap S_i:i<\omega_1\}.
		\]
	\end{proof}
	
	For the reminder of the paper, we fix a set of large models $\mathcal{L}$ appropriate for $\mathcal{S}$.

	The next result is very fundamental in the understanding of the structure of side conditions of models of two types and their limitations. It was already noticed by Neeman \cite{Neeman2014Forcingwithsequencesofmodelsoftwotypes}, and as we will see in Proposition \ref{prop-res-gaps}, the structure of our two-type symmetric systems crucially depends on the following fact.
	
	\begin{proposition}\label{propgaps}
		Let $M\in\mathcal{S}$ and $X\in\mathcal{L}\cap M$. Then, for every $Q\in\mathcal{S}\cup\mathcal{L}$ such that $\varepsilon_{X\cap M}\leq\varepsilon_Q<\varepsilon_X$, $Q\notin M$.
	\end{proposition}
	\begin{proof}
		Suppose, towards a contradiction, that there is some $Q\in\mathcal{S}\cup\mathcal{L}$ such that $\varepsilon_{X\cap M}\leq\varepsilon_Q<\varepsilon_X$ and $Q\in M$. First, note that since $X\in\mathcal{L}$, then $\varepsilon_X=X\cap\omega_2$ is an ordinal in $\omega_2$. Therefore, $\varepsilon_Q=\sup(Q\cap\omega_2)$ is an ordinal in $X\cap\omega_2$, and hence $\varepsilon_Q\in M\cap (X\cap\omega_2)$. So we can conclude that $\varepsilon_Q<\varepsilon_{X\cap M}$, which contradicts the assumption $\varepsilon_{X\cap M}\leq\varepsilon_Q$.
	\end{proof}

	We will next prove some basic results on the interplay between the operations of taking intersections and taking $\omega_1$-isomorphic copies of small and large elementary submodels, in the way described in Definition \ref{def-appropriate}. In particular, we will see that, in our setting and in some cases, these two operations commute. 
	
	\begin{proposition}\label{prop20}
		Let $Q_0, Q_1,Q_1'\in\mathcal{S}\cup\mathcal{L}$ be such that $Q_0\in Q_1[\omega_1]$ and $Q_1[\omega_1]\cong Q_1'[\omega_1]$. Let $Q_0'$ denote $\Psi_{Q_1[\omega_1],Q_1'[\omega_1]}(Q_0)$. Then $$Q_0'[\omega_1]=\Psi_{Q_1[\omega_1],Q_1'[\omega_1]}(Q_0[\omega_1]),$$ and $\Psi_{Q_1[\omega_1],Q_1'[\omega_1]}\upharpoonright Q_0[\omega_1]$ witnesses $Q_0[\omega_1]\cong Q_0'[\omega_1]$.
	\end{proposition}
	\begin{proof}
		Denote the isomorphism $\Psi_{Q_1[\omega_1],Q_1'[\omega_1]}$ by $\Psi$. For the first part, it suffices to note that for any $f(\gamma)\in Q_0[\omega_1]$ with $f\in Q_0$ and $\gamma<\omega_1$, we have
		\[
		\Psi_{Q_1[\omega_1],Q_1'[\omega_1]}(f(\gamma))=\Psi_{Q_1[\omega_1],Q_1'[\omega_1]}(f)(\gamma),
		\] 
		which is an element of $Q_0'[\omega_1]$, because $\Psi_{Q_1[\omega_1],Q_1'[\omega_1]}(f)\in Q_0'$.
		
		
		For the second part, note that it follows easily that $\Psi\upharpoonright Q_0[\omega_1]$ is an isomorphism between the structures $(Q_0[\omega_1];\in,Q_0,\lhd)$ and $(Q_0'[\omega_1];\in,Q_0',\lhd)$. To see that it is the identity on $Q_0[\omega_1]\cap Q_0'[\omega_1]$, we only need to note that $Q_0[\omega_1]\cap Q_0'[\omega_1]$ is contained in $Q_1[\omega_1]\cap Q_1'[\omega_1]$.
	\end{proof}
	
	\begin{corollary}\label{corollary8}
		Let $M_0,M_1\in\mathcal{S}$ be such that $M_0[\omega_1]\cong M_1[\omega_1]$. Suppose that $X_0\in\mathcal{L}\cap M_0$ and denote $\Psi_{M_0[\omega_1],M_1[\omega_1]}(X_0)$ by $X_1$. Then, the map $\Psi_{M_0[\omega_1],M_1[\omega_1]}\upharpoonright(X_0\cap M_0)[\omega_1]$ witnesses $(X_0\cap M_0)[\omega_1]\cong(X_1\cap M_1)[\omega_1]$.
	\end{corollary}
	
	\begin{proposition}\label{prop9}
		Let $M\in\mathcal{S}$ and $X_0,X_1\in\mathcal{L}\cap M$. Suppose that $X_0\cong X_1$. Then, $\Psi_{X_0,X_1}\upharpoonright(X_0\cap M)[\omega_1]$ witnesses $(X_0\cap M)[\omega_1]\cong(X_1\cap M)[\omega_1]$.
	\end{proposition}
	
	\begin{proof}
		Note that, in light of Proposition \ref{prop20}, it suffices to check that $\Psi_{X_0,X_1}"(X_0\cap M)$ equals $X_1\cap M$. But this easily follows from the fact that $\Psi_{X_0,X_1}$ is definable in $M$ with parameters $X_0$ and $X_1$, and therefore $\Psi_{X_0,X_1}\in M$. The image $\Psi_{X_0,X_1}(X_0\cap M)\in X_1$ of $X_0\cap M$ coincides with $\Psi_{X_0,X_1}"(X_0\cap M)$, as $\Psi_{X_0,X_1}(X_0\cap M)\subseteq X_1$.
	\end{proof}	
	
	\begin{proposition}\label{prop10}
		Let $M_0,M_1\in\mathcal{S}$ be such that $M_0[\omega_1]\cong M_1[\omega_1]$. Let $X_0,X_1\in\mathcal{L}$ be such that $X_0\cong X_1$, $X_0\in M_0$, and $X_1\in M_1$. Then, the map $\Psi_{X_0,X_1}\upharpoonright(X_0\cap M_0)[\omega_1]$ witnesses $(X_0\cap M_0)[\omega_1]\cong(X_1\cap M_1)[\omega_1]$.
	\end{proposition}
	\begin{proof}
		Let $X_0'$ denote $\Psi_{M_0[\omega_1],M_1[\omega_1]}(X_0)$. Note that $X_0'\in M_1$ and that $\Psi_{X_0,X_0'}:=\Psi_{M_0[\omega_1],M_1[\omega_1]}\upharpoonright X_0$ is the unique isomorphism between the models $(X_0;\in,\lhd)$ and $(X_0';\in,\lhd)$. By Proposition \ref{prop20}, $X_0$ and $X_0'$ are isomorphic witnessed by $\Psi_{M_0[\omega_1],M_1[\omega_1]}\restr X_0$. Since $X_0$ and $X_1$ are isomorphic by assumption, $X_0'$ and $X_1$ are isomorphic witnessed by $\Psi_{X_0',X_1}=\Psi_{X_0,X_1}\circ\Psi_{X_0',X_0}$. Thus we have $\Psi_{X_0,X_1}=\Psi_{X_0',X_1}\circ\Psi_{X_0,X_0'}$. It follows from Corollary \ref{corollary8} that the map $\Psi_{X_0,X_0'}\upharpoonright\big((X_0\cap M_0)[\omega_1]\big)$, which equals $\Psi_{M_0[\omega_1],M_1[\omega_1]}\upharpoonright\big((X_0\cap M_0)[\omega_1]\big)$, is an isomorphism between $$\big((X_0\cap M_0)[\omega_1];\in,X_0\cap M_0,\lhd\big)$$ and $$\big((X_0'\cap M_1)[\omega_1];\in,X_0'\cap M_1,\lhd\big).$$ Hence, we only need to show that $\Psi_{X_0',X_1}\upharpoonright\big((X_0'\cap M_1)[\omega_1]\big)$ is the unique isomorphism between $(X_0'\cap M_1)[\omega_1]$ and $(X_1\cap M_1)[\omega_1]$, but this follows from the last proposition.
	\end{proof}	
	
	\begin{proposition}\label{prop21}
		Let $M_0,M_1\in\mathcal{S}$ such that $M_0[\omega_1]\cong M_1[\omega_1]$. If $X\in\mathcal{L}$ is a member of $M_0\cap M_1$, then $X\cap M_0=X\cap M_1$.
	\end{proposition}
	\begin{proof}
		Since $X\in M_0[\omega_1]\cap M_1[\omega_1]$, the model $X$ is fixed by $\Psi_{M_0[\omega_1],M_1[\omega_1]}$. Therefore, as $X\cap M_0\in X$, and thus $X\cap M_0\in M_0[\omega_1]\cap M_1[\omega_1]$, we have $X\cap M_0=\Psi_{M_0[\omega_1],M_1[\omega_1]}(X\cap M_0)=X\cap M_1$.
	\end{proof}	
	
	\begin{lemma}\label{agreement}
		For any two models $Q_0,Q_1\in\mathcal{S}\cup\mathcal{L}$ such that $Q_0[\omega_1]\cong Q_1[\omega_1]$, the following hold:
		\begin{enumerate}
			\item $Q_0\cap Q_1\cap\omega_3$ is an initial segment of both $Q_0\cap\omega_3$ and $Q_1\cap\omega_3$.
			
				
			
			\item $Q_0[\omega_1]\cap Q_1[\omega_1]\cap\omega_3$ is an initial segment of both $Q_0[\omega_1]\cap\omega_3$ and $Q_1[\omega_1]\cap\omega_3$.
			
				
		\end{enumerate}
	\end{lemma}
	\begin{proof}
		Note that (2) entails (1) via $Q_0[\omega_1]\cong Q_1[\omega_1]$. Hence, it suffices to show that for every $\beta\in Q_0[\omega_1]\cap Q_1[\omega_1]\cap\omega_3$, if $\alpha\in Q_0[\omega_1]\cap\beta$, then $\alpha\in Q_1[\omega_1]\cap\beta$. Fix a surjection $\pi_\beta:|\beta|\to\beta$ such that $\pi_\beta\in Q_0[\omega_1]\cap Q_1[\omega_1]$.
		
		Let $\xi\in Q_0[\omega_1]\cap\omega_2$ be such that $\pi_\beta(\xi)=\alpha$, and note that $\Psi_{Q_0[\omega_1],Q_1[\omega_1]}(\xi)$ is a member of $Q_1[\omega_1]\cap\omega_2$. However, since $\Psi_{Q_0[\omega_1],Q_1[\omega_1]}$ is the identity on $Q_0[\omega_1]\cap\omega_2$, we have $\Psi_{Q_0[\omega_1],Q_1[\omega_1]}(\xi)=\xi\in Q_1[\omega_1]\cap\omega_2$. Therefore, $\alpha=\pi_\beta(\xi)$ is a member of $Q_1[\omega_1]$ as we wanted.
	\end{proof}

	\section{Two-type symmetric systems}\label{section-pure}
	
	In this section, we will introduce symmetric systems of models of two types, by naturally combining the ideas of Neeman's chains of models of two types and Asperó-Mota-Todor\v cevi\'c's symmetric systems of elementary submodels, and we will prove their main properties. Not surprisingly, we will see that the forcing consisting of two-type symmetric systems is strongly $\mathcal{S}$-proper and strongly $\mathcal{L}$-proper, has the $\aleph_3$-Knaster condition if $2^{\aleph_1}=\aleph_2$, and preserves $2^{\aleph_1}=\aleph_2$.
	
	We will finish the section with two applications of this forcing notion. First, we will show how it adds an $\omega_1$-club subset of $\omega_2$. Then, building on ideas of Kuzeljevi\'c and Todor\v cevi\'c (\cite{KuzeljevicTodorcevic2017Forcingwithmatrices}), we will show that it also adds a Kurepa tree on $\omega_2$.
	
	For any set $\mathcal{M}\subseteq\mathcal{S}\cup\mathcal{L}$, we let $\dom(\mathcal{M})$ denote the set of $\omega_2$-heights $\{\varepsilon_Q:Q\in\mathcal{M}\}$. If $\varepsilon\in\dom(\mathcal{M})$, let $\mathcal{M}(\varepsilon)$ be the set of all models $Q\in\mathcal{M}$ such that $\varepsilon_Q=\varepsilon$. Moreover, we let $\mathcal{M}(<\varepsilon)$ denote the set of models $Q\in\mathcal{M}$ such that $\varepsilon_Q<\varepsilon$. We define similarly $\mathcal{M}(\leq\varepsilon)$, $\mathcal{M}(>\varepsilon)$, and $\mathcal{M}(\geq\varepsilon)$. If $\varepsilon_0<\varepsilon_1$ are two ordinals in $\dom(\mathcal{M})$, we let $\mathcal{M}\restr[\varepsilon_0,\varepsilon_1)$ denote the set of models $Q\in\mathcal{M}$ such that $\varepsilon_0\leq\varepsilon_Q<\varepsilon_1$. Lastly, we let $\mathcal{M}[\omega_1]$ denote the set $\{Q[\omega_1]:Q\in\mathcal{M}\}$.
	
	\begin{definition}\label{defSSM2T}
		Let $\mathcal{M}$ be a finite set of members of $H(\kappa)$. We say that $\mathcal{M}$ is an \emph{$(\mathcal{S},\mathcal{L})$-symmetric system} if and only if the following holds:
		\begin{enumerate}[label=(\Alph*)]
			\item Every $Q\in\mathcal{M}$ is an element of $\mathcal{S}\cup\mathcal{L}$.
			
			\item For any two $Q_0,Q_1\in\mathcal{M}$, if $\varepsilon_{Q_0}=\varepsilon_{Q_1}$, then $Q_0[\omega_1]\cong Q_1[\omega_1]$. 
			
			\item If $\varepsilon_0$ is the immediate predecessor of $\varepsilon_1$ in $\dom(\mathcal{M})$, then for every $Q_0\in\mathcal{M}(\varepsilon_0)$, there is some $Q_1\in\mathcal{M}(\varepsilon_1)$ such that $Q_0\in Q_1$.
			
			\item For all $Q_0,Q_1,Q_1'\in\mathcal{M}$ such that $Q_0\in Q_1$ and $\varepsilon_{Q_1}=\varepsilon_{Q_1'}$, $$\Psi_{Q_1[\omega_1],Q_1'[\omega_1]}(Q_0)\in\mathcal{M}.$$
			
			\item For every $X\in\mathcal{M}\cap\mathcal{L}$ and every $M\in\mathcal{M}\cap\mathcal{S}$, if $X\in M$, then $X\cap M\in\mathcal{M}$.
		\end{enumerate}
	\end{definition}
	
	\begin{remark}\label{remark-lie-below}
		If $\mathcal{M}$ is an $(\mathcal{S},\mathcal{L})$-symmetric system, for any two $Q,P\in\mathcal{M}$, having $Q\in P[\omega_1]$ is equivalent to the existence of a maximal $\in$-chain of models that connects $Q$ and $P$ in $\mathcal{M}$. More precisely, there exists a sequence $\la R_i:i\leq n\ra$ in $\mathcal{M}$, weakly linearly ordered by $\in$,\footnote{A binary relation $<$ is a \emph{weak linear order} on a set $L$ if the transitive closure of $<$ is a linear order on $L$.} such that $R_0=Q$ and $R_n=P$, and is maximal in the sense that
		\[
		\{\varepsilon_{R_i}:i\leq n\}=\dom(\mathcal{M})\cap[\varepsilon_Q,\varepsilon_{P}].
		\]
		In this situation, we will say that \emph{$Q$ lies below $P$.} In particular, $\mathcal{M}\cap P[\omega_1]$ is the set of all members of $\mathcal{M}$ that can reach $P$ by a maximal $\in$-chain of members of $\mathcal{M}$ in the above sense.
		
		Note that for any two $Q_0,Q_1\in\mathcal{M}$, if $Q_0\in Q_1$, then $Q_0\in Q_1[\omega_1]$. However, the converse is not in general true (see Proposition \ref{propgaps}). 
	\end{remark}
	
	\begin{remark}
		If $\mathcal{M}$ is an $(\mathcal{S},\mathcal{L})$-symmetric system and $\varepsilon_0<\varepsilon_1$ belong to $\dom(\mathcal{M})$, then for every $Q_0\in\mathcal{M}(\varepsilon_0)$, we can find a model $Q_1\in\mathcal{M}(\varepsilon_1)$ such that $Q_0\in Q_1[\omega_1]$, by successive applications of clause (C).
	\end{remark}
	
	We will refer to clause (C) and the conclusion of the last remark, indistinctly, as the \emph{shoulder axiom for $\mathcal{M}$}. 
	
	Let us start by proving some basic facts about $(\mathcal{S},\mathcal{L})$-symmetric systems. 
	
	\begin{proposition}\label{prop18}
		Let $\mathcal{M}$ be an $(\mathcal{S},\mathcal{L})$-symmetric system. The following hold for any two models $Q_0,Q_1\in\mathcal{M}$:
		\begin{enumerate}[label=(\arabic*)]
			\item If $Q_0\in Q_1$ and $|Q_0|\leq|Q_1|$, then $Q_0\subseteq Q_1$.
			\item If $Q_1$ is a small model, $Q_0\in Q_1[\omega_1]$, and there is no large model $X\in\mathcal{M}$ such that $\varepsilon_{Q_0}<\varepsilon_X<\varepsilon_{Q_1}$, then $Q_0\in Q_1$.
		\end{enumerate}
	\end{proposition}
	\begin{proof}
		The first item is an immediate consequence of the basic facts about elementary submodels from Section \ref{section-prelim}. Let us show (2). If $\varepsilon_{Q_1}$ is the successor of $\varepsilon_{Q_0}$ in $\dom(\mathcal{M})$, there is $Q_1'\in\mathcal{M}(\varepsilon_{Q_1})$ such that $Q_0\in Q_1'$, by the shoulder axiom. Therefore, since $Q_0\in Q_1'\cap Q_1[\omega_1]$, by (B), $$Q_0=\Psi_{Q_1'[\omega_1],Q_1[\omega_1]}(Q_0)\in Q_1.$$
		Suppose now that $\dom(\mathcal{M})\cap(\varepsilon_{Q_0},\varepsilon_{Q_1})$ is nonempty, and let $\la\varepsilon_i:i\leq n\ra$ enumerate this set in increasing order. By $n+2$ applications of the shoulder axiom, we may find $M_0,\dots, M_n, Q_1'\in\mathcal{M}$ such that 
		\begin{itemize}
			\item $Q_0\in M_0\in\dots\in M_n\in Q_1'$,
			\item $\varepsilon_{M_i}=\varepsilon_i$, for all $i\leq n$, and
			\item $\varepsilon_{Q_1'}=\varepsilon_{Q_1}$.
		\end{itemize}
		Since $M_0,\dots,M_n$ are small models by hypothesis, item (1) implies that $Q_0\in Q_1'$. Therefore, $Q_0\in Q_1$, by the exact same argument as above, using (B). 
	\end{proof}
	
	The following result strengthens clause (D) from the definition of $(\mathcal{S},\mathcal{L})$-symmetric system. Throughout the paper, when we refer to the \emph{symmetry} of an $(\mathcal{S},\mathcal{L})$-symmetric system, we will usually mean the following stronger form of closure under isomorphisms within an $(\mathcal{S},\mathcal{L})$-symmetric system.
	
	\begin{proposition}\label{prop6}
		Let $\mathcal{M}$ be an $(\mathcal{S},\mathcal{L})$-symmetric system, and let $Q_0,Q_1,Q_1'$ be in $\mathcal{M}$. If $Q_0\in Q_1[\omega_1]$ and $\varepsilon_{Q_1}=\varepsilon_{Q_1'}$, then $\Psi_{Q_1[\omega_1],Q_1'[\omega_1]}(Q_0)\in\mathcal{M}$.
	\end{proposition}
	\begin{proof}
		If $Q_1,Q_1'$ are large models, the conclusion follows trivially. Hence, we may assume that $Q_1,Q_1'$ are countable. If there are no large models $X\in\mathcal{L}$ such that $\varepsilon_{Q_0}<\varepsilon_X<\varepsilon_{Q_1}$, then $Q_0\in Q_1$ by the last proposition, and the conclusion follows again trivially. Suppose now the contrary, and let $\varepsilon$ be the greatest $\omega_2$-height of any large model in $\mathcal{M}$ lying strictly between $\varepsilon_{Q_0}$ and $\varepsilon_{Q_1}$. By two applications of the shoulder axiom, we can find two models $X,Q_1^*\in\mathcal{M}$ such that $\varepsilon_X=\varepsilon$, $\varepsilon_{Q_1^*}=\varepsilon_{Q_1}$, $Q_0\in X$, and $X\in Q_1^*[\omega_1]$. By the last proposition and the maximality of $\varepsilon$, we have $X\in Q_1^*$. Apply clause (D) from Definition \ref{defSSM2T} twice to get $X_1:=\Psi_{Q_1^*[\omega_1],Q_1[\omega_1]}(X)$ and $X_1':=\Psi_{Q_1[\omega_1],Q_1'[\omega_1]}(X_1)$ in $\mathcal{M}$. Since $Q_0\in Q_1^*[\omega_1]\cap Q_1[\omega_1]$, we have 
		\[
		Q_0=\Psi_{Q_1^*[\omega_1],Q_1[\omega_1]}(Q_0)=\Psi_{X,X_1}(Q_0)\in X_1,
		\]
		Therefore, by clause (D), we can conclude that
		\[
		\Psi_{Q_1[\omega_1],Q_1'[\omega_1]}(Q_0)=\Psi_{X_1,X_1'}(Q_0)\in\mathcal{M}.
		\]
	\end{proof}
	
	\begin{proposition}\label{prop13}
		Let $\mathcal{M}$ be an $(\mathcal{S},\mathcal{L})$-symmetric system and let $Q_0,Q_1$ be two models in $\mathcal{M}$ such that $Q_0\in Q_1[\omega_1]$. If there is $\varepsilon\in\dom(\mathcal{M})$ such that $\varepsilon_{Q_0}<\varepsilon<\varepsilon_{Q_1}$, then there is some $P\in\mathcal{M}(\varepsilon)$ such that $Q_0\in P[\omega_1]$ and $P\in Q_1[\omega_1]$. 
	\end{proposition}
	\begin{proof}
		Apply the shoulder axiom twice to find $P',Q_1'\in\mathcal{M}$ such that $\varepsilon_{P'}=\varepsilon$, $\varepsilon_{Q_1'}=\varepsilon_{Q_1}$, $Q_0\in P'[\omega_1]$ and $P'\in Q_1'[\omega_1]$. Let us denote $P:=\Psi_{Q_1'[\omega_1],Q_1[\omega_1]}(P')$, and note that $P\in\mathcal{M}(\varepsilon)\cap Q_1[\omega_1]$ by the symmetry of $\mathcal{M}$. Hence, since $Q_0\in Q_1'[\omega_1]\cap Q_1[\omega_1]$, by clause (B) of the definition of $(\mathcal{S},\mathcal{L})$-symmetric system,
		\[
		Q_0=\Psi_{Q_1'[\omega_1],Q_1[\omega_1]}(Q_0)=\Psi_{P'[\omega_1],P[\omega_1]}(Q_0)\in P[\omega_1].
		\]
	\end{proof}
	
	
	The following proposition, which tells us that an $\omega_1$-isomorphic copy of an $(\mathcal{S},\mathcal{L})$-symmetric system is again an $(\mathcal{S},\mathcal{L})$-symmetric system, is an easy exercise.
		
	\begin{proposition}\label{isocopy}
		Let $\mathcal{M}$ be an $(\mathcal{S},\mathcal{L})$-symmetric system. Suppose that $Q_0,Q_1\in\mathcal{S}\cup\mathcal{L}$ are $\omega_1$-isomorphic. If $\mathcal{M}\subseteq Q_0[\omega_1]$, then $\Psi_{Q_0[\omega_1],Q_1[\omega_1]}(\mathcal{M})$ is an $(\mathcal{S},\mathcal{L})$-symmetric system.
	\end{proposition}


	\subsection{Amalgamation lemmas}\label{subsection-amalgamation-lemmas} In this subsection we will prove some standard results, in the context of forcing with side conditions, for amalgamating $(\mathcal{S},\mathcal{L})$-symmetric systems. As we shall see in the next subsection, the preservation theorems for the pure side condition forcing follow, by standard arguments, from these amalgamation lemmas.
	
	\begin{lemma}\label{amal-ontop}
		Let $\mathcal{M}$ be an $(\mathcal{S},\mathcal{L})$-symmetric system and let $Q\in\mathcal{S}\cup\mathcal{L}$ be such that $\mathcal{M}\subseteq Q$. Then there is an $(\mathcal{S},\mathcal{L})$-symmetric system $\mathcal{M}_Q$ such that $\mathcal{M}\cup\{Q\}\subseteq\mathcal{M}_Q$ 
	\end{lemma}
	\begin{proof}
		If $Q$ is a large model, it is easy to check that $\mathcal{M}_Q=\mathcal{M}\cup\{Q\}$ is an $(\mathcal{S},\mathcal{L})$-symmetric system. So, suppose that $Q$ is a small model, and let us denote it by $M$. We claim that
		\[
		\mathcal{M}_M=\mathcal{M}\cup\{M\}\cup\{X\cap M:X\in\mathcal{M}\cap\mathcal{L}\}
		\]
		is an $(\mathcal{S},\mathcal{L})$-symmetric system.
		
		In order to see that $\mathcal{M}_M$ satisfies clause (B) of Definition \ref{defSSM2T}, we only need to note that for every $X_0,X_1\in\mathcal{M}\cap\mathcal{L}$ such that $\varepsilon_{X_0}=\varepsilon_{X_1}$, we have $(X_0\cap M)[\omega_1]\cong (X_1\cap M)[\omega_1]$ by Proposition \ref{prop9}.
		
		Clause (C) follows from the fact that if $X\in\mathcal{M}\cap\mathcal{L}$ and $P\in\mathcal{M}\cap X$, then $P\in X\cap M$.
		
		Let us now show clause (D). Let $P,Q_0,Q_1\in\mathcal{M}_M$ be such that $P\in Q_0$ and $\varepsilon_{Q_0}=\varepsilon_{Q_1}$. We need to check that $\Psi_{Q_0[\omega_1],Q_1[\omega_1]}(P)\in\mathcal{M}_M$. There are three relevant cases:
		
		\textbf{Case 1:} $P=X\cap M$ for some $X\in\mathcal{M}\cap\mathcal{L}$ and for each $i\in\{0,1\}$, $Q_i=X_i\cap M$ for some $X_i\in\mathcal{M}\cap\mathcal{L}$. Note that this case is impossible because $X\cap M$ cannot be an element of $M$.
		
		\textbf{Case 2:} $P\in\mathcal{M}$, and for each $i\in\{0,1\}$, $Q_i=X_i\cap M$ for some large model $X_i$ in  $\mathcal{M}$. Since $(X_0\cap M)[\omega_1]\cong(X_1\cap M)[\omega_1]$ is witnessed by the isomorphism $\Psi_{X_0,X_1}$ by Proposition \ref{prop9}, we have $\Psi_{Q_0[\omega_1],Q_1[\omega_1]}(P)=\Psi_{X_0,X_1}(P)$, which is an element of $\mathcal{M}$ by symmetry.
		
		\textbf{Case 3:} $P=X\cap M$ for some $X\in\mathcal{M}\cap\mathcal{L}$, and $Q_0,Q_1\in\mathcal{M}$. If $\varepsilon_{Q_0}=\varepsilon_{Q_1}=\varepsilon_X$, as $P\in X\cap Q_0$, we have $$X\cap M=\Psi_{X,Q_0}(X\cap M)=Q_0\cap M$$ by clause (B). Therefore, by Proposition \ref{prop9},
		\[
		\Psi_{Q_0,Q_1}(X\cap M)=\Psi_{Q_0,Q_1}(Q_0\cap M)=Q_1\cap M,
		\]
		which is an element of $\mathcal{M}_M$ by construction. Suppose now that $\varepsilon_X<\varepsilon_{Q_0}=\varepsilon_{Q_1}$. Using the shoulder axiom for $\mathcal{M}$, find a model $Q_2\in\mathcal{M}(\varepsilon_{Q_0})$ such that $X\in Q_2[\omega_1]$. Let $Y:=\Psi_{Q_2[\omega_1],Q_0[\omega_1]}(X)$, which is an element of $\mathcal{M}$ by symmetry, and note that, as $P\in Q_2[\omega_1]\cap Q_0$, by clause (B) we have $P=\Psi_{Q_2[\omega_1],Q_0[\omega_1]}(P)\in Y$. Hence, by Proposition \ref{prop9}, and since $P\in X\cap Y$, $$P=X\cap M=Y\cap M.$$ Let now $Z:=\Psi_{Q_0[\omega_1],Q_1[\omega_1]}(Y)$, which is an element of $\mathcal{M}$ by symmetry, and note that 
		\begin{align*}
			\Psi_{Q_0[\omega_1],Q_1[\omega_1]}(P)&=\Psi_{Q_0[\omega_1],Q_1[\omega_1]}(Y\cap M)\\
			&=\Psi_{Y,Z}(Y\cap M)=Z\cap M,
		\end{align*}
		which is an element of $\mathcal{M}_M$ by construction.
		
		Lastly, to check clause (E), it suffices to note that if $X,Y$ are two large models in $\mathcal{M}$ such that $X\in(Y\cap M)$, then $X\cap(Y\cap M)=X\cap M$.
	\end{proof}
	
	If $\mathcal{M}$ is an $(\mathcal{S},\mathcal{L})$-symmetric system and $M$ is a small model in $\mathcal{M}$, we will denote by $\vect{\mathcal{E}}_{\mathcal{M},M}^\mathcal{L}$ the increasing enumeration $\langle\varepsilon_i^\mathcal{L}:i<n\rangle$ of $\dom(\mathcal{M}\cap M\cap\mathcal{L})$. Moreover, we will denote by $\vect{\mathcal{E}}_{\mathcal{M},M}^\mathcal{S}$ the increasing enumeration $\langle\varepsilon_i^\mathcal{S}:i<n\rangle$ of the set $$\{\varepsilon_{X\cap M}:X\in\mathcal{M}\cap M\cap\mathcal{L}\}.$$
	
	Note that if $X_0\in\dots\in X_{n-1}$ and $Y_0\in\dots\in Y_{n-1}$ are two chains in $\mathcal{M}\cap M\cap\mathcal{L}$ such that, for every $i<n$, $\varepsilon_{X_i}=\varepsilon_{Y_i}=\varepsilon_i^\mathcal{L}$, then, in light of Proposition \ref{prop9}, $\varepsilon_{X_i\cap M}=\varepsilon_{Y_i\cap M}=\varepsilon_i^\mathcal{S}$. 
	
	\begin{proposition}\label{prop-res-gaps}
		Let $\mathcal{M}$ be an $(\mathcal{S},\mathcal{L})$-symmetric system and $M\in\mathcal{M}\cap\mathcal{S}$. Suppose that $\vect{\mathcal{E}}_{\mathcal{M},M}^\mathcal{L}=\langle\varepsilon_i^\mathcal{L}:i<n\rangle$ and $\vect{\mathcal{E}}_{\mathcal{M},M}^\mathcal{S}=\langle\varepsilon_i^\mathcal{S}:i<n\rangle$, for some $n<\omega$. Then, $\dom(\mathcal{M}\cap M)=\mathcal{E}_0\cup\mathcal{E}_1\cup\mathcal{E}_2$, where
		\begin{itemize}
			\item $\mathcal{E}_0=\dom(\mathcal{M})\cap\varepsilon_0^\mathcal{S}$,
			\item $\mathcal{E}_1=\bigcup_{i<n-1}\big(\dom(\mathcal{M})\cap[\varepsilon_i^\mathcal{L},\varepsilon_{i+1}^\mathcal{S})\big)$, and
			\item $\mathcal{E}_2=\dom(\mathcal{M})\cap[\varepsilon_{n-1}^\mathcal{L},\varepsilon_M)$.
		\end{itemize}
	\end{proposition}
	\begin{proof}
		Let us first show the inclusion $\dom(\mathcal{M}\cap M)\subseteq\mathcal{E}_0\cup\mathcal{E}_1\cup\mathcal{E}_2$. Let $\varepsilon\in\dom(\mathcal{M}\cap M)$ and suppose, towards a contradiction, that $\varepsilon\notin\mathcal{E}_0\cup\mathcal{E}_1\cup\mathcal{E}_2$. Then $$\varepsilon\in\dom(\mathcal{M})\cap \bigcup_{i<n}[\varepsilon_i^\mathcal{S},\varepsilon_i^\mathcal{L}),$$ but this is impossible by Proposition \ref{propgaps}.
		
		Let us now prove the other inclusion. It is immediate to see, by definition of $\vect{\mathcal{E}}_{\mathcal{M},M}^\mathcal{L}$ and $\vect{\mathcal{E}}_{\mathcal{M},M}^\mathcal{S}$, that every model $Q\in\mathcal{M}$ such that either
		\begin{itemize}
			\item $\varepsilon_Q<\varepsilon_0^\mathcal{S}$,
			\item $\varepsilon_Q\in(\varepsilon_i^\mathcal{L},\varepsilon_{i+1}^\mathcal{S})$ for some $i<n-1$, or
			\item $\varepsilon_Q\in (\varepsilon_{n-1}^\mathcal{L},\varepsilon_M)$,
		\end{itemize} 
		is a small model. Let $P\in\mathcal{M}$ be such that $\varepsilon_P\in[\varepsilon_i^\mathcal{L},\varepsilon_{i+1}^\mathcal{S})$ for some $i<n-1$. By the symmetry of $\mathcal{M}$, we may further assume that $P$ is a member of $(X_{i+1}\cap M)[\omega_1]$, for some $X_{i+1}\in\mathcal{M}\cap M\cap\mathcal{L}$ such that $\varepsilon_{X_{i+1}}=\varepsilon_{i+1}^\mathcal{L}$ (and hence $\varepsilon_{X_{i+1}\cap M}=\varepsilon_{i+1}^\mathcal{S}$). Therefore, $P\in X_{i+1}\cap M$ by Proposition \ref{prop18}. Similarly, if $\varepsilon_P<\varepsilon_0^\mathcal{S}$ and $P\in (X_0\cap M)[\omega_1]$ for some $X_0\in\mathcal{M}\cap M\cap\mathcal{L}$ such that $\varepsilon_{X_0}=\varepsilon_0^\mathcal{L}$, or $\varepsilon_P\in[\varepsilon_{n-1}^\mathcal{L},\varepsilon_M)$ and $P\in M[\omega_1]$, in both cases $P$ must be a member of $M$. So, we can conclude that $\varepsilon_P\in\dom(\mathcal{M}\cap M)$.
	\end{proof}
	
	Recall from Remark \ref{remark-lie-below} that if $Q$ is either a large model or of the form $M[\omega_1]$, for some small $M\in\mathcal{M}$, then $\mathcal{M}\cap Q$ coincides with the set of all models in $\mathcal{M}$ lying below $Q$. However, if $Q$ is a small model, there may be models in $\mathcal{M}$ lying below $Q$ which are not members of $Q$. 
	
	The next proposition, that extends Proposition \ref{prop-res-gaps}, offers a very precise picture of the geometry of $\mathcal{M}\cap Q$. Although it might not be explicitly referenced, we encourage the reader to keep it mind for a better understanding of the proofs of the next amalgamation lemmas.
	
	
	
	\begin{proposition}\label{prop-res-gaps2}
		Let $\mathcal{M}$ be an $(\mathcal{S},\mathcal{L})$-symmetric system and $M\in\mathcal{M}\cap\mathcal{S}$. Suppose that $$\vect{\mathcal{E}}_{\mathcal{M},M}^\mathcal{L}=\langle\varepsilon_i^\mathcal{L}:i<n\rangle$$ and $$\vect{\mathcal{E}}_{\mathcal{M},M}^\mathcal{S}=\langle\varepsilon_i^\mathcal{S}:i<n\rangle,$$ for some $n<\omega$. Let $Q\in\mathcal{M}\cap M[\omega_1]$. Then, the following hold:
		\begin{enumerate}
			\item If $\varepsilon_i^\mathcal{S}\leq\varepsilon_Q<\varepsilon_i^\mathcal{L}$ for some $i<n$, then $Q$ is not a member of $M$.
			
			\item If $\varepsilon_Q<\varepsilon_0^\mathcal{S}$, then $Q$ is a small model and there is $P\in\mathcal{M}\cap M$ such that $\varepsilon_{P}=\varepsilon_Q$. Moreover, there is some $X_0\in\mathcal{M}\cap M$ such that $\varepsilon_{X_0}=\varepsilon_0^\mathcal{L}$ and $P\in X_0\cap M$.
			
			\item If $\varepsilon_i^\mathcal{L}\leq\varepsilon_Q<\varepsilon_{i+1}^\mathcal{S}$ for some $i<n-1$, then:
			\begin{enumerate}
				\item Either $\varepsilon_Q=\varepsilon_i^\mathcal{L}$ and $Q$ is a large model, or $\varepsilon_i^\mathcal{L}<\varepsilon_Q<\varepsilon_{i+1}^\mathcal{S}$ and $Q$ is a small model.
				
				\item There is $P\in\mathcal{M}\cap M$ such that $\varepsilon_{P}=\varepsilon_Q$. Moreover, there is some $X_{i+1}\in\mathcal{M}\cap M$ such that $\varepsilon_{X_{i+1}}=\varepsilon_{i+1}^\mathcal{L}$ and $P\in X_{i+1}\cap M$.
			\end{enumerate}
			
			\item If $\varepsilon_{n-1}^\mathcal{L}\leq\varepsilon_Q<\varepsilon_M$, then $Q\in M$ and either $\varepsilon_Q=\varepsilon_{n-1}^\mathcal{L}$ and $Q$ is a large model, or $\varepsilon_{n-1}^\mathcal{L}<\varepsilon_Q<\varepsilon_M$ and $Q$ is a small model.
		\end{enumerate}
	\end{proposition}
	\begin{proof}
		Item (1) follows directly from Proposition \ref{propgaps}. Let us move to the proof of item (3), and skip the proofs of items (2) and (4), which are very similar but simpler. 
		
		Since $\vect{\mathcal{E}}_{\mathcal{M},M}^\mathcal{L}=\langle\varepsilon_i^\mathcal{L}:i<n\rangle$ enumerates $\dom(\mathcal{M}\cap M\cap\mathcal{L})$, clause (a) is straightforward. Let us show clause (b). Let $N_{i+1}\in\mathcal{M}$ be such that $\varepsilon_{N_{i+1}}=\varepsilon_{i+1}^\mathcal{S}$ and $Q\in N_{i+1}[\omega_1]$, given by the shoulder axiom. Since all models lying between $Q$ and $N_{i+1}$ are small by clause (a), Proposition \ref{prop18} ensures that $Q\in N_{i+1}$. Now, let $X_{i+1}\in\mathcal{M}\cap M$ be such that $\varepsilon_{X_{i+1}}=\varepsilon_{i+1}^\mathcal{L}$. Then, by symmetry, $P:=\Psi_{N_{i+1}[\omega_1],(X_{i+1}\cap M)[\omega_1]}(Q)$ is a member of $\mathcal{M}\cap (X_{i+1}\cap M)$.
	\end{proof}
	
	As a summary, the last proposition tells us that either $Q$ is such that $\varepsilon_i^\mathcal{S}\leq\varepsilon_Q<\varepsilon_i^\mathcal{L}$ for some $i<n$, and then $Q$ is not a member of $M$, or $Q$ is not such that $\varepsilon_i^\mathcal{S}\leq\varepsilon_Q<\varepsilon_i^\mathcal{L}$ for some $i<n$, and then, even though $Q$ might not be a member of $M$, there is an isomorphic copy of $Q$ which is.
	
	\begin{lemma}\label{amal-rest}
		Let $\mathcal{M}$ be an $(\mathcal{S},\mathcal{L})$-symmetric system and let $Q$ be a member of $\mathcal{M}\cup\mathcal{M}[\omega_1]$. Then $\mathcal{M}\cap Q$ is an $(\mathcal{S},\mathcal{L})$-symmetric system.
	\end{lemma}
	\begin{proof}
		If $Q\in\mathcal{L}$ or $Q\in\mathcal{M}[\omega_1]$, then $\mathcal{M}\cap Q$ coincides with the set of all models in $\mathcal{M}$ lying below $Q$, and the verification of clauses (A)-(E) from Definition \ref{defSSM2T} follows from a simpler version of the argument that we will present for $Q$ small.
		
		So, let us assume that $Q$ is a small model and denote it by $M$. Clauses (A), (B), (D), and (E) are trivial or easy to check. Let us show clause (C). Let $Q_0\in\mathcal{M}\cap M$ and suppose that $\varepsilon_1$ is the immediate successor of $\varepsilon_0:=\varepsilon_{Q_0}$ in $\dom(\mathcal{M}\cap M)$. We need to find $Q_1\in\mathcal{M}\cap M$ such that $\varepsilon_{Q_1}=\varepsilon_1$ and $Q_0\in Q_1$. Suppose that $$\vect{\mathcal{E}}_{\mathcal{M},M}^\mathcal{L}=\langle\varepsilon_i^\mathcal{L}:i<n\rangle$$ and $$\vect{\mathcal{E}}_{\mathcal{M},M}^\mathcal{S}=\langle\varepsilon_i^\mathcal{S}:i<n\rangle,$$ for some $n<\omega$. We will divide the proof into two cases, based on the decomposition of $\dom(\mathcal{M}\cap M)$ given by Proposition \ref{prop-res-gaps}:
		
		\textbf{Case 1.} Suppose first that $\varepsilon_0,\varepsilon_1\in[\varepsilon_i^\mathcal{L},\varepsilon_{i+1}^\mathcal{S})$ for some $i<n-1$. If we can find a model $X_{i+1}\in\mathcal{M}\cap M\cap\mathcal{L}$ such that $Q_0\in X_{i+1}\cap M$ and $\varepsilon_{X_{i+1}}=\varepsilon_{i+1}^\mathcal{L}$, then by Proposition \ref{prop13} there is $Q_1\in\mathcal{M}$ such that 
		\begin{itemize}
			\item $\varepsilon_{Q_1}=\varepsilon_1$,
			
			\item $Q_0\in Q_1[\omega_1]$, and
			
			\item $Q_1\in (X_{i+1}\cap M)[\omega_1]$,
		\end{itemize}
		and from this, it is easy to see that $Q_0\in Q_1\in X_{i+1}\cap M$. 
		
		In order to find a model $X_{i+1}$ as above, we need to perform a recursive construction of an $\in$-chain
		\[
		Q_0\in X_{i+1}\cap M\in X_{i+1}\in\dots\in X_j\cap M\in X_j\in\dots\in X_{n-1}\cap M\in X_{n-1}\in M,
		\]
		where $X_{i+1},\dots,X_j,\dots,X_{n-1}$ are all in $\mathcal{M}\cap\mathcal{L}\cap M$, and
		\[
		\varepsilon_{X_{i+1}}=\varepsilon_{i+1}^{\mathcal{L}},\dots,\varepsilon_{X_{j}}=\varepsilon_{j}^{\mathcal{L}},\dots,\varepsilon_{X_{n-1}}=\varepsilon_{n-1}^{\mathcal{L}}.
		\]
		For the pair of models $(Q_0,M)$, use Proposition \ref{prop13} to find a large $X_{n-1}\in\mathcal{M}$ such that $Q_0\in X_{n-1}\in M[\omega_1]$ and $\varepsilon_{X_{n-1}}=\varepsilon_{n-1}^{\mathcal{L}}$, and note that $X_{n-1}\in M$ by Proposition \ref{prop18}. Therefore, $Q_0\in X_{n-1}\cap M\in\mathcal{M}\cap X_{n-1}$. By recursively applying the same argument to the pairs $(Q_0,X_j\cap M)$, for $i+2\leq j<n$, we construct an $\in$-chain with the properties listed above and, in particular, we obtain a model $X_{i+1}\in\mathcal{M}\cap M\cap\mathcal{L}$ such that $Q_0\in X_{i+1}\cap M$ and $\varepsilon_{X_{i+1}}=\varepsilon_{i+1}^\mathcal{L}$, as we wanted.		
		
		The cases $\varepsilon_0,\varepsilon_1<\varepsilon_0^\mathcal{S}$ and $\varepsilon_0,\varepsilon_1\in[\varepsilon_{n-1}^\mathcal{L},\varepsilon_M)$ can be argued similarly. 
		
		\textbf{Case 2.} Suppose now that $\varepsilon_0<\varepsilon_i^\mathcal{S}<\varepsilon_i^\mathcal{L}\leq\varepsilon_1$, for some $i<n$. Note that then $\varepsilon_0=\max(\dom(\mathcal{M})\cap\varepsilon_i^\mathcal{S})$ and $\varepsilon_1=\varepsilon_i^\mathcal{L}$. As in the last case, we may find some $X_i\in\mathcal{M}\cap M\cap\mathcal{L}$ such that $Q_0\in X_i\cap M$ and $\varepsilon_{X_i}=\varepsilon_i^\mathcal{L}$. Hence, we are done simply by letting $Q_1=X_i$.
	\end{proof}
	
	\begin{notation}
		If $\mathcal{M}$ is an $(\mathcal{S},\mathcal{L})$-symetric system, recall that we let $\mathcal{M}[\omega_1]$ denote the set $\{Q[\omega_1]:Q\in\mathcal{M}\}$. If $\mathcal{M}_0$ and $\mathcal{M}_1$ are two $(\mathcal{S},\mathcal{L})$-symmetric systems, we write $\mathcal{M}_0\cong\mathcal{M}_1$ to denote that, for some $n<\omega$, there are enumerations $\langle Q_i^0:i<n\rangle$ and $\langle Q_i^1:i<n\rangle$ of $\mathcal{M}_0$ and $\mathcal{M}_1$, respectively, such that  $(\bigcup\mathcal{M}_0[\omega_1];\in,Q_i^0,\lhd)_{i<n}$ and $(\bigcup\mathcal{M}_1[\omega_1];\in,Q_i^1,\lhd)_{i<n}$ are isomorphic via an isomorphism which is the identity on $\bigcup\mathcal{M}_0[\omega_1]\cap\bigcup\mathcal{M}_1[\omega_1]$.
	\end{notation}
	
	\begin{lemma}\label{pureamalgamation2}
		Let $n<\omega$ and suppose that $\mathcal{M}_0,\dots,\mathcal{M}_n$ are $(\mathcal{S},\mathcal{L})$-symmetric systems such that $\mathcal{M}_i\cong\mathcal{M}_j$, for all $i,j\leq n$. Then $\bigcup_{i\leq n}\mathcal{M}_i$ is an $(\mathcal{S},\mathcal{L})$-symmetric system.
	\end{lemma}
	\begin{proof}
		Obviously, $\bigcup_{i\leq n}\mathcal{M}_i$ is a finite subset of $\mathcal{S}\cup\mathcal{L}$. We will check that this union of $(\mathcal{S},\mathcal{L})$-symmetric systems satisfies clauses (B)-(E) from Definition \ref{defSSM2T}. For any two $i,j\leq n$, denote the isomorphism witnessing $\mathcal{M}_i\cong\mathcal{M}_j$ by $\Psi_{i,j}$, and the intersection $\bigcup\mathcal{M}_i[\omega_1]\cap\bigcup\mathcal{M}_j[\omega_1]$ by $\Delta_{i,j}$. 
		
		Let us start by showing clause (B). Let $Q_0\in\mathcal{M}_i$ and $Q_1\in\mathcal{M}_j$ be two models such that $\varepsilon_{Q_0}=\varepsilon_{Q_1}$, for some $i,j\leq n$. If we let $Q_2$ denote the model $\Psi_{i,j}(Q_0)$ of $\mathcal{M}_j$, then $$(\Psi_{Q_2[\omega_1],Q_1[\omega_1]}\circ\Psi_{i,j})\restr Q_0[\omega_1]:(Q_0[\omega_1];\in,Q_0,\lhd)\to(Q_1[\omega_1];\in,Q_1,\lhd)$$ is an isomorphism. Suppose now that $x\in Q_0[\omega_1]\cap Q_1[\omega_1]$. Then $x$ must be an element of $\Delta_{i,j}$, and hence $\Psi_{i,j}(x)=x$. In particular, this implies that $x\in Q_2[\omega_1]$, and as $\mathcal{M}_j$ is an $(\mathcal{S},\mathcal{L})$-symmetric system, $\Psi_{Q_2[\omega_1],Q_1[\omega_1]}(x)=x$. Hence, $(\Psi_{Q_2[\omega_1],Q_1[\omega_1]}\circ\Psi_{i,j})\restr Q_0[\omega_1]$ is an $\omega_1$-isomorphism between $Q_0[\omega_1]$ and $Q_1[\omega_1]$.
		
		To show clause (C), it suffices to note that the existence of $\Psi_{i,j}$ implies that $\dom(\mathcal{M}_i)=\dom(\mathcal{M}_j)$, for any two $i,j\leq n$.
		
		Let $i_0,i_1,i_2\leq n$ and let $Q_0\in\mathcal{M}_{i_0}$, $Q_1\in\mathcal{M}_{i_1}$ and $Q_2\in\mathcal{M}_{i_2}$, such that $Q_0\in Q_1$ and $\varepsilon_{Q_1}=\varepsilon_{Q_2}$. In order to show clause (D), we must verify that $\Psi_{Q_1[\omega_1],Q_2[\omega_1]}(Q_0)$ belongs to $\bigcup_{i\leq n}\mathcal{M}_i$. Since $Q_0\in Q_1$, it follows that $Q_0\in \Delta_{i_0,i_1}$, and hence that $Q_0=\Psi_{i_0,i_1}(Q_0)\in\mathcal{M}_{i_1}$. Let $Q_0'=\Psi_{i_1,i_2}(Q_0)$ and $Q_1'=\Psi_{i_1,i_2}(Q_1)$, which are both elements of $\mathcal{M}_{i_2}$. Then $$\Psi_{Q_1[\omega_1],Q_2[\omega_1]}(Q_0)=\Psi_{Q_1'[\omega_1],Q_2[\omega_1]}(Q_0'),$$ which is an element of $\mathcal{M}_{i_2}$ because $Q_0',Q_1',Q_2\in\mathcal{M}_{i_2}$ and $\mathcal{M}_{i_2}$ is an $(\mathcal{S},\mathcal{L})$-symmetric system.
		
		Lastly, we will show clause (E). Let $X\in\mathcal{M}_i\cap\mathcal{L}$ and $M\in\mathcal{M}_j\cap\mathcal{S}$ be such that $X\in M$, for some $i,j\leq n$. We need to check that $X\cap M\in\bigcup_{i\leq n}\mathcal{M}_i$. Note that as $X\in\mathcal{M}_i\cap M$, in particular $X\in\Delta_{i,j}$. Hence, $X$ is fixed by the isomorphism $\Psi_{i,j}$. Therefore, $X,M\in\mathcal{M}_j$, and so $X\cap M\in\mathcal{M}_j$.
	\end{proof}
	
	\begin{definition}
		Let $\mathcal{M},\mathcal{W}$ be two $(\mathcal{S},\mathcal{L})$-symmetric systems. We denote by $\mathcal{W}\land\mathcal{M}$ the set of all models of the form $X\cap M$, where $M\in\mathcal{M}\cap\mathcal{S}$ and $X\in\mathcal{W}\cap\mathcal{L}\cap M$.
	\end{definition}
	
	\begin{proposition}\label{prop-pure-amalgamation-1}
		If $\mathcal{M}$ is an $(\mathcal{S},\mathcal{L})$-symmetric system, $X\in\mathcal{M}$ is a large model, and $\mathcal{W}$ is an $(\mathcal{S},\mathcal{L})$-symmetric system such that $\mathcal{M}\cap X\subseteq\mathcal{W}\subseteq X$, then $\mathcal{W}\land\mathcal{M}\subseteq\mathcal{W}$.
	\end{proposition}
	\begin{proof}
		Let $Y\in\mathcal{W}\cap\mathcal{L}$ and $M\in\mathcal{M}\cap\mathcal{S}$ be such that $Y\in M$. We have to show that $Y\cap M\in\mathcal{W}$. 
		
		\textbf{Case 1.} Suppose first that $\varepsilon_M<\varepsilon_X$. If $M\in X$, then $M\in\mathcal{W}$, and hence $X\cap M\in\mathcal{W}$ because $\mathcal{W}$ is an $(\mathcal{S},\mathcal{L})$-symmetric system. Suppose that $M\notin X$. Then there is some $Z\in\mathcal{M}(\varepsilon_X)$ such that $M\in Z$. Let $N:=\Psi_{Z,X}(M)$, which is a member of $\mathcal{M}\cap X\subseteq\mathcal{W}$ by symmetry, and note that as $Y\in Z\cap X$, then $Y=\Psi_{Z,X}(Y)\in N$ by clause (B). Hence, $Y\cap N\in\mathcal{W}$. But then we can conclude that $Y\cap M\in\mathcal{W}$, as $Y\cap N=Y\cap M$ by Proposition \ref{prop21}. This case is depicted in Figure \ref{fig1}.
				
		\begin{figure}[ht]
			\centering
				\begin{tikzpicture} 
					\tikzset{thin/.style = {line width= 1pt}} 
					
					\coordinate[] (x) at (2,0){};
					\coordinate[] (y) at (0,0.5){};
					
					\coordinate[] (YMYN) at (0,0){};
					\node[right=3pt] at (YMYN){$Y\cap M=\color{blue}{Y\cap N}$};
					\draw[blue, thin] (-0.14,0)--(0.14,0);
					
					\coordinate[] (Y) at (0,1){};
					\node[right=3pt] at (Y){{\color{blue}$Y$}};
					\draw[thin] (YMYN)--(Y);
					\draw[blue, thin] (-0.14,1)--(0.14,1);
					
					\coordinate[] (M) at (-1,2){};
					\node[left=3pt] at (M){$M$};
					\draw[thin] (Y)--(M);
					\draw[thin] (-1.14,2)--(-0.86,2);
					
					\coordinate[] (Z) at (-1,3){};
					\node[left=3pt] at (Z){$Z$};
					\draw[thin, -|] (M)--(Z);
					
					\coordinate[] (N) at (1,2){};
					\node[right=3pt] at (N){{\color{blue}$N=\Psi_{Z,X}(M)$}};
					\draw[thin] (Y)--(N);
					\draw[blue, thin] (0.86,2)--(1.14,2);
					
					\coordinate[] (X) at (1,3){};
					\node[right=3pt] at (X){$X$};
					\draw[thin, -|] (N)--(X);
					
				\end{tikzpicture}
				\caption{The final configuration of Case 1. The models $Y\cap N, Y, N=\Psi_{Z,X}(M)$ are members of $\mathcal{W}$. The rest of the models are members of $\mathcal{M}$.}
				\label{fig1}
			\end{figure}
		
		\textbf{Case 2.} Suppose now that $\varepsilon_X<\varepsilon_M$. We will argue by induction that this case can be reduced to the first one. We will need the following claim.
		
		\begin{claim}\label{claim-pure-amalgamation-1}
			If $N$ is a small model in $\mathcal{M}$ such that $\varepsilon_X<\varepsilon_N$ and $Y\in N$, then we can find some $Z\in\mathcal{M}\cap\mathcal{L}$ such that $\varepsilon_X\leq\varepsilon_Z$ and $Y\in Z\in N$.
		\end{claim}
		\begin{claimproof}
			Let $\varepsilon$ be the maximum of $\dom(\mathcal{M}\cap\mathcal{L})\cap\varepsilon_N$, which is greater than or equal to $\varepsilon_X$. If $\varepsilon=\varepsilon_X$, as $X,N\in\mathcal{M}$, we can use the shoulder axiom to find $N'\in\mathcal{M}(\varepsilon_N)$ such that $X\in N'$ (note that by the maximality of $\varepsilon$ there are no large models of $\omega_2$-height strictly between $\varepsilon_X$ and $\varepsilon_N$, so we can apply Proposition \ref{prop18}). Hence, if we let $Z:=\Psi_{N'[\omega_1],N[\omega_1]}(X)$, which is a member of $\mathcal{M}$ by symmetry, we have $$Y=\Psi_{N'[\omega_1],N[\omega_1]}(Y)\in Z\in N,$$ by clause (B) for $\mathcal{M}$, as $Y\in N'[\omega_1]\cap N$. Suppose now that $\varepsilon>\varepsilon_X$. Let $Z',N'\in\mathcal{M}$ be such that $\varepsilon_{Z'}=\varepsilon$, $\varepsilon_{N'}=\varepsilon_N$, and $X\in Z'\in N'[\omega_1]$, by two applications of the shoulder axiom for $\mathcal{M}$, and note that $Z'\in N'$ by the maximality of $\varepsilon$ and Proposition \ref{prop18}. Let $Z:=\Psi_{N'[\omega_1],N[\omega_1]}(Z')$, which is an element of $\mathcal{M}$ by symmetry, and note that as $Y\in N'[\omega_1]\cap N$, then $Y=\Psi_{N'[\omega_1],N[\omega_1]}(Y)\in Z\in N$ by clause (B) (see Figure \ref{fig2}).
			
			\begin{figure}[ht]
				\centering
					\begin{tikzpicture} 
						\tikzset{thin/.style = {line width= 1pt}} 
						
						\coordinate[] (x) at (2,0){};
						\coordinate[] (y) at (0,0.5){};
						
						\node at (4.7,0) {};
						
						\coordinate[] (Y) at (0,0){};
						\node[right=3pt] at (Y){$Y$};
						\draw[thin] (-0.14,0)--(0.14,0);
						
						\coordinate[] (Z) at (-1,2){};
						\node[left=3pt] at (Z){$\Psi_{N'[\omega_1],N[\omega_1]}(Z')=Z$};
						\draw[thin] (Y)--(Z);
						\draw[thin] (-1.14,2)--(-0.86,2);
						
						\coordinate[] (N) at (-1,3){};
						\node[left=3pt] at (N){$N$};
						\draw[thin, -|] (Z)--(N);
						
						\coordinate[] (X) at (1,1){};
						\node[right=3pt] at (X){$X$};
						\draw[thin] (Y)--(X);
						\draw[thin] (0.86,1)--(1.14,1);
						
						\coordinate[] (Z') at (1,2){};
						\node[right=3pt] at (Z'){$Z'$};
						\draw[thin] (X)--(Z');
						\draw[thin] (0.86,2)--(1.14,2);
						
						\coordinate[] (N') at (1,3){};
						\node[right=3pt] at (N'){$N'$};
						\draw[thin, -|] (Z')--(N');
					\end{tikzpicture}
					\caption{The final configuration in the proof of Claim \ref{claim-pure-amalgamation-1}}
					\label{fig2}
				\end{figure}
		\end{claimproof}
		
		Let $Z\in\mathcal{M}\cap\mathcal{L}$ such that $\varepsilon_X\leq\varepsilon_Z$ and $Y\in Z\in M$, given by the last claim. Then $Y\in Z\cap M\in\mathcal{M}$. If $\varepsilon_X<\varepsilon_{Z\cap M}$, note that we can use the claim again to find some $Z'\in\mathcal{M}\cap\mathcal{L}$ such that $\varepsilon_X\leq\varepsilon_{Z'}$ and $Y\in Z'\in Z\cap M$. Therefore, $$Y\in Z'\cap (Z\cap M)=Z'\cap M\in\mathcal{M}.$$ Hence, by successive applications of the last claim we can find a large model $Z_0\in\mathcal{M}$ such that
		\begin{itemize}
			\item $\varepsilon_X\leq\varepsilon_{Z_0}$,
			
			\item $Y\in Z_0\in M$, and
			
			\item $\varepsilon_{Z_0\cap M}<\varepsilon_X$.
		\end{itemize}   
		Therefore, by the first case, $Y\cap(Z_0\cap M)=Y\cap M\in\mathcal{W}$.
	\end{proof}
	
	\begin{proposition}\label{prop-pure-amalgamation-2}
		Let $\mathcal{M}$ be an $(\mathcal{S},\mathcal{L})$-symmetric system and let $M$ be a small model in $\mathcal{M}$. Let $\mathcal{V}$ be another $(\mathcal{S},\mathcal{L})$-symmetric system such that $\mathcal{M}\cap M[\omega_1]\subseteq\mathcal{V}\subseteq M[\omega_1]$ and $\mathcal{V}(\max\dom(\mathcal{V}))\subseteq M$. Define
		\[
		\varepsilon^-:=\max\dom(\mathcal{M}(<\varepsilon_M)\cap\mathcal{L})
		\]
		and 
		\[
		\varepsilon^+:=\min\dom(\mathcal{M}(>\varepsilon_M)\cap\mathcal{L}),
		\]
		if they exist, otherwise let $\varepsilon^-$ be any ordinal smaller than $\min\dom(\mathcal{M})$ and let $\varepsilon^+$ be any ordinal greater than $\max\dom(\mathcal{M})$. If $$\mathcal{V}(>\varepsilon^-)\land\mathcal{M}\restr[\varepsilon_M,\varepsilon^+)\subseteq\mathcal{V},$$ then $\mathcal{V}\land\mathcal{M}\subseteq\mathcal{V}$.
	\end{proposition}
	\begin{proof}
		Let $N\in\mathcal{M}\cap\mathcal{S}$ and $X\in\mathcal{V}\cap\mathcal{L}$ such that $X\in N$. We have to show that $X\cap N\in\mathcal{V}$. We divide the proof into three cases.
		
		\textbf{Case 1.} Suppose that $\varepsilon_N<\varepsilon_M$. Use the shoulder axiom to find a model $M'\in\mathcal{M}(\varepsilon_M)$ such that $N\in M'[\omega_1]$. Define $N_0:=\Psi_{M'[\omega_1],M[\omega_1]}(N)$, and note that $N_0\in\mathcal{M}\cap M[\omega_1]\subseteq\mathcal{V}$. Moreover, since $X\in N\subseteq M'[\omega_1]$ and $X\in\mathcal{V}\subseteq M[\omega_1]$, the isomorphism $\Psi_{M'[\omega_1],M[\omega_1]}$ fixes $X$, and hence $X\in N_0$. Since both $X$ and $N_0$ are members of $\mathcal{V}$, their intersection $X\cap N_0$ is also a member of $\mathcal{V}$. Therefore, by Proposition \ref{prop21}, we can conclude that $X\cap N=X\cap N_0\in\mathcal{V}$. This case is pictured in Figure \ref{fig3}.
		
		\begin{figure}[ht]
			\centering
				\begin{tikzpicture} 
					\tikzset{thin/.style = {line width= 1pt}} 
					
					\coordinate[] (x) at (2,0){};
					\coordinate[] (y) at (0,0.5){};
					
					\coordinate[] (XNXN0) at (0,0){};
					\node[right=3pt] at (XNXN0){$X\cap N=\color{blue}{X\cap N_0}$};
					\draw[blue, thin] (-0.14,0)--(0.14,0);
					
					\coordinate[] (X) at (0,1){};
					\node[right=3pt] at (X){{\color{blue}$X$}};
					\draw[thin] (YMYN)--(X);
					\draw[blue, thin] (-0.14,1)--(0.14,1);
					
					\coordinate[] (N) at (-1,2){};
					\node[left=3pt] at (N){$N$};
					\draw[thin] (X)--(N);
					\draw[thin] (-1.14,2)--(-0.86,2);
					
					\coordinate[] (M') at (-1,3){};
					\node[left=3pt] at (M'){$M'$};
					\draw[thin, -|] (N)--(M');
					
					\coordinate[] (N0) at (1,2){};
					\node[right=3pt] at (N0){{\color{blue}$N_0=\Psi_{M'[\omega_1],M[\omega_1]}(N)$}};
					\draw[thin] (X)--(N0);
					\draw[blue, thin] (0.86,2)--(1.14,2);
					
					\coordinate[] (M) at (1,3){};
					\node[right=3pt] at (M){$M$};
					\draw[thin, -|] (N0)--(M);
					
				\end{tikzpicture}
				\caption{The final configuration of Case 1. The models $X\cap N_0,X,N_0$ are members of $\mathcal{V}$. The rest of the models are members of $\mathcal{M}$.}
				\label{fig3}
			\end{figure}
		
		\textbf{Case 2.} Suppose that $\varepsilon_M\leq\varepsilon_N<\varepsilon^+$. If $\varepsilon_X>\varepsilon^-$, the conclusion follows from the hypothesis
		\[
		\mathcal{V}(>\varepsilon^-)\land\mathcal{M}\restr[\varepsilon_M,\varepsilon^+)\subseteq\mathcal{V}.
		\]
		So assume that $\varepsilon_X\leq\varepsilon^-$. Apply the shoulder axiom to find $N_0\in\mathcal{M}(\varepsilon_N)$ such that $M\in N_0[\omega_1]$, and note that $M\in N_0$, by the minimality of $\varepsilon^+$ and Proposition \ref{prop18}. Moreover, note that since $X\in N\cap N_0[\omega_1]$, the isomorphism $\Psi_{N[\omega_1],N_0[\omega_1]}$ fixes $X$, and hence $X\in N_0$. Now, let $Y\in\mathcal{V}\cap\mathcal{L}$ be such that $\varepsilon_Y=\max\dom(\mathcal{V}\cap\mathcal{L})$ and $X\in Y$, given by the shoulder axiom for $\mathcal{V}$, and note that, as $\mathcal{M}\cap M[\omega_1]\subseteq\mathcal{V}$, then $\varepsilon_Y\geq\varepsilon^-$. Hence, since there are no large models between $Y$ and $M$ (because of the maximality of $\varepsilon_Y$), and $\mathcal{V}\subseteq M[\omega_1]$ and $\mathcal{V}(\max\dom(\mathcal{V}))\subseteq M$, by Proposition \ref{prop18} we have $Y\in M\subseteq N_0$. Therefore, by the hypothesis 
		\[
		\mathcal{V}(>\varepsilon^-)\land\mathcal{M}\restr[\varepsilon_M,\varepsilon^+)\subseteq\mathcal{V},
		\]
		we have $Y\cap N_0\in\mathcal{V}$. Lastly, since $X\in Y\cap N_0$, and both $X$ and $Y\cap N_0$ are members of $\mathcal{V}$, $X\cap(Y\cap N_0)=X\cap N_0$ is also a member of $\mathcal{V}$. Therefore, by Proposition \ref{prop21}, we can conclude that $X\cap N=X\cap N_0\in\mathcal{V}$. This case is pictured in Figure \ref{fig4}
		
		\begin{figure}[ht]
			\centering
				\begin{tikzpicture} 
					\tikzset{thin/.style = {line width= 1pt}} 
					
					\coordinate[] (x) at (2,0){};
					\coordinate[] (y) at (0,0.5){};
					
					\coordinate[] (e+) at (4,4.5){};
					\node[right=2pt] at (e+){$\varepsilon^+$};
					\draw[thin] (3.86,4.42)--(4.14,4.42);
					\draw[dotted, thin] (-2.5,4.42)--(4.14,4.42);
					
					\coordinate[] (e-) at (4,1.5){};
					\node[right=2pt] at (e-){$\varepsilon^-$};
					\draw[thin] (3.86,1.42)--(4.14,1.42);
					\draw[dotted, thin] (-2.5,1.42)--(4.14,1.42);
					
					\coordinate[] (XNXN0) at (0,0){};
					\node[right=3pt] at (XNXN0){$X\cap N=\color{blue}{X\cap N_0}=X\cap(Y\cap N_0)$};
					\draw[blue, thin] (-0.14,0)--(0.14,0);
					
					\coordinate[] (X) at (0,1){};
					\node[right=3pt] at (X){{\color{blue}$X$}};
					\draw[thin] (YMYN)--(X);
					\draw[blue, thin] (-0.14,1)--(0.14,1);
					
					\coordinate[] (N) at (-1,4){};
					\node[left=3pt] at (N){$N$};
					\draw[thin] (X)--(N);
					\draw[thin] (-1.14,4)--(-0.86,4);
					
					\coordinate[] (YN0) at (1,2){};
					\node[right=3pt] at (YN0){{{\color{blue}$Y\cap N_0$}}};
					\draw[thin] (X)--(YN0);
					\draw[blue, thin] (0.86,2)--(1.14,2);
					
					\coordinate[] (Y) at (1,2.5){};
					\node[right=3pt] at (Y){{\color{blue}{$Y$}}};
					\draw[thin] (YN0)--(Y);
					\draw[blue, thin] (0.86,2.5)--(1.14,2.5);
										
					\coordinate[] (M) at (1,3){};
					\node[right=5pt] at (M){{$M$}};
					\draw[thin] (Y)--(M);
					\draw[blue, thin] (0.75,3)--(1.25,3);
					
					\coordinate[] (N0) at (1,4){};
					\node[right=3pt] at (N0){$N_0$};
					\draw[thin, -|] (M)--(N0);
					
				\end{tikzpicture}
				\caption{The final configuration of Case 2. The models $X\cap N_0,X,Y\cap N_0,Y$ are members of $\mathcal{V}$. The rest of the models are members of $\mathcal{M}$.}
				\label{fig4}
			\end{figure}
				
		\textbf{Case 3.} Suppose that $\varepsilon^+<\varepsilon_N$. This case requires an inductive argument similar to the one from the proof of the Case 2 of Proposition \ref{prop-pure-amalgamation-1}. Let us first show the following claim.
		
		\begin{claim}\label{claim-pure-amalgamation-2}
			If $N_0$ is a small model in $\mathcal{M}$ such that $\varepsilon^+<\varepsilon_{N_0}$ and $X\in N_0$, then we can find some $Z\in\mathcal{M}\cap\mathcal{L}$ such that $\varepsilon^+\leq\varepsilon_Z$ and $X\in Z\in N_0$.
		\end{claim}
		\begin{claimproof}
			Let $\varepsilon$ be the maximum of $\dom(\mathcal{M}\cap\mathcal{L})\cap\varepsilon_{N_0}$, and note that $\varepsilon\geq\varepsilon^+$. Let $Z',N_0'\in\mathcal{M}$ such that
			\begin{itemize}
				\item $\varepsilon_{Z'}=\varepsilon$,
				\item $\varepsilon_{N_0'}=\varepsilon_{N_0}$, and
				\item $M\in Z'\in N_0'[\omega_1]$,
			\end{itemize} 
			by two applications of the shoulder axiom for $\mathcal{M}$, and note that $Z'\in N_0'$ by the maximality of $\varepsilon$ and Proposition \ref{prop18}. Let $Z:=\Psi_{N_0'[\omega_1],N_0[\omega_1]}(Z')$, which is an element of $\mathcal{M}\cap N_0$ because $N_0',N_0,Z'\in\mathcal{M}$. Since $X\in M[\omega_1]\subseteq N_0'[\omega_1]$, and $X\in N_0$ by assumption, $X$ is fixed by the isomorphism $\Psi_{N_0'[\omega_1],N_0[\omega_1]}$. Hence, $X=\Psi_{N_0'[\omega_1],N_0[\omega_1]}(X)\in Z\in N_0$.
		\end{claimproof}
		
		Let $Z\in\mathcal{M}\cap\mathcal{L}$ such that $\varepsilon^+\leq\varepsilon_Z$ and $X\in Z\in N$, given by the last claim, and note that $X\in Z\cap N\in\mathcal{M}$. If $\varepsilon^+<\varepsilon_{Z\cap N}$, we can use the claim again to find some $Z'\in\mathcal{M}\cap\mathcal{L}$ such that $\varepsilon^+\leq\varepsilon_{Z'}$ and $X\in Z'\in Z\cap N$, and hence $$X\in Z'\cap(Z\cap N)=Z'\cap N\in\mathcal{M}.$$ Therefore, by successive applications of the last claim we can find a large model $Z_0\in\mathcal{M}$ such that
		\begin{itemize}
			\item $\varepsilon^+\leq\varepsilon_{Z_0}$,
			\item $X\in Z_0\in N$, and
			\item $\varepsilon_{Z_0\cap N}<\varepsilon^+$.
		\end{itemize} 
		So, by the first two cases, $X\cap (Z_0\cap N)=X\cap N\in\mathcal{V}$.
	\end{proof}
	
	\begin{proposition}\label{prop-pure-amalgamation-3}
		Let $\mathcal{M}$ be an $(\mathcal{S},\mathcal{L})$-symmetric system, and let $Q\in\mathcal{M}$. Let $\mathcal{W}$ be another $(\mathcal{S},\mathcal{L})$-symmetric system such that
		\begin{enumerate}
			\item $\mathcal{M}\cap Q[\omega_1]\subseteq\mathcal{W}\subseteq Q[\omega_1]$,
			
			\item $\mathcal{W}(\max\dom(\mathcal{W}))\subseteq Q$, and
			
			\item $\mathcal{W}\land\mathcal{M}\subseteq\mathcal{W}$.
		\end{enumerate}   
		Then,
		\[
		\mathcal{V}:=\mathcal{M}(\geq\varepsilon_Q)\cup\{\Psi_{Q[\omega_1],P[\omega_1]}(W):W\in\mathcal{W},P\in\mathcal{M}(\varepsilon_Q)\}
		\]
		is an $(\mathcal{S},\mathcal{L})$-symmetric system such that $\mathcal{M}\cup\mathcal{W}\subseteq\mathcal{V}$.
	\end{proposition}
	\begin{proof}
		First, note that 
		\[
		\mathcal{V}(<\varepsilon_Q)=\{\Psi_{Q[\omega_1],P[\omega_1]}(W):W\in\mathcal{W},P\in\mathcal{M}(\varepsilon_Q)\}
		\]
		is an $(\mathcal{S},\mathcal{L})$-symmetric system by Proposition \ref{isocopy} and Lemma \ref{pureamalgamation2}, which extends $\mathcal{M}(<\varepsilon_Q)$ and $\mathcal{W}$. It is straightforward to check that $\mathcal{V}$ satisfies clauses (A)-(C) of the definition of $(\mathcal{S},\mathcal{L})$-symmetric system. Hence, we only need to check that $\mathcal{V}$ satisfies clauses (D) and (E).
		
		Let us first check clause (D). Let $V_0,V_1,V_1'\in\mathcal{V}$ be such that $V_0\in V_1$ and $\varepsilon_{V_1}=\varepsilon_{V_1'}$. We have to show that $\Psi_{V_1[\omega_1],V_1'[\omega_1]}(V_0)\in\mathcal{V}$. Note that if $\varepsilon_{V_1}=\varepsilon_{V_1'}<\varepsilon_Q$, the conclusion follows from the fact that $\mathcal{V}(<\varepsilon_Q)$ is an $(\mathcal{S},\mathcal{L})$-symmetric system, and if $\varepsilon_Q\leq\varepsilon_{V_0}$, the conclusion follows from $\mathcal{V}(\geq\varepsilon_Q)=\mathcal{M}(\geq\varepsilon_Q)$. Therefore, we may assume that $\varepsilon_{V_0}<\varepsilon_Q\leq\varepsilon_{V_1}=\varepsilon_{V_1'}$. Since $V_0\in\mathcal{V}(<\varepsilon_Q)$, there must be some $W_0\in\mathcal{W}$ and some $P\in\mathcal{M}(\varepsilon_Q)$ such that $V_0=\Psi_{Q[\omega_1],P[\omega_1]}(W_0)$. Let us divide the proof into two cases.
		
		\textbf{Case 1.} Suppose first that $\varepsilon_Q=\varepsilon_{V_1}=\varepsilon_{V_1'}$. Since $V_0\in P[\omega_1]\cap V_1$, the model $V_0$ is fixed by the isomorphism $\Psi_{V_1[\omega_1],P[\omega_1]}$. Therefore,
		\begin{align*}
			\Psi_{V_1[\omega_1],V_1'[\omega_1]}(V_0)&=\Psi_{Q[\omega_1],V_1'[\omega_1]}\circ\Psi_{P[\omega_1],Q[\omega_1]}\circ\Psi_{V_1[\omega_1],P[\omega_1]}(V_0)\\
			&=\Psi_{Q[\omega_1],V_1'[\omega_1]}\circ\Psi_{P[\omega_1],Q[\omega_1]}(V_0)\\
			&=\Psi_{Q[\omega_1],V_1'[\omega_1]}(W_0),
		\end{align*}
		which is an element of $\mathcal{V}$ by definition. 
		
		\textbf{Case 2.} Suppose now that $\varepsilon_Q<\varepsilon_{V_1}=\varepsilon_{V_1'}$. Since $P,V_1\in\mathcal{M}$, there must be some $V_2\in\mathcal{M}(\varepsilon_{V_1})$ such that $P\in V_2[\omega_1]$ by the shoulder axiom. Let $P_1:=\Psi_{V_2[\omega_1],V_1[\omega_1]}(P)$, which is an element of $\mathcal{M}$ by symmetry, and note that as $V_0\in V_2[\omega_1]\cap V_1$, then $V_0=\Psi_{V_2[\omega_1],V_1[\omega_1]}(V_0)\in P_1[\omega_1]$. Let $P_1':=\Psi_{V_1[\omega_1],V_1'[\omega_1]}(P_1)$, which is an element of $\mathcal{M}$ because $V_1,V_1',P_1\in\mathcal{M}$, and note that $\Psi_{V_1[\omega_1],V_1'[\omega_1]}(V_0)=\Psi_{P_1[\omega_1],P_1'[\omega_1]}(V_0)$ (see Figure \ref{fig5}). Hence, we will be done if we show that $\Psi_{P_1[\omega_1],P_1'[\omega_1]}(V_0)\in\mathcal{V}$. However, note that this is not immediate from Case 1, since a priori $V_0$ does not need to be an element of $P_1$. 
		
		\begin{figure}[ht]
			\centering
				\begin{tikzpicture} 
					\tikzset{thin/.style = {line width= 1pt}} 
					
					\coordinate[] (x) at (2,0){};
					\coordinate[] (y) at (0,0.5){};
					
					\node at (4.7,0) {};
					
					\coordinate[] (V0) at (0,0){};
					\node[right=3pt] at (V0){$V_0$};
					\draw[thin] (-0.14,0)--(0.14,0);
					
					\coordinate[] (P) at (-1.5,1){};
					\node[left=3pt] at (P){$P$};
					\draw[thin,dotted] (V0)--(P);
					\draw[thin] (-1.64,1)--(-1.36,1);
					
					\coordinate[] (V2) at (-1.5,2){};
					\node[left=3pt] at (V2){$V_2$};
					\draw[thin, dotted, -|] (P)--(V2);
					
					\coordinate[] (P1) at (1.5,1){};
					\node[right=3pt] at (P1){$P_1=\Psi_{V_2[\omega_1],V_1[\omega_1]}(P)$};
					\draw[thin,dotted] (V0)--(P1);
					\draw[thin] (1.36,1)--(1.64,1);
										
					\coordinate[] (V1) at (1.5,2){};
					\node[above=3pt] at (V1){$V_1$};
					\draw[thin,dotted,-|] (P1)--(V1);
					\draw[thin,->] (-1.3,2) to node [above] {$\Psi_{V_2[\omega_1],V_1[\omega_1]}$} (1.3,2);
					
					\coordinate[] (V0') at (6,0){};
					\node[right=3pt] at (V0'){$\Psi_{V_1[\omega_1],V_1'[\omega_1]}(V_0)$};
					\draw[thin] (5.86,0)--(6.14,0);
					
					\coordinate[] (P1') at (6,1){};
					\node[right=3pt] at (P1'){$P_1'=\Psi_{V_1[\omega_1],V_1'[\omega_1]}(P_1)$};
					\draw[thin,dotted,-|] (V0')--(P1');
					
					\coordinate[] (V1') at (6,2){};
					\node[right=3pt] at (V1'){$V_1'$};
					\draw[thin,dotted,-|] (P1')--(V1');
					\draw[thin,->] (1.8,2) to node [above] {$\Psi_{V_1[\omega_1],V_1'[\omega_1]}$} (5.7,2);
					
				\end{tikzpicture}
				\caption{The dotted lines denote that the bottom models are members of the hulls of the top models.}
				\label{fig5}
			\end{figure}
		
		If there is no large model $X\in\mathcal{V}$ such that $\varepsilon_{V_0}<\varepsilon_X<\varepsilon_Q$, since $\mathcal{W}(\max\dom(\mathcal{W}))\subseteq Q$ by hypothesis, then $V_0\in P_1$, and hence we get the conclusion from the first case. Suppose now that there is some large model $X\in\mathcal{V}$ such that $\varepsilon_{V_0}<\varepsilon_X<\varepsilon_Q$. Since $\mathcal{V}(<\varepsilon_M)$ is an $(\mathcal{S},\mathcal{L})$-symmetric system, we can find such a model $X_0$ of maximal $\omega_2$-height and such that $V_0\in X_0\in P_1[\omega_1]$. Note that, by the maximality of $\varepsilon_{X_0}$, and since $\mathcal{W}(\max\dom(\mathcal{W}))\subseteq Q$, in fact, $X_0\in P_1$. Therefore, by the first case, $X_0':=\Psi_{P_1[\omega_1],P_1'[\omega_1]}(X_0)\in\mathcal{V}$. Note that $\Psi_{P_1[\omega_1],P_1'[\omega_1]}(V_0)=\Psi_{X_0,X_0'}(V_0)$. Hence, the conclusion follows from $V_0,X_0,X_0'\in\mathcal{V}(<\varepsilon_Q)$ and the fact that $\mathcal{V}(<\varepsilon_Q)$ is an $(\mathcal{S},\mathcal{L})$-symmetric system.
		
		Let us now show that $\mathcal{V}$ satisfies clause (E). Let $X\in\mathcal{V}\cap\mathcal{L}$ and $M\in\mathcal{V}\cap\mathcal{S}$ such that $X\in M$. We need to check that $X\cap M\in\mathcal{V}$. Once again, note that if $\varepsilon_M<\varepsilon_Q$, the conclusion follows from the fact that $\mathcal{V}(<\varepsilon_Q)$ is an $(\mathcal{S},\mathcal{L})$-symmetric system, and if $\varepsilon_Q\leq\varepsilon_X$, the conclusion follows from $\mathcal{V}(\geq\varepsilon_Q)=\mathcal{M}(\geq\varepsilon_Q)$. So, we may assume that $\varepsilon_X<\varepsilon_Q\leq\varepsilon_M$. Since $X\in\mathcal{V}(<\varepsilon_Q)$, there are $W\in\mathcal{W}$ and $P\in\mathcal{M}(\varepsilon_Q)$ such that $X=\Psi_{Q[\omega_1],P[\omega_1]}(W)$. We divide the proof in two cases.
		
				
		\textbf{Case 1.} There exists no large model $Z\in\mathcal{M}(\geq\varepsilon_Q)\cap M$ such that $X\in Z$. In this case, there exists no large model $Z\in\mathcal{M}$ such that $\varepsilon_Q\leq\varepsilon_Z<\varepsilon_M$. By the shoulder axiom and Proposition \ref{prop18}, we can find two models $N,M'\in\mathcal{M}(\varepsilon_M)$ such that $P\in M'$ and $Q\in N$. Since $X\in M\cap M'[\omega_1]$, and $M$ and $M'$ are $\omega_1$-isomorphic, $X\in M'$. Let $P_0=\Psi_{M'[\omega_1],N[\omega_1]}(P)$ and $X_0=\Psi_{M'[\omega_1],N[\omega_1]}(X)$. Since $X\in M'$, and $M'$ and $N$ are $\omega_1$-isomorphic, we have $X_0\in N$. Since $P_0,Q\in N$, the isomorphism $\Psi_{P_0[\omega_1],Q[\omega_1]}$ is also a member of $N$. Therefore, $W=\Psi_{P_0[\omega_1],Q[\omega_1]}(X_0)\in N$. Now, since $W\cap N$ is a member of $\mathcal{W}\land\mathcal{M}\subseteq\mathcal{W}$, by an application of Proposition \ref{prop10},
		\[
			X\cap M'=\Psi_{W,X}(W\cap N)=\Psi_{Q[\omega_1],P[\omega_1]}(W\cap N)\in\mathcal{V}(<\varepsilon_Q).
		\]		
		Hence, we can conclude that $X\cap M=X\cap M'\in \mathcal{V}(<\varepsilon_Q)$ by Proposition \ref{prop21}.
		
		\textbf{Case 2.} There exists a large model $Z\in\mathcal{M}(\geq\varepsilon_Q)\cap M$ such that $X\in Z$. Let us assume that $Z$ is of minimal height among these models. We have $Z\cap M\in\mathcal{M}$.
		
		\textbf{Subcase 2.1.} $\varepsilon_Q<\varepsilon_{Z\cap M}$. In this case there exists no large model $Y\in\mathcal{M}$ with $\varepsilon_Q\leq\varepsilon_Y<\varepsilon_{Z\cap M}$. Now we argue as in Case 1, with $Z\cap M$ in place of $M$, to conclude that $X\cap(Z\cap M)\in\mathcal{V}(<\varepsilon_Q)$. Since $X\in Z\cap M$,
		\[
		X\cap M=X\cap(Z\cap M)\in\mathcal{V}(<\varepsilon_Q).
		\]
		
		\textbf{Subcase 2.2.} $\varepsilon_Q=\varepsilon_{Z\cap M}$. Since $W\cap Q$ is a member of $\mathcal{W}\land\mathcal{M}\subseteq\mathcal{W}$, by an application of Proposition \ref{prop10},
		\[
		X\cap P=\Psi_{W,X}(W\cap Q)=\Psi_{Q[\omega_1],P[\omega_1]}(W\cap Q)\in\mathcal{V}(<\varepsilon_Q).
		\]
		Hence, we can conclude that $X\cap M=X\cap (Z\cap M)=X\cap P\in \mathcal{V}(<\varepsilon_Q)$ by Proposition \ref{prop21}.
		
		\textbf{Subcase 2.3.} $\varepsilon_Q>\varepsilon_{Z\cap M}$. Since $Z\cap M\in\mathcal{M}(<\varepsilon_Q)\subseteq\mathcal{V}(<\varepsilon_Q)$, both models $X$ and $Z\cap M$ are members of $\mathcal{V}(<\varepsilon_Q)$. Hence, since $\mathcal{V}(<\varepsilon_Q)$ is an $(\mathcal{S},\mathcal{L})$-symmetric system, we have $X\cap (Z\cap M)\in\mathcal{V}(<\varepsilon_Q)$. Thus, we can conclude that $X\cap M=X\cap (Z\cap M)\in\mathcal{V}(<\varepsilon_Q)$.
	\end{proof}

	\begin{lemma}\label{lemma-amal-L}
		Let $\mathcal{M}$ be an $(\mathcal{S},\mathcal{L})$-symmetric system, and let $X\in\mathcal{M}\cap\mathcal{L}$. Let $\mathcal{W}$ be another $(\mathcal{S},\mathcal{L})$-symmetric system such that $\mathcal{M}\cap X\subseteq\mathcal{W}\subseteq X$. Then there is an $(\mathcal{S},\mathcal{L})$-symmetric system $\mathcal{V}$ such that $\mathcal{M}\cup\mathcal{W}\subseteq\mathcal{V}$. 
	\end{lemma}
	\begin{proof}
		By Proposition \ref{prop-pure-amalgamation-1} we have $\mathcal{W}\land\mathcal{M}\subseteq\mathcal{W}$, and since $X$ is a large model, $\mathcal{W}(\max\dom(\mathcal{W}))\subseteq X$. Therefore, by Proposition \ref{prop-pure-amalgamation-3},
		\[
		\mathcal{V}:=\mathcal{M}(\geq\varepsilon_X)\cup\{\Psi_{X,Y}(W):W\in\mathcal{W},Y\in\mathcal{M}(\varepsilon_X)\}
		\]
		is an $(\mathcal{S},\mathcal{L})$-symmetric system such that $\mathcal{M}\cup\mathcal{W}\subseteq\mathcal{V}$.
	\end{proof}

	\begin{lemma}\label{lemma-amal-S}
		Let $\mathcal{M}$ be an $(\mathcal{S},\mathcal{L})$-symmetric system, and let $M\in\mathcal{M}\cap\mathcal{S}$. Let $\mathcal{W}$ be another $(\mathcal{S},\mathcal{L})$-symmetric system such that $\mathcal{M}\cap M\subseteq\mathcal{W}\subseteq M$. Then there is an $(\mathcal{S},\mathcal{L})$-symmetric system $\mathcal{U}$ such that $\mathcal{M}\cup\mathcal{W}\subseteq\mathcal{U}$.
	\end{lemma}
	\begin{proof}
		Fix a maximal $\in$-chain $X_0\in\dots\in X_{n-1}$ of large models in $\mathcal{M}\cap M$ and let $X_n$ be a large model containing $\mathcal{M}\cup\mathcal{W}$.\footnote{The only use of the model $X_n$ is in making the notation more homogeneous. There will be no point in the proof in which we make explicit use of the properties of $X_n$.} For all $i\leq n$, we denote
		\begin{itemize}
			\item $\mathcal{M}_i=\mathcal{M}\cap X_i$,
			\item $\mathcal{W}_i=\mathcal{W}\cap X_i=\mathcal{W}\cap(X_i\cap M)$ (note that $\mathcal{W}_i\subseteq X_i\cap M$), and
			\item $M_i=X_i\cap M$.
		\end{itemize}    
		Observe that $\mathcal{M}_n=\mathcal{M}$, $\mathcal{W}_n=\mathcal{W}$ and $M_n=M$. Moreover, note that $\mathcal{M}_i$ and $\mathcal{W}_i$ are $(\mathcal{S},\mathcal{L})$-symmetric systems according to Lemma \ref{amal-rest}, and that $\mathcal{M}_i\cap M_i\subseteq\mathcal{W}_i\subseteq M_i$. We will build, by induction on $i\leq n$, an $\subseteq$-increasing sequence of $(\mathcal{S},\mathcal{L})$-symmetric systems $\langle \mathcal{U}_i:i\leq n\rangle$ such that $\mathcal{M}_i\cup\mathcal{W}_i\subseteq\mathcal{U}_i\subseteq X_i$, and in each step of the induction, we will build three $(\mathcal{S},\mathcal{L})$-symmetric systems $\mathcal{V}_i\subseteq\mathcal{V}_i^*\subseteq\mathcal{U}_i$. The $(\mathcal{S},\mathcal{L})$-symmetric system that witnesses the compatibility of $\mathcal{M}$ and $\mathcal{W}$ will be $\mathcal{U}_n$.
		
		Before starting with the induction, let us remark two properties of the union of $(\mathcal{S},\mathcal{L})$-symmetric systems $\mathcal{W}\cup\mathcal{M}$, which can be easily derived from propositions \ref{prop-res-gaps} and \ref{prop-res-gaps2}:
		\begin{enumerate}[label=(\roman*)]
			\item Every model $Q\in(\mathcal{W}\cup\mathcal{M})\restr[\varepsilon_{M_{i}},\varepsilon_{X_{i}})$, for some $i\leq n$, must be a member of $\mathcal{M}\setminus M$, by Proposition \ref{propgaps}.
			
			\item $Q\in\mathcal{S}\cup\mathcal{L}$ is a member of $\mathcal{W}$ if and only if either $Q\in\Psi_{X_0,Y}"(\mathcal{W}_0)$ for some $Y\in\mathcal{W}(\varepsilon_{X_0})$, or $Q\in\Psi_{X_{i+1},Y}"(\mathcal{W}_{i+1}(\geq\varepsilon_{X_i}))$ for some $i<n$ and some $Y\in\mathcal{W}(\varepsilon_{X_{i+1}})$. Note that in this case, we have
			\begin{align*}
				\Psi_{X_0,Y}"(\mathcal{W}_0)&=\Psi_{M_0[\omega_1],(Y\cap M)[\omega_1]}"(\mathcal{W}_0),\\
				\Psi_{X_{i+1},Y}"(\mathcal{W}_{i+1}(\geq\varepsilon_{X_i}))&=\Psi_{M_{i+1}[\omega_1],(Y\cap M)[\omega_1]}"(\mathcal{W}_{i+1}(\geq\varepsilon_{X_i})).
			\end{align*}
		\end{enumerate}
		These two facts will be constantly used throughout the proof, sometimes without mention.
		
		Let us start the induction now. Let $i<n$ and let $\varepsilon_i^-$ denote either $\varepsilon_{X_{i-1}}$, in case $i>0$, or an arbitrary ordinal smaller than the minimum of $\dom(\mathcal{W})$, otherwise. Note that if $i>0$, then $\varepsilon_i^-$ coincides also with the maximal $\omega_2$-height of any large model in $\dom(\mathcal{M}_i(<\varepsilon_{M_i}))$. This will be relevant when we apply Proposition \ref{prop-pure-amalgamation-2}. Suppose that $\mathcal{V}_i$ is an $(\mathcal{S},\mathcal{L})$-symmetric system with the following properties:
		\begin{enumerate}[label=(\alph*)]
			\item $\mathcal{M}_i\cap M_i[\omega_1]\subseteq\mathcal{V}_i\subseteq M_i[\omega_1]$.
			\item $\mathcal{W}_i\subseteq \mathcal{V}_i$.
			\item $\mathcal{V}_i(\geq\varepsilon_i^-)=\mathcal{W}_i(\geq\varepsilon_i^-)$.
		\end{enumerate}
		If $i=0$, we simply let $\mathcal{V}_0=\mathcal{W}_0$. This is pictured in Figure \ref{fig7}
		
		\begin{figure}[ht]
			\centering
				\begin{tikzpicture} 
					\tikzset{thin/.style = {line width= 1pt}} 
					
					\coordinate[] (x) at (2,0){};
					\coordinate[] (y) at (0,0.5){};
					
					
					
					
					\path[draw,thin,red,fill=red!10, pattern=north east lines,pattern color=red] (-1.2,0) -- (0,1.47) -- (1.2,0) -- cycle;
					
					\path[draw,thin,blue,fill=blue!10,pattern=north west lines,pattern color=blue] (-2.5,-1.5) -- (-1.23,-0.03) -- (1.23,-0.03) -- (2.5,-1.5);

					\coordinate[] (ei-1) at (-3.5,0){};
					\node[left=3pt] at (ei-1){$\varepsilon_i^-=\varepsilon_{X_{i-1}}$};
					\draw[thin,dotted] (ei-1) -- (3.5,0);
					
					\coordinate[] (Mi) at (0,1.5){};
					\node[right=3pt] at (Mi){$M_{i}=X_{i}\cap M$};
					\draw[thin, -|] (0,0)--(Mi);
					
					\coordinate[] (Xi) at (0,3){};
					\node[right=3pt] at (Xi){$X_{i}$};
					\draw[thin, -|] (Mi)--(Xi);
					
					
					\coordinate[] (Wical) at (3,0.5){};
					\node[] at (Wical){$\textcolor{red}{\mathcal{W}_i(\geq\varepsilon_i^-)}=\textcolor{blue}{\mathcal{V}_i(\geq\varepsilon_i^-)}$};
					
					\coordinate[] (Vical) at (3,-0.8){};
					\node[] at (Vical){\textcolor{blue}{$\mathcal{V}_i(<\varepsilon_i^-)$}};
					
					

				\end{tikzpicture}
				\caption{The configuration given by the induction hypothesis. The models $M_i,X_i$ are members of $\mathcal{M}$.}
				\label{fig7}
			\end{figure}
		
		We will first extend $\mathcal{V}_i$ to an $(\mathcal{S},\mathcal{L})$-symmetric system $\mathcal{V}_i^*$ that satisfies $\mathcal{V}_i^*\land\mathcal{M}_i\subseteq\mathcal{V}_i^*$. Then, by Proposition \ref{prop-pure-amalgamation-3}, 
		\[
		\mathcal{U}_i:=\mathcal{M}_i(\geq\varepsilon_{M_i})\cup\{\Psi_{M_i[\omega_1],N_i[\omega_1]}(V):V\in\mathcal{V}_i^*,N_i\in\mathcal{M}_i(\varepsilon_{M_i})\}
		\]
		will be an $(\mathcal{S},\mathcal{L})$-symmetric system such that $\mathcal{M}_i\cup\mathcal{V}_i^*\subseteq\mathcal{U}_i\subseteq X_i$. Hence, in particular, $\mathcal{W}_i\subseteq\mathcal{U}_i$. Lastly, if $i<n$, we will have to define $\mathcal{V}_{i+1}$ so that the induction can go through. Since $\mathcal{W}_{i+1}\cap X_i=\mathcal{W}_i\subseteq\mathcal{U}_i\subseteq X_i$, we will have $\mathcal{U}_i\land\mathcal{W}_{i+1}\subseteq\mathcal{U}_i$ by Proposition \ref{prop-pure-amalgamation-1}. Hence, by Proposition \ref{prop-pure-amalgamation-3},
		\[
		\mathcal{V}_{i+1}:=\mathcal{W}_{i+1}(\geq\varepsilon_{X_i})\cup\{\Psi_{X_i,Y_i}(U):U\in\mathcal{U}_i,Y_i\in\mathcal{W}_{i+1}(\varepsilon_{X_i})\}
		\]
		will be an $(\mathcal{S},\mathcal{L})$-symmetric system such that $\mathcal{W}_{i+1}\cup\mathcal{U}_i\subseteq\mathcal{V}_{i+1}$. We will finish the proof by checking that $\mathcal{V}_{i+1}$ has the properties (a)-(c) that we assumed for $\mathcal{V}_i$ at the beginning of the induction. 
		
		Let us start with the construction of $\mathcal{V}_i^*$. The argument is based on the proof of Claim 2.22 of \cite{Neeman2014Forcingwithsequencesofmodelsoftwotypes} and we will use a very similar notation, of course adapted to the current context.
		
		Let $W$ be any large model in $\mathcal{V}_i(>\varepsilon_i^-)=\mathcal{W}_i(>\varepsilon_i^-)$, and note that $W\notin\mathcal{M}$ due to the maximality of the sequence $\langle X_i:i<n\rangle$. Define
		\[
		\varepsilon_i^+:=\min\dom(\mathcal{M}_i(>\varepsilon_{M_i})\cap\mathcal{L}).
		\] 
		That is, $\varepsilon_i^+$ is the least $\omega_2$-height of any large model in $\mathcal{M}_i(>\varepsilon_{M_i})$ (it is possible that $\varepsilon_i^+=\varepsilon_{X_i}$), and note that every model $Q\in\mathcal{M}_i\restr[\varepsilon_{M_i},\varepsilon_i^+)$ is a small model. Let $E_W$ be the set of all models $N\in\mathcal{M}_i\restr[\varepsilon_{M_i},\varepsilon_i^+)$ such that $W\in N$. 
		
		Define $F_W$ as the set $\{W\cap N:N\in E_W\}$. The following properties of $F_W$ are easy to show:
		\begin{enumerate}[label=(p.\arabic*)]
			\item If $N_0,N_1\in E_W$ have the same $\omega_2$-height, then $W\cap N_0=W\cap N_1$ by Proposition \ref{prop21}.
			
			\item If $N_0',N_1'\in E_W$ are such that $N_0'\in N_1'$, then $W\cap N_0'\in W\cap N_1'$.
			
			\item As a consequence of the first two items, $F_W$ forms an $\in$-chain, and as $M_i\in E_W$, it has $W\cap M_i=W\cap (M\cap X_i)=W\cap M$ as its minimal element.
		\end{enumerate}  
		Fix an $\in $-chain $N_0\in\dots\in N_k$ of maximal length of $\mathcal{M}_i\restr[\varepsilon_{M_i},\varepsilon_i^+)$ such that $N_0=M_i$, and denote it by $E_i$. We can conclude from the three observations above that $F_W=\{W\cap N:N\in E_i\}$ for every $W\in\mathcal{W}_i(>\varepsilon_i^-)\cap\mathcal{L}$.
		
		As a consequence, we can derive the following additional properties of $F_W$:
		\begin{enumerate}
			\item[(p.4)] If $W_0,W_1\in\mathcal{W}_i(>\varepsilon_i^-)\cap\mathcal{L}$ have the same $\omega_2$-height, then $\Psi_{W_0,W_1}$ witnesses $(W_0\cap N)[\omega_1]\cong(W_1\cap N)[\omega_1]$, for every $N\in E_i$, by Proposition \ref{prop9}.
			
			\item[(p.5)] Consequently, if $W\in\mathcal{W}_i(>\varepsilon_i^-)\cap\mathcal{L}$ and $\varepsilon_0$ is the immediate predecessor of $\varepsilon_{W}$ in $\dom(\mathcal{V}_i)$, then the ordinals $\{\varepsilon_{W\cap N}:W\cap N\in F_W\}$ lie strictly between $\varepsilon_0$ and $\varepsilon_W$. 
			
			\item[(p.6)] Let $W\in\mathcal{W}_i(>\varepsilon_i^-)\cap\mathcal{L}$ and $V\in\mathcal{V}_i\cap W$ such that $\varepsilon_V$ is the immediate predecessor of $\varepsilon_W$ in $\dom(\mathcal{V}_i)$. Then, by (c) in the definition of $\mathcal{V}_i$, $V$ must be a member of $\mathcal{W}_i$. Therefore, as $\mathcal{W}\subseteq M$, $V\in W\cap M$. 
		\end{enumerate} 
		
		\begin{figure}[ht]
			\centering
				\begin{tikzpicture} 
					\tikzset{thin/.style = {line width= 1pt}} 
					
					\coordinate[] (x) at (2,0){};
					\coordinate[] (y) at (0,0.5){};
					
					
					
					\path[draw,thin,red,fill=red!2] (-4.8,-1.8) -- (0,3.47) -- (4.8,-1.8);
					
					\path[draw,thin,green,fill=green!2] (-2.4,-1.8) -- (0,0.77) -- (2.4,-1.8);

					\coordinate[] (ei-1) at (-5.2,-1.8){};
					\node[left=3pt] at (ei-1){$\varepsilon_i^-$};
					\draw[thin,dotted] (ei-1) -- (5.2,-1.8);
					
					\coordinate[] (ei+1) at (-5.2,2.4){};
					\node[left=3pt] at (ei+1){$\varepsilon_i^+$};
					\draw[thin,dotted] (ei+1) -- (5.2,2.4);

					\coordinate[] (W) at (-0.3,0.1){};
					\node[right=4pt] at (W){\textcolor{green}{$W$}};
					\draw[thin,green] (-0.50,0.1)--(-0.10,0.1);
					\coordinate[] (WNk) at (-0.3,-0.3){};
					\node[right=2pt] at (WNk){\textcolor{purple}{$W\cap N_k$}};
					\draw[thin,purple] (-0.45,-0.3)--(-0.15,-0.3);
					\draw[thin,purple] (-0.45,-0.5)--(-0.15,-0.5);
					\draw[thin,purple] (-0.45,-0.7)--(-0.15,-0.7);
					\coordinate[] (WNj) at (-0.3,-0.8){};
					\node[right=2pt] at (WNj){\textcolor{purple}{$W\cap N_j$}};
					\draw[thin,purple] (-0.45,-0.9)--(-0.15,-0.9);
					\draw[thin,purple] (-0.45,-1.1)--(-0.15,-1.1);
					\coordinate[] (WMi) at (-0.3,-1.3){};
					\node[right=2pt] at (WMi){\textcolor{purple}{$W\cap M_i=W\cap M$}};
					\draw[thin,purple] (-0.45,-1.3)--(-0.15,-1.3);

					\coordinate[] (Mi) at (0,0.8){};
					\node[right=5pt] at (Mi){$M_i=\textcolor{orange}{N_0}$};
					\draw[thin,orange] (-0.3,0.8)--(0.3,0.8);
					
					\draw[thin,orange] (-0.15,1)--(0.15,1);
					\draw[thin,orange] (-0.15,1.2)--(0.15,1.2);
					\coordinate[] (Nj) at (0,1.4){};
					\node[right=2pt] at (Nj){\textcolor{orange}{$N_j$}};
					\draw[thin,orange] (-0.15,1.4)--(0.15,1.4);
					\draw[thin,orange] (-0.15,1.6)--(0.15,1.6);
					\draw[thin,orange] (-0.15,1.8)--(0.15,1.8);
					\coordinate[] (Nk) at (0,2){};
					\node[right=2pt] at (Nk){\textcolor{orange}{$N_k$}};
					\draw[thin,orange] (-0.15,2)--(0.15,2);
										
					\coordinate[] (Xi) at (0,3.5){};
					\node[right=5pt] at (Xi){$X_{i}$};
					\draw[thin] (-0.3,3.5)--(0.3,3.5);
					
					\coordinate[] (Mical) at (-3.3,1){};
					\node[right=3pt] at (Mical){\textcolor{red}{$\mathcal{M}_i$}};
					
					\coordinate[] (Wical) at (-3.4,-0.5){};
					\node[] at (Wical){$\textcolor{blue}{\mathcal{V}_i(\geq\varepsilon_i^-)}=\textcolor{green}{\mathcal{W}_i(\geq\varepsilon_i^-)}$};

				\end{tikzpicture}
				\caption{The construction of $\mathcal{V}_i^*$. The inner triangle represents $\mathcal{V}_i(\geq\varepsilon_i^-)=\mathcal{W}_i(\geq\varepsilon_i^-)$. The outer triangle represents $\mathcal{M}_i$. The models $W\cap M_i, W\cap N_j,W\cap N_k$ are members of $F_W$. The model $W$ is a member of $\mathcal{W}_i(\geq\varepsilon_i^-)$. The models $N_0,N_j,N_k$ are members of $E_i$. The models $M_i,X_i$ are members of $\mathcal{M}$.}
				\label{fig8}
			\end{figure}
		
		Let $\mathcal{V}_i^*$ be the union of $\mathcal{V}_i$ and all the sets $F_W$, for $W\in\mathcal{W}_i(>\varepsilon_i^-)\cap\mathcal{L}$. Note that $\mathcal{V}_i^*$ coincides with the union of $\mathcal{V}_i$ and the set
		\[
		\mathcal{V}_i(>\varepsilon_i^-)\land\mathcal{M}_i\restr[\varepsilon_{M_i},\varepsilon_i^+).
		\]
		Therefore, if we show that $\mathcal{V}_i^*$ is an $(\mathcal{S},\mathcal{L})$-symmetric system, by Proposition \ref{prop-pure-amalgamation-2}, $\mathcal{V}_i^*\land\mathcal{M}_i\subseteq\mathcal{V}_i^*$.
		
		Before showing that $\mathcal{V}_i^*$ is an $(\mathcal{S},\mathcal{L})$-symmetric system, let us isolate some of the properties of $\mathcal{V}_i^*$, which follow easily from the properties (p.1)-(p.6) of the sets $F_W$. 
		\begin{enumerate}
			\item[(P.1)] If $V_0,V_1\in\mathcal{V}_i^*$ are such that $\varepsilon_{V_0}=\varepsilon_{V_1}$, then $V_0\in\mathcal{V}_i$ iff $V_1\in\mathcal{V}_i$.
			
			\item[(P.2)] If $V_0,V_1\in\mathcal{V}_i^*$ are such that $\varepsilon_{V_0}=\varepsilon_{V_1}$, and $V_0\in F_{W_0}$ and $V_1\in F_{W_1}$ for some $W_0,W_1\in\mathcal{W}_i(>\varepsilon_i^-)\cap\mathcal{L}$, then $\varepsilon_{W_0}=\varepsilon_{W_1}$. 
			
			\item[(P.3)] If $V\in\mathcal{W}_i\cap W$ for some $W\in\mathcal{W}_i(>\varepsilon_i^-)\cap\mathcal{L}$, then $V\in W\cap M$. In particular, if $V\in\mathcal{V}_i(>\varepsilon_i^-)\cap W$ for some $W\in\mathcal{W}_i(>\varepsilon_i^-)\cap\mathcal{L}$, then $V\in W\cap M$.
		\end{enumerate} 
		
		
		\begin{claim}\label{claim-amalS-1}
			$\mathcal{V}_i^*$ is an $(\mathcal{S},\mathcal{L})$-symmetric system.
		\end{claim}
		\begin{claimproof}
			It is clear that $\mathcal{V}_i^*$ satisfies clause (A) from definition \ref{defSSM2T}.
			
			Let us now check clause (B). So, let $V_0,V_1\in\mathcal{V}_i^*$ such that $\varepsilon_{V_0}=\varepsilon_{V_1}$. We must show that $V_0[\omega_1]\cong V_1[\omega_1]$. If $V_0,V_1\in\mathcal{V}_i$, the conclusion follows from the fact that $\mathcal{V}_i$ is an $(\mathcal{S},\mathcal{L})$-symmetric system. By (P.1), the only other possibility is that $V_0\in F_{W_0}$ and $V_1\in F_{W_1}$, where $W_0,W_1\in\mathcal{W}_i(>\varepsilon_i^-)\cap\mathcal{L}$. By (P.2), $W_0$ and $W_1$ must have the same $\omega_2$-height. Therefore, $V_0[\omega_1]\cong V_1[\omega_1]$ is witnessed by $\Psi_{W_0,W_1}$, by (p.4).
			
			Let us now check that $\mathcal{V}_i^*$ satisfies clause (C). Let $V_0\in\mathcal{V}_i^*$ and suppose that $\varepsilon_1$ is the immediate successor of $\varepsilon_{V_0}$ in $\dom(\mathcal{V}_i^*)$. We need to find a model $V_1\in\mathcal{V}_i^*(\varepsilon_1)$ such that $V_0\in V_1$. Since $\mathcal{V}_i$ is an $(\mathcal{S},\mathcal{L})$-symmetric system by assumption, it suffices to check (C) for $\varepsilon_{V_0}\geq\varepsilon_i^-$. Let us divide the proof in three cases:
			
			\textbf{Case (C.1):} $V_0\in\mathcal{V}_i$ and there is no $W\in\mathcal{W}_i(>\varepsilon_i^-)\cap\mathcal{L}$ such that $\varepsilon_{V_0}$ is the immediate predecessor of $\varepsilon_W$ in $\dom(\mathcal{V}_i)$. Then $\varepsilon_1\in\dom(\mathcal{V}_i)$, and the conclusion follows from the fact that $\mathcal{V}_i$ is an $(\mathcal{S},\mathcal{L})$-symmetric system. 
			
			\textbf{Case (C.2):} $V_0\in\mathcal{V}_i\cap W$ for some $W\in\mathcal{W}_i(>\varepsilon_i^-)\cap\mathcal{L}$, and $\varepsilon_{V_0}$ is the immediate predecessor of $\varepsilon_W$ in $\dom(\mathcal{V}_i)$. Then $\varepsilon_1=\varepsilon_{W\cap M}$ by (p.5), and by (P.3), $V_0\in W\cap M$. Hence, it suffices to let $V_1=W\cap M$. 
			
			\textbf{Case (C.3):} $V_0\in F_W$ for some $W\in\mathcal{W}_i(>\varepsilon_i^-)\cap\mathcal{L}$. If $\varepsilon_{V_0}$ is the top element of $\dom(F_W)$, let $V_1=W$. Otherwise, $\varepsilon_1\in\dom(F_W)$, and the conclusion follows from the fact that the elements of $F_W$ form an $\in$-chain (property (p.3)).
			
			The next step is to prove that $\mathcal{V}_i^*$ satisfies clause (D). Let $V_0,V_1,V_2\in\mathcal{V}_i^*$ such that $V_0\in V_1$ and $\varepsilon_{V_1}=\varepsilon_{V_2}$. We need to prove that $\Psi_{V_1[\omega_1],V_2[\omega_1]}(V_0)$ is a member of $\mathcal{V}_i^*$. If $V_0,V_1,V_2\in\mathcal{V}_i$, then the conclusion follows from the fact that $\mathcal{V}_i$ is an $(\mathcal{S},\mathcal{L})$-symmetric system. In light of (P.1), we can divide the proof in three cases:
			
			\textbf{Case (D.1):} $V_0\in\mathcal{V}_i$, $V_1\in F_{W_1}$ and $V_2\in F_{W_2}$, for some $W_1,W_2$ in $\mathcal{W}_i(>\varepsilon_i^-)\cap\mathcal{L}$. Then the conclusion follows from the fact that, by (p.4), $V_1[\omega_1]\cong V_2[\omega_1]$ is witnessed by $\Psi_{W_1,W_2}$.
			
			\textbf{Case (D.2):} $V_0\in F_{W_0}$ for some $W_0\in\mathcal{W}_i(>\varepsilon_i^-)\cap\mathcal{L}$, and $V_1,V_2\in\mathcal{V}_i$. If $\varepsilon_{W_0}=\varepsilon_{V_1}$, then $V_0\in F_{V_1}$, and hence, $\Psi_{V_1,V_2}(V_0)\in F_{V_2}$. If $\varepsilon_{W_0}<\varepsilon_{V_1}$, as $W_0,V_1\in\mathcal{W}_i$, we can find $V_1'\in\mathcal{W}_i(\varepsilon_{V_1})$ such that $W_0\in V_1'[\omega_1]$. Let $W_1:=\Psi_{V_1'[\omega_1],V_1[\omega_1]}(W_0)$, which is an element of $\mathcal{W}_i$ by Proposition \ref{prop6}. Note that as $V_0\in V_1\cap V_1'[\omega_1]$, 
			\[
			V_0=\Psi_{V_1'[\omega_1],V_1[\omega_1]}(V_0)=\Psi_{W_0,W_1}(V_0).
			\]
			Therefore, $V_0\in F_{W_1}$. Let $W_2:=\Psi_{V_1[\omega_1],V_2[\omega_1]}(W_1)$, which is an element of $\mathcal{W}_i$ by similar reasons as before. Then $\Psi_{V_1[\omega_1],V_2[\omega_1]}(V_0)=\Psi_{W_1,W_2}(V_0)$, and $\Psi_{W_1,W_2}(V_0)$ is an element of $F_{W_2}$. This is pictured in Figure \ref{fig9}
			
			\begin{figure}[ht]
				\centering
					\begin{tikzpicture} 
						\tikzset{thin/.style = {line width= 1pt}} 
						
						\coordinate[] (x) at (2,0){};
						\coordinate[] (y) at (0,0.5){};
						
						\node at (4.7,0) {};
						
						\coordinate[] (V0) at (0,0){};
						\node[right=3pt] at (V0){$V_0$};
						\draw[thin] (-0.14,0)--(0.14,0);
						
						\coordinate[] (W0) at (-1.5,1){};
						\node[left=3pt] at (W0){$W_0$};
						\draw[thin,dotted] (V0)--(W0);
						\draw[thin] (-1.64,1)--(-1.36,1);
						
						\coordinate[] (V1') at (-1.5,2){};
						\node[left=3pt] at (V1'){$V_1'$};
						\draw[thin, dotted, -|] (W0)--(V1');
						
						\coordinate[] (W1) at (1.5,1){};
						\node[right=3pt] at (P1){$W_1=\Psi_{V_1'[\omega_1],V_1[\omega_1]}(W_0)$};
						\draw[thin,dotted] (V0)--(W1);
						\draw[thin] (1.36,1)--(1.64,1);
						
						\coordinate[] (V1) at (1.5,2){};
						\node[above=3pt] at (V1){$V_1$};
						\draw[thin,dotted,-|] (W1)--(V1);
						\draw[thin,->] (-1.3,2) to node [above] {$\Psi_{V_1'[\omega_1],V_1[\omega_1]}$} (1.3,2);
						
						\coordinate[] (V0') at (6,0){};
						\node[right=3pt] at (V0'){$\Psi_{V_1[\omega_1],V_2[\omega_1]}(V_0)$};
						\draw[thin] (5.86,0)--(6.14,0);
						
						\coordinate[] (W2) at (6,1){};
						\node[right=3pt] at (W2){$W_2=\Psi_{V_1[\omega_1],V_2[\omega_1]}(W_1)$};
						\draw[thin,dotted,-|] (V0')--(W2);
						
						\coordinate[] (V2) at (6,2){};
						\node[right=3pt] at (V2){$V_2$};
						\draw[thin,dotted,-|] (W2)--(V2);
						\draw[thin,->] (1.8,2) to node [above] {$\Psi_{V_1[\omega_1],V_2[\omega_1]}$} (5.7,2);
						
					\end{tikzpicture}
					\caption{The final configuration of Case (D.2).}
					\label{fig9}
				\end{figure}
			
			\textbf{Case (D.3):} $V_j\in F_{W_j}$ for some $W_j\in\mathcal{W}_i(>\varepsilon_i^-)\cap\mathcal{L}$, for $j\in\{0,1,2\}$. Then $V_1[\omega_1]\cong V_2[\omega_1]$ is witnessed by $\Psi_{W_1,W_2}$, and hence $\Psi_{V_1[\omega_1],V_2[\omega_1]}(V_0)=\Psi_{W_1,W_2}(V_0)$ is an element of $\mathcal{V}_i^*$ by case (D.2).
			
			Lastly, we need to check that $\mathcal{V}_i^*$ satisfies clause (E). Let $Y\in\mathcal{V}_i^*\cap\mathcal{L}$ and $N\in\mathcal{V}_i^*\cap\mathcal{S}$ such that $Y\in N$. We have to show that $Y\cap N\in\mathcal{V}_i^*$. In fact, since $\mathcal{V}_i^*$ is the result of adding new countable models coming from the sets $F_W$ to $\mathcal{V}_i$, we can restrict ourselves to the case $Y\in\mathcal{V}_i$ and $N\in F_W$, for some $W\in\mathcal{W}_i(>\varepsilon_i^-)\cap\mathcal{L}$. Therefore, we can assume that $N$ is of the form $W\cap N'$, where $N'$ is a small model in $E_i$. Recall from the definition of $F_W$, that $E_i$ was chosen to be an $\in$-chain $N_0\in\dots \in N_k$ of $\mathcal{M}_i\restr[\varepsilon_{M_i},\varepsilon_i^+)$ of maximal length such that $N_0=M_i$. 
			
			Let us suppose that $i>0$, since the case $i=0$ is simpler, and divide the proof in three cases:
			 
			\textbf{Case (E.1):} $\varepsilon_Y>\varepsilon_i^-$. Then $Y\cap N=Y\cap(W\cap N')=Y\cap N'$, which is an element of $F_Y$, and hence $Y\cap N\in\mathcal{V}_i^*$.
			
			\textbf{Case (E.2):} $\varepsilon_Y=\varepsilon_i^-=\varepsilon_{X_{i-1}}$. Note that, since there are no large models in $\mathcal{M}_i\restr(\varepsilon_{X_{i-1}},\varepsilon_i^+)$ by the definition of $\varepsilon_i^+$, $X_{i-1}\in N'$. Therefore $X_{i-1}\cap N'\in\mathcal{M}_i$. In fact, since $\mathcal{V}_i$ extends $\mathcal{M}_i\cap M_i[\omega_1]$ by the inductive hypothesis, $X_{i-1}\cap N'$ is also an element of $\mathcal{V}_i$. Hence, as $\mathcal{V}_i$ is an $(\mathcal{S},\mathcal{L})$-symmetric system, $\Psi_{X_{i-1},Y}(X_{i-1}\cap N')\in\mathcal{V}_i$, and we are done by noting that
			\[
			\Psi_{X_{i-1},Y}(X_{i-1}\cap N')=Y\cap N'=Y\cap(W\cap N')=Y\cap N
			\]
			by Proposition \ref{prop9}.
			
			\textbf{Case (E.3):} $\varepsilon_Y<\varepsilon_i^-=\varepsilon_{X_{i-1}}$. Since both $Y$ and $X_{i-1}$ are elements of $\mathcal{V}_i$, there must be some $Z\in\mathcal{V}_i(\varepsilon_i^-)$ such that $Y\in Z$ by the shoulder axiom. Since $\mathcal{V}_i(\geq\varepsilon_i^-)=\mathcal{W}_i(\geq\varepsilon_i^-)$ and $\mathcal{V}_i\subseteq M_i[\omega_1]$ by the inductive hypothesis, $Z$ must be an element of $M_i$, and hence, a member of $N'$ as well (recall that $N'\in E_i$). Note that $Z\cap N'$ is an element of $\mathcal{V}_i$ by case (E.2). Therefore, since $\mathcal{V}_i$ is an $(\mathcal{S},\mathcal{L})$-symmetric system and $Y\in Z\cap N'$, $$Y\cap N=Y\cap(W\cap N')=Y\cap N'=Y\cap (Z\cap N')\in\mathcal{V}_i.$$	
		\end{claimproof}
		
		Therefore, as we anticipated, now that we have shown that $\mathcal{V}_i^*$ is an $(\mathcal{S},\mathcal{L})$-symmetric system such that 
		\[
		\mathcal{V}_i^*(>\varepsilon_i^-)\land\mathcal{M}_i\restr[\varepsilon_{M_i},\varepsilon_i^+)\subseteq\mathcal{V}_i^*,
		\]
		Proposition \ref{prop-pure-amalgamation-2} ensures that $\mathcal{V}_i^*\land\mathcal{M}_i\subseteq\mathcal{V}_i^*$. Therefore, 
		\[
		\mathcal{U}_i:=\mathcal{M}_i(\geq\varepsilon_{M_i})\cup\{\Psi_{M_i[\omega_1],N_i[\omega_1]}(V):V\in\mathcal{V}_i^*,N_i\in\mathcal{M}_i(\varepsilon_{M_i})\}
		\]
		is an $(\mathcal{S},\mathcal{L})$-symmetric system such that $\mathcal{M}_i\cup\mathcal{V}_i^*\subseteq\mathcal{U}_i\subseteq X_i$ by Proposition \ref{prop-pure-amalgamation-3}, and in particular, $\mathcal{W}_i\subseteq\mathcal{U}_i$.
		
		If $i=n$, we have obtained an $(\mathcal{S},\mathcal{L})$-symmetric system $\mathcal{U}_n$ that extends $\mathcal{M}_n=\mathcal{M}$ and $\mathcal{W}_n=\mathcal{W}$, and we are done.
		
		If $i<n$, since $\mathcal{W}_{i+1}\cap X_i=\mathcal{W}_i\subseteq\mathcal{U}_i\subseteq X_i$, we have $\mathcal{U}_i\land\mathcal{W}_{i+1}\subseteq\mathcal{U}_i$ by Proposition \ref{prop-pure-amalgamation-1}. So, by Proposition \ref{prop-pure-amalgamation-3},
		\[
		\mathcal{V}_{i+1}:=\mathcal{W}_{i+1}(\geq\varepsilon_{X_i})\cup\{\Psi_{X_i,Y_i}(U):U\in\mathcal{U}_i,Y_i\in\mathcal{W}_{i+1}(\varepsilon_{X_i})\}
		\]
		is an $(\mathcal{S},\mathcal{L})$-symmetric system such that $\mathcal{W}_{i+1}\cup\mathcal{U}_i\subseteq\mathcal{V}_{i+1}$. Hence, we will be done by showing that $\mathcal{V}_{i+1}$ has the properties (a)-(c) that we assumed for $\mathcal{V}_i$ at the beginning of the inductive construction of the $(\mathcal{S},\mathcal{L})$-symmetric systems $\mathcal{U}_i$. That is, we have to check that $\mathcal{V}_{i+1}$ satisfies the following:
		\begin{enumerate}[label=(\alph*)]
			\item $\mathcal{M}_{i+1}\cap M_{i+1}[\omega_1]\subseteq\mathcal{V}_{i+1}\subseteq M_{i+1}[\omega_1]$.
			
			\item $\mathcal{W}_{i+1}\subseteq\mathcal{V}_{i+1}$.
			
			\item $\mathcal{V}_{i+1}(\geq\varepsilon_{X_i})=\mathcal{W}_{i+1}(\geq\varepsilon_{X_i})$.
		\end{enumerate}
		
		Item (a) follows from $\mathcal{M}_{i+1}(\geq\varepsilon_{X_i})\cap M_{i+1}[\omega_1]\subseteq\mathcal{W}_{i+1}(\geq\varepsilon_{X_i})$, and the fact that, since $\mathcal{M}_i\subseteq\mathcal{U}_i$, we have
		\begin{align*}
			\mathcal{M}_{i+1}(<\varepsilon_{X_i})\cap M_{i+1}[\omega_1]&=\{\Psi_{X_i,Y_i}(Q):Q\in\mathcal{M}_i,Y_i\in\mathcal{M}_{i+1}(\varepsilon_{X_i})\}\\
			&\subseteq\{\Psi_{X_i,Y_i}(U):U\in\mathcal{U}_i,Y_i\in\mathcal{W}_{i+1}(\varepsilon_{X_i})\}\\
			&\subseteq\mathcal{V}_{i+1}.
		\end{align*}
		Now, since $\mathcal{W}_i\subseteq\mathcal{U}_i$, we have 
		\begin{align*}
			\mathcal{W}_{i+1}(<\varepsilon_{X_i})&=\{\Psi_{X_i,Y_i}(W):W\in\mathcal{W}_i,Y_i\in\mathcal{W}_{i+1}(\varepsilon_{X_i})\}\\
			&\subseteq\{\Psi_{X_i,Y_i}(U):U\in\mathcal{U}_i,Y_i\in\mathcal{W}_{i+1}(\varepsilon_{X_i})\}\\
			&\subseteq\mathcal{V}_{i+1},
		\end{align*}
		which implies item (b). Lastly, item (c) follows directly from the definition of $\mathcal{V}_{i+1}$.	
	\end{proof}

	\subsection{Preservation lemmas}\label{subsection-pure-preservation}
	
	In this subsection, we will show how the results from the last section can be used to prove that, under the right assumptions, all cardinals and $2^{\aleph_1}=\aleph_2$ are preserved when forcing with $(\mathcal{S},\mathcal{L})$-symmetric systems.
	
	\begin{definition}
		Let $\mathbb{M}(\mathcal{S},\mathcal{L})$ be the forcing notion whose conditions are $(\mathcal{S},\mathcal{L})$-symmetric systems and the order is reverse inclusion.
	\end{definition}
	
	\begin{theorem}\label{pureproperS}
		The forcing $\mathbb{M}(\mathcal{S},\mathcal{L})$ is strongly $\mathcal{S}$-proper.
	\end{theorem}
	\begin{proof}
		
		Let $\mathcal{M}\in\mathbb{M}(\mathcal{S},\mathcal{L})$ and $M\in\mathcal{S}$ such that $\mathcal{M}\in M$. Then by Lemma \ref{amal-ontop}, there is $\mathcal{M}_M\in\mathbb{M}(\mathcal{S},\mathcal{L})$ such that $M\in\mathcal{M}_M$ and $\mathcal{M}_M\supseteq\mathcal{M}$. We claim that $\mathcal{M}_M$ is strongly $(M,\mathbb{M}(\mathcal{S},\mathcal{L}))$-generic. Let $D\subseteq\mathbb{M}(\mathcal{S},\mathcal{L})\cap M$ be a dense subset, and let $\mathcal{N}\in\mathbb{M}(\mathcal{S},\mathcal{L})$ such that $\mathcal{N}\supseteq\mathcal{M}_M$. From Lemma \ref{amal-rest} it follows that $\mathcal{N}\cap M$ is a condition in $\mathbb{M}(\mathcal{S},\mathcal{L})\cap M$. Since $D$ is dense, there is some $\mathcal{W}\in D$ such that $\mathcal{W}\supseteq\mathcal{N}\cap M$, and as $D\subseteq M$, we have that $\mathcal{W}\in M$ (and $\mathcal{W}\subseteq M$ because $\mathcal{W}$ is finite). Therefore, $\mathcal{N}$ and $\mathcal{W}$ are compatible by Lemma \ref{lemma-amal-S}, and therefore we can conclude that $\mathcal{M}_M$ is strongly $(M,\mathbb{M}(\mathcal{S},\mathcal{L}))$-generic as we wanted.
	\end{proof}
	
	The proof of the following lemma is exactly the same as the proof of Lemma \ref{pureproperS}, but uses Lemma \ref{lemma-amal-L} instead of Lemma \ref{lemma-amal-S}.
	
	\begin{theorem}\label{pureproperL}
		The forcing $\mathbb{M}(\mathcal{S},\mathcal{L})$ is strongly $\mathcal{L}$-proper.
	\end{theorem}
	
	From the proof of the last two theorems, we can extract the following result.
	
	\begin{lemma}\label{lemma-pure-genericity}
		Any condition $\mathcal{M}\in\mathbb{M}(\mathcal{S},\mathcal{L})$ is strongly $(Q,\mathbb{M}(\mathcal{S},\mathcal{L}))$-generic for every $Q\in\mathcal{M}$.
	\end{lemma}
	
	\begin{theorem}\label{purecc}
		If $2^{\aleph_1}=\aleph_2$, then the forcing $\mathbb{M}(\mathcal{S},\mathcal{L})$ has the $\aleph_3$-Knaster condition.
	\end{theorem}
	\begin{proof}
		For every $\alpha<\omega_3$, let $\mathcal{M}_\alpha$ be an $(\mathcal{S},\mathcal{L})$-symmetric system and let $M_\alpha\in\mathcal{S}$ with $\mathcal{M}_\alpha\in M_\alpha$. By $2^{\aleph_1}=\aleph_2$, we may assume that the set $\{M_\alpha[\omega_1]:\alpha<\omega_3\}$ forms a $\Delta$-system with root $\mathcal{X}$, and that for any two distinct $\alpha,\beta\in\omega_3$, there is a unique $\omega_1$-isomorphism
		\[
		\Psi_{M_\alpha[\omega_1],M_\beta[\omega_1]}:(M_\alpha[\omega_1];\in,M_\alpha,\mathcal{M}_\alpha,\lhd)\to(M_\beta[\omega_1];\in,M_\beta,\mathcal{M}_\beta,\lhd).
		\]
		For any two $\alpha,\beta<\omega_3$, let $\Psi_{\alpha,\beta}$ denote the map $\Psi_{M_\alpha[\omega_1],M_\beta[\omega_1]}\restr\bigcup\mathcal{M}_\alpha[\omega_1]$. Note that $\Psi_{\alpha,\beta}$ is the identity on $\mathcal{X}$. The reason is that since there is a definable bijection between $H(\omega_2)$ and $\omega_2$ (this also follows from $2^{\aleph_1}=\aleph_2$), $\Psi_{\alpha,\beta}$ fixes $\mathcal{X}$ if and only if it fixes $\mathcal{X}\cap\omega_2$. Hence, $\Psi_{\alpha,\beta}$ witnesses $\mathcal{M}_\alpha\cong\mathcal{M}_\beta$. Therefore, $\mathcal{M}_\alpha\cup\mathcal{M}_\beta$ is an $(\mathcal{S},\mathcal{L})$-symmetric system according to Lemma \ref{pureamalgamation2} and, in particular, $\mathcal{M}_\alpha$ and $\mathcal{M}_\beta$ are compatible in $\mathbb{M}(\mathcal{S},\mathcal{L})$.
	\end{proof}
	
	Therefore, combining the previous results with Lemmas \ref{preservacio-proper} and \ref{preservation-cc}, we obtain the following preservation theorem.
	
	\begin{corollary}
		If $2^{\aleph_1}=\aleph_2$ holds and $\mathcal{L}$ is stationary in $H(\kappa)$, then $\mathbb{M}(\mathcal{S},\mathcal{L})$ preserves all cardinals.
	\end{corollary}

	\begin{theorem}\label{purepreservationCH}
		If $\mathcal{L}$ is stationary in $H(\kappa)$, then $\mathbb{M}(\mathcal{S},\mathcal{L})$ preserves $2^{\aleph_1}=\aleph_2$.
	\end{theorem}
	\begin{proof}
		Assume that $2^{\aleph_1}=\aleph_2$ holds and let $\langle\tau_\alpha:\alpha<\omega_3\rangle$ be a sequence of
		$\mathbb{M}(\mathcal{S},\mathcal{L})$-names for subsets of $\omega_1$. Suppose that $\mathcal{M}\in\mathbb{M}(\mathcal{S},\mathcal{L})$ is a condition forcing that $\langle\tau_\alpha:\alpha<\omega_3\rangle$ is a sequence of pairwise different subsets of $\omega_1$. For every $\alpha<\omega_3$, let $X_\alpha^*$ be an elementary submodel of a large enough $H(\theta)$, for $\theta>\kappa$, such that $\mathbb{M}(\mathcal{S},\mathcal{L}),\tau_\alpha,\mathcal{M}\in X_\alpha^*$ and $X_\alpha:=X_\alpha^*\cap H(\kappa)\in\mathcal{L}$. We may assume that there are two different $\alpha,\beta<\omega_3$ for which the structures $(X_\alpha^*;\in,\tau_\alpha)$ and $(X_\beta^*;\in,\tau_\beta)$ are isomorphic, and the corresponding isomorphism fixes $X_\alpha^*\cap X_\beta^*$ and sends $\tau_\alpha$ to $\tau_\beta$. This follows from the fact that, as $2^{\aleph_1}=\aleph_2$, there are only $\aleph_2$-many isomorphism types for such structures. Let $\mathcal{M}_{\alpha,\beta}:=\mathcal{M}\cup\{X_\alpha,X_\beta\}$, which is easily seen to be an $(\mathcal{S},\mathcal{L})$-symmetric system. We only need to notice that since $\mathcal{M}\in X_\alpha\cap X_\beta$, then $\Psi_{X_\alpha^*,X_\beta^*}(\mathcal{M})=\mathcal{M}$. Moreover, according to Lemma \ref{lemma-pure-genericity}, $\mathcal{M}_{\alpha,\beta}$ is strongly $(X_\alpha,\mathbb{M}(\mathcal{S},\mathcal{L}))$-generic and strongly $(X_\beta,\mathbb{M}(\mathcal{S},\mathcal{L}))$-generic. 
		
		Now, recall that $\mathcal{M}$ forces that $\tau_\alpha$ and $\tau_\beta$ are different. Hence, there must be an $(\mathcal{S},\mathcal{L})$-symmetric system $\mathcal{N}\supseteq\mathcal{M}_{\alpha,\beta}$ and an ordinal $\gamma<\omega_1$, such that $\mathcal{N}\Vdash\check{\gamma}\in\tau_\alpha\setminus\tau_\beta$. Let $D\subseteq\mathbb{M}(\mathcal{S},\mathcal{L})\cap X_\alpha$ be the dense set of conditions deciding whether $\check{\gamma}$ is an element of $\tau_\alpha$ or not. Since $\mathcal{M}_{\alpha,\beta}$ is strongly $(X_\alpha,\mathbb{M}(\mathcal{S},\mathcal{L}))$-generic, $\mathcal{N}$ is also strongly $(X_\alpha,\mathbb{M}(\mathcal{S},\mathcal{L}))$-generic, and hence $D$ must be predense below $\mathcal{N}$. Therefore, there are conditions $\mathcal{W}\in D$ and $\mathcal{U}\in\mathbb{M}(\mathcal{S},\mathcal{L})$ such that $\mathcal{U}$ extends $\mathcal{W}\cup\mathcal{N}$. On the one hand, note that $\mathcal{W}$ forces $\check{\gamma}\in\tau_\alpha$. Hence, since $\Psi_{X_\alpha^*,X_\beta^*}$ is an isomorphism between large models, $\check{\gamma}$ is a name for an ordinal of $\omega_1$, and $\mathcal{W}\in X_\alpha$, we have 
		\[
		\Psi_{X_\alpha^*,X_\beta^*}(\mathcal{W})\Vdash\check{\gamma}\in\Psi_{X_\alpha^*,X_\beta^*}(\tau_\alpha)=\tau_\beta.
		\]
		On the other hand, note that since $\mathcal{U}$ extends $\mathcal{N}$, the models $X_\alpha$ and $X_\beta$ must be members of $\mathcal{U}$. Therefore, since $\mathcal{W}\subseteq\mathcal{U}$ and $\mathcal{W}\in X_\alpha$, then $\Psi_{X_\alpha^*,X_\beta^*}(\mathcal{W})\subseteq\mathcal{U}$ by the symmetry of $\mathcal{U}$. Therefore, $\mathcal{U}\Vdash\check{\gamma}\in\tau_\beta$. But this is impossible because $\mathcal{U}$ is an extension of $\mathcal{N}$, and $\mathcal{N}$ forces that $\check{\gamma}\in\tau_\alpha\setminus\tau_\beta$. Therefore, we can conclude that there is no condition $\mathcal{M}$ forcing that $\langle\tau_\alpha:\alpha<\omega_3\rangle$ is a sequence of pairwise different subsets of $\omega_1$.				
	\end{proof}	
	
	\subsection{First application: Forcing an $\omega_1$-club subset of $\omega_2$}\label{subsection-pure-first-application}
	
	As a first application of the forcing $\mathbb{M}(\mathcal{S},\mathcal{L})$, we will show how it adds a two-type symmetric system that covers all of $H(\kappa)^V$ and an $\omega_1$-club subset of $\omega_2$.
	
	Let $G$ be an $\mathbb{M}(\mathcal{S},\mathcal{L})$-generic filter over $V$, and define the sets
	\[
	\mathcal{M}_G:=\bigcup\{\mathcal{M}\in\mathbb{M}(\mathcal{S},\mathcal{L}):\mathcal{M}\in G\}
	\] 
	and
	\[
	C_G:=\{\varepsilon_X:X\in\mathcal{M}_G\cap\mathcal{L}\}=\dom(\mathcal{M}_G\cap\mathcal{L}).
	\]
	
	It is not too hard to see that $\mathcal{M}_G$ satisfies clauses (A)-(E) from the definition of $(\mathcal{S},\mathcal{L})$-symmetric system. Moreover, for every $x\in H(\kappa)^V$, the set
	\[
	D_x=\{\mathcal{M}\in\mathbb{M}(\mathcal{S},\mathcal{L}):\exists Q\in\mathcal{M}(x\in Q)\}
	\]
	is dense in $\mathbb{M}(\mathcal{S},\mathcal{L})$, and hence, by genericity, $\mathcal{M}_G$ is unbounded in $H(\kappa)^V$ (in the sense of Definition \ref{def-unbd-stat}). In particular, since $\mathcal{S}$ is a club subset of $[H(\kappa)]^\omega$, $\mathcal{M}_G\cap\mathcal{S}$ is always unbounded in $H(\kappa)^V$. Moreover, if $\mathcal{L}$ is unbounded (in particular, stationary) in $H(\kappa)^V$, then $\mathcal{M}_G\cap\mathcal{L}$ is also unbounded in $H(\kappa)^V$.  
	
	
	It is also worth mentioning that, even though the poset $\mathbb{M}(\mathcal{S},\mathcal{L})$ and Neeman's two-type side condition forcing $\mathbb{P}_{\omega,\mathcal{S},\mathcal{L},H(\kappa)}$ (\cite{Neeman2014Forcingwithsequencesofmodelsoftwotypes}) have some similarities, they produce different generic objects in general. For instance, while Neeman's generic $G'$ for $\mathbb{P}_{\omega,\mathcal{S},\mathcal{L},H(\kappa)}$ would entail $|H(\kappa)^V|=\aleph_2^{V[G']}$, and thus collapse $\kappa$, the poset $\mathbb{M}(\mathcal{S},\mathcal{L})$ preserves all cardinals, under the assumptions that $\mathcal{L}$ is stationary and $2^{\aleph_1}=\aleph_2$. Therefore, $G'$ cannot exist in the generic extension by $\mathbb{M}(\mathcal{S},\mathcal{L})$.


	\begin{definition}
		We say that $C\subseteq\omega_2$ is an \emph{$\omega_1$-club of $\omega_2$} if it is unbounded in $\omega_2$ and closed at ordinals of cofinality $\omega_1$, that is, for every $\nu\in\omega_2$, if $C\cap\nu$ is cofinal in $\nu$ and $\cf(\nu)=\omega_1$, then $\nu\in C$.
	\end{definition}
	
	The proof of the following theorem is standard in the context of forcing with side conditions. You can find a proof of the same result for Neeman's forcing with two-type chains of models as Claim 2.4 in \cite{Neeman2017Twoapplicationsoffinitesideconditionsatomega2}. 
	
	\begin{theorem}\label{theorem-pure-club}
		If $\mathcal{L}$ is stationary in $H(\kappa)$, then $C_G$ is an $\omega_1$-club subset of $\omega_2$.
	\end{theorem}
	\begin{proof}
		$C_G$ is clearly unbounded by the unboundedness of $\mathcal{L}$ and the genericity of $G$. 
		
		By Lemma \ref{amal-ontop} and the unboundedness of $\mathcal{L}$, for every $\varepsilon<\omega_2$, the set of conditions $\mathcal{M}\in\mathbb{M}(\mathcal{S},\mathcal{L})$ such that $\varepsilon\in X$, for some $X\in\mathcal{M}\cap\mathcal{L}$, is dense. Hence, by genericity there is some $\varepsilon'\in C_G$ such that $\varepsilon'>\varepsilon$. It follows that $C_G$ is unbounded in $\omega_2$.
		
		Let $\nu\in\omega_2$ be such that $C_G\cap\nu$ is cofinal in $\nu$ and $cf(\nu)=\omega_1$. We will show that $\nu\in C_G$. Let $\varepsilon_0$ be the least ordinal in $\dom(\mathcal{M}_G)$ such that $\nu\leq\varepsilon_0$. Let $Q\in\mathcal{M}_G(\varepsilon_0)$. If $Q$ was a small model, the set $$\dom(\mathcal{M}_G\cap Q)=\{\varepsilon_P:P\in\mathcal{M}_G\cap Q\}$$ would be countable, and hence bounded in $\nu$. Note that $$\dom(\mathcal{M}_G\cap Q[\omega_1])=\dom(\mathcal{M}_G)\cap\varepsilon_0$$ is unbounded in $\nu$. 
		
		We will reach a contradiction by showing that $\dom(\mathcal{M}_G\cap Q)$ is cofinal in $\dom(\mathcal{M}_G\cap Q[\omega_1])$. Let $P\in \mathcal{M}_G\cap Q[\omega_1]$, and suppose that $P,Q\in\mathcal{M}$, where $\mathcal{M}\in G$. If $P\in Q$, we are done. Hence, suppose that $P\notin Q$. If $\varepsilon_Q$ was the successor of $\varepsilon_P$ in $\dom(\mathcal{M})$, then $P$ would be a member of $Q$. So, if we let $\varepsilon$ be the maximum of $\dom(\mathcal{M})\cap\varepsilon_Q$, we have $\varepsilon_P<\varepsilon<\varepsilon_Q$. By Proposition \ref{prop13}, we can find $R\in\mathcal{M}(\varepsilon)$ such that $P\in R[\omega_1]$ and $R\in Q[\omega_1]$. But since $\varepsilon_Q$ is the immediate successor of $\varepsilon_R$ in $\dom(\mathcal{M})$, we have $R\in Q$. Therefore, we have found $R\in\mathcal{M}_G\cap Q$ such that $\varepsilon_P<\varepsilon_R<\varepsilon_Q=\varepsilon_0$. This shows that $\dom(\mathcal{M}_G\cap Q)$ is cofinal in $\dom(\mathcal{M}_G\cap Q[\omega_1])$ as we wanted, and hence we reach a contradiction. Therefore, we can conclude that $Q$ has to be a large model.
		
		If $\varepsilon_Q=\nu$, then $\nu\in C_G$ and we are done. Suppose towards a contradiction that $\nu<\varepsilon_0$. Since $Q$ is a large model, $\varepsilon_0=\varepsilon_Q=Q\cap\omega_2$, and hence $\nu\in Q$. Let $M\in G$ be a small model such that $\nu,Q\in M$. It exists by stationarity of $\mathcal{S}$, and genericity of $G$. By the closure under intersections of $(\mathcal{S},\mathcal{L})$-symmetric systems, $Q\cap M\in\mathcal{M}_G$. But then $\nu<\varepsilon_{Q\cap M}<\varepsilon_Q=\varepsilon_0$, which contradicts the minimality of the ordinal $\varepsilon_0$.
	\end{proof}
	
	As we will see in Section \ref{section-decorations}, if we want to improve the conclusion of the last theorem by adding a club subset of $\omega_2$, we need to modify the forcing $\mathbb{M}(\mathcal{S},\mathcal{L})$ by adding working parts to it, which in Neeman's terminology are called \emph{decorations} (\cite{Neeman2014Forcingwithsequencesofmodelsoftwotypes}).

	\subsection{Second application: Forcing a Kurepa tree on $\omega_2$}
	
	Recall that a \emph{Kurepa tree on $\omega_2$} is a tree on $\omega_2$ such that each level has cardinality less than $\aleph_2$, but has at least $\aleph_3$ cofinal branches. In this subsection, we will show how $\mathbb{M}(\mathcal{S},\mathcal{L})$ also adds a Kurepa tree on $\omega_2$ with $\kappa$-many branches. The argument follows the same ideas that were used by Kuzeljevi\'c and Todor\v cevi\'c in \cite{KuzeljevicTodorcevic2017Forcingwithmatrices} to show that the forcing consisting of finite matrices of countable elementary submodels adds a Kurepa tree (on $\omega_1$).
	
	It is worth mentioning that the forcing $\mathbb{M}^{mor}(\mathcal{S},\mathcal{L})$ from Section \ref{section-morass}, adds a simplified $(\omega_2,1)$-morass, which in particular implies the existence of a Kurepa tree on $\omega_2$. However, we wanted to add this result here to illustrate how the most simple version of the forcing with $(\mathcal{S},\mathcal{L})$-symmetric systems is interesting in its own right.
	
	For the remainder of the section suppose that $\mathcal{L}$ is stationary in $H(\kappa)$ and that $2^{\aleph_1}=\aleph_2$ holds in $V$. For every $Q\in\mathcal{S}\cup\mathcal{L}$, let $\Psi_{Q}$ denote the transitive collapsing map of $Q$. Let $G$ be an $\mathbb{M}(\mathcal{S},\mathcal{L})$-generic filter over $V$. Work in $V[G]$ for the rest of the subsection. 
	
	\begin{definition}
		For every $\alpha<\kappa$, define the function $b_\alpha:\omega_2\to\omega_2$ by letting $b_\alpha(\varepsilon)=\nu$, if there is some $Q\in \mathcal{M}_G$ such that $\varepsilon_Q=\varepsilon$, $\alpha\in Q$ and $\Psi_Q(\alpha)=\nu$, or $b_\alpha(\varepsilon)=0$, otherwise.
	\end{definition} 
	
	\begin{lemma}
		For every $\alpha<\kappa$, $b_\alpha$ is a well-defined function.
	\end{lemma}
	\begin{proof}
		Let $Q_0,Q_1\in\mathcal{M}_G$ be such that $\varepsilon_{Q_0}=\varepsilon_{Q_1}$. Then $Q_0$ and $Q_1$ are $\omega_1$-isomorphic witnessed by $\Psi_{Q_0[\omega_1],Q_1[\omega_1]}$. In particular, if $\alpha\in Q_0\cap Q_1\cap\kappa$, then $\Psi_{Q_0[\omega_1],Q_1[\omega_1]}(\alpha)=\alpha$, and so
		\[
		\Psi_{Q_0[\omega_1],Q_1[\omega_1]}"(Q_0\cap\alpha)=Q_1\cap\alpha.
		\]
		We can conclude that $\Psi_{Q_0}(\alpha)=\Psi_{Q_1}(\alpha)$. 
	\end{proof}
	
	\begin{lemma}
		If $\alpha<\beta<\kappa$, then $b_\alpha\neq b_\beta$.
	\end{lemma}
	\begin{proof}
		Let $Q\in\mathcal{S}\cup\mathcal{L}$ such that $\alpha,\beta\in Q$. Then $b_\alpha(\varepsilon_Q)=\Psi_Q(\alpha)\neq\Psi_Q(\beta)= b_\beta(\varepsilon_Q)$.
	\end{proof}
	
	Define $B=\{b_\alpha:\alpha<\kappa\}$. For every $\alpha<\kappa$, $B_\alpha=\{b_\alpha\restr\varepsilon:\varepsilon<\omega_2\}$ will be the $\alpha$-th branch of the Kurepa tree. Hence, the Kurepa tree will be given by $T_G:=\bigcup_{\varepsilon<\omega_2}T_\varepsilon$, where $T_\varepsilon=B\restr\varepsilon=\{b_\alpha\restr\varepsilon:\alpha<\kappa\}$ will denote the $\varepsilon$-th level of the tree, for each $\varepsilon<\omega_2$. 
	
	\begin{theorem}
		$T_G$ is a Kurepa tree on $\omega_2$.
	\end{theorem}
	\begin{proof}
		Note that, by the previous results, we are only missing to prove that, for every $\varepsilon<\omega_2$, the level $T_\varepsilon$ has size less than $\aleph_2$. Fix some $\varepsilon<\omega_2$. By the stationarity of $\mathcal{L}$ and genericity, there is an elementary submodel $X^*$ of a large enough $H(\theta)$, for $\theta>\kappa$, such that $\mathbb{M}(\mathcal{S},\mathcal{L}),\varepsilon\in X^*$ and $X:=X^*\cap H(\kappa)\in\mathcal{M}_G\cap\mathcal{L}$. We will show the equality
		\[
		T_\varepsilon=\{b_\alpha\restr\varepsilon:\alpha\in X\cap\kappa\},
		\]
		which will imply, in turn, that $T_\varepsilon$ has size at most $|X|=\aleph_1$. 
		
		The inclusion $T_\varepsilon\supseteq\{b_\alpha\restr\varepsilon:\alpha\in X\cap\kappa\}$ is clear. Let us show the other direction. Let $\beta<\kappa$ such that $b_\beta\restr\varepsilon\in T_\varepsilon$. We claim that there is $\alpha\in X\cap\kappa$ such that  $b_\beta\restr\varepsilon=b_\alpha\restr\varepsilon$. Let us divide the proof in two cases.
		
		\textbf{Case 1.} Suppose that for every $\varepsilon'<\varepsilon$, $b_\beta(\varepsilon')=0$. We will show that $b_\beta\restr\varepsilon=b_0\restr\varepsilon$, by showing that $b_0(\varepsilon')=0$ for every $\varepsilon'<\varepsilon$. Note that every model $Q\in\mathcal{S}\cup\mathcal{L}$ contains $0$ as an element. Therefore, if there is some $Q\in\mathcal{M}_G$ such that $\varepsilon_{Q}=\varepsilon'$, then $\Psi_{Q}(0)=0=b_0(\varepsilon')$, and if there is no $Q\in\mathcal{M}_G$ such that $\varepsilon_{Q}=\varepsilon'$, then $b_0(\varepsilon')=0$, by definition of $b_0$.
		
		\textbf{Case 2.} Suppose that there is some $\varepsilon'<\varepsilon$ such that $b_\beta(\varepsilon')\neq 0$. Let $\varepsilon_0$ be the least such $\varepsilon'$. Let $Q\in\mathcal{M}_G$ such that $\varepsilon_Q=\varepsilon_0$, $\beta\in Q$ and $b_\beta(\varepsilon_0)=\Psi_Q(\beta)$. Since $\varepsilon_Q=\varepsilon_0<\varepsilon<\varepsilon_X$, there must be some $Y\in\mathcal{M}_G$ such that $Q\in Y$, and as $Y$ is a large model, $\beta\in Q\subseteq Y$.\footnote{Note that the proof would break here if we had chosen a small model $M$, instead of the large model $X$. The ``transitivity" of the model $Y$, given by Proposition \ref{prelim-trans}, is crucially used at this point.} Denote $\alpha:=\Psi_{Y,X}(\beta)$. We claim that $b_\beta\restr\varepsilon=b_\alpha\restr\varepsilon$.
		
		Let first $\varepsilon_1<\varepsilon_0$. We show that $b_\beta(\varepsilon_1)=0=b_\alpha(\varepsilon_1)$. Suppose that there is some $Q'\in\mathcal{M}_G(\varepsilon_1)$ such that $\alpha\in Q'$. Then we can find a model $X'\in\mathcal{M}_G(\varepsilon_X)$ such that $Q'\in X'$ (recall that $X$ is a large model). Let $Q''=\Psi_{X',X}(Q')$, and note that since $\alpha\in X\cap X'$, we have $\alpha=\Psi_{X',X}(\alpha)\in Q''$. Therefore, $\beta\in\Psi_{X,Y}(Q'')$. But that would contradict the minimality of $\varepsilon_0$, since $\varepsilon_{\Psi_{X,Y}(Q'')}=\varepsilon_1<\varepsilon_0$.
	
		Let now $\varepsilon_2\geq\varepsilon_0$ be such that $\varepsilon_2<\varepsilon$. Suppose first that there is some $P'\in\mathcal{M}_G(\varepsilon_2)$ such that $\beta\in P'$, and observe that $b_\beta(\varepsilon_2)=\Psi_{P'}(\beta)$. By the shoulder axiom, there is some $Y'\in\mathcal{M}_G(\varepsilon_X)$ such that $P'\in Y'$, and hence $P:=\Psi_{Y',Y}(P')\in\mathcal{M}_G\cap Y$ and $\beta=\Psi_{Y',Y}(\beta)\in P$ (note that $\beta\in Y'\cap Y$). Therefore, $b_\beta(\varepsilon_2)=\Psi_{P'}(\beta)=\Psi_P(\beta)$. Now, if we let $P^*:=\Psi_{Y,X}(P)$, which is an element of $\mathcal{M}_G\cap X$, we have $\alpha=\Psi_{Y,X}(\beta)\in P^*$. Since $\alpha=\Psi_{Y,X}(\beta)$ and $P^*=\Psi_{Y,X}(P)$, it is clear that $\Psi_P(\beta)=\Psi_{P^*}(\alpha)$, and hence
		\[
		b_\beta(\varepsilon_2)=\Psi_P(\beta)=\Psi_{P^*}(\alpha)=b_\alpha(\varepsilon_2).
		\]
		Suppose now that there is no $P''\in\mathcal{M}_G(\varepsilon_2)$ such that $\beta\in P''$, and observe that $b_\beta(\varepsilon_2)=0$. Suppose, aiming for a contradiction, that there is a model $P^{**}\in\mathcal{M}_G(\varepsilon_2)$ such that $\alpha\in P^{**}$. Using the symmetry of $\mathcal{M}_G$, we may further assume that $P^{**}\in X$. However, this is impossible, because then we would have $\beta=\Psi_{X,Y}(\alpha)\in\Psi_{X,Y}(P^{**})$. Hence, there cannot be a model $P^{**}\in\mathcal{M}_G(\varepsilon_2)$ such that $\alpha\in P^{**}$, and we can conclude that
		\[
		b_\beta(\varepsilon_2)=0=b_\alpha(\varepsilon_2).
		\]
		
		This finishes the proof of $b_\beta\restr\varepsilon=b_\alpha\restr\varepsilon$ of the second case.
	\end{proof}

	\section{Forcing with decorations}\label{section-decorations}
	
	In this section we will introduce a variant of the forcing $\mathbb{M}(\mathcal{S},\mathcal{L})$, which we will denote $\mathbb{M}^{dec}(\mathcal{S},\mathcal{L})$, inspired by Neeman's decorated version of his two-type side conditions poset (see Definition 2.35 from \cite{Neeman2014Forcingwithsequencesofmodelsoftwotypes}). As we mentioned at the end of Subsection \ref{subsection-pure-first-application}, the main reason to introduce the poset $\mathbb{M}^{dec}(\mathcal{S},\mathcal{L})$ is to improve the result from Theorem \ref{theorem-pure-club}. Namely, by proving that the decorated forcing notion adds a club subset of $\omega_2$, instead of just an $\omega_1$-club subset of $\omega_2$.
	
	The first part of the section will be devoted to proving that the forcing $\mathbb{M}^{dec}(\mathcal{S},\mathcal{L})$ satisfies the same amalgamation lemmas as the pure side condition forcing $\mathbb{M}(\mathcal{S},\mathcal{L})$, and hence that it has the same properties. From the amalgamation lemmas will follow that $\mathbb{M}^{dec}(\mathcal{S},\mathcal{L})$ is strongly $\mathcal{S}$-proper and strongly $\mathcal{L}$-proper, that it has the $\aleph_3$-Knaster condition, and that it preserves $2^{\aleph_1}=\aleph_2$. The second part of the section will be devoted to proving the two main applications of the forcing $\mathbb{M}^{dec}(\mathcal{S},\mathcal{L})$. Namely, we will prove that the poset $\mathbb{M}^{dec}(\mathcal{S},\mathcal{L})$ adds a club subset of $\omega_2$, denoted $C_G^{dec}$, and a function on $\omega_2$ that bounds every canonical function below $\omega_3$ on the generic club $C_G^{dec}$.
	
	\begin{definition}\label{def-deco}
		Let $\mathbb{M}^{dec}(\mathcal{S},\mathcal{L})$ be the poset whose conditions are pairs $(\mathcal{M},d)$ satisfying the following properties:
		\begin{enumerate}
			\item $\mathcal{M}\in\mathbb{M}(\mathcal{S},\mathcal{L})$.
			\item $d:\dom(\mathcal{M})\to[H(\kappa)]^{<\omega}$ is a function, called a \emph{decoration on $\mathcal{M}$}, with the following properties:
			\begin{enumerate}
				\item[(2.a)] If $\varepsilon_0$ is the immediate predecessor of $\varepsilon_1$ in $\dom(\mathcal{M})$, then for every $x\in d(\varepsilon_0)$, there is some $Q\in\mathcal{M}(\varepsilon_1)$ such that $x\in Q$.
				
				\item[(2.b)] Let $Q_0,Q_1\in\mathcal{M}$ and $\varepsilon\in\dom(\mathcal{M})$ such that $\varepsilon<\varepsilon_{Q_0}=\varepsilon_{Q_1}$. If $x\in d(\varepsilon)\cap Q_0$, then $\Psi_{Q_0[\omega_1],Q_1[\omega_1]}(x)\in d(\varepsilon)$.
				
				
			\end{enumerate} 
		\end{enumerate}
		If $(\mathcal{M},d),(\mathcal{N},f)\in\mathbb{M}^{dec}(\mathcal{S},\mathcal{L})$, we let $(\mathcal{N},f)\leq(\mathcal{M},d)$ if and only if
		\begin{enumerate}
			\item $\mathcal{N}\supseteq\mathcal{M}$, and
			\item $f(\varepsilon)\supseteq d(\varepsilon)$, for every $\varepsilon\in\dom(\mathcal{M})$.
		\end{enumerate}
	\end{definition}
	
	The next proposition, which strengthens the symmetry given by clause (2.b) in the definition of $\mathbb{M}^{dec}(\mathcal{S},\mathcal{L})$, is proven exactly like Proposition \ref{prop6}.

	\begin{proposition}
		Let $(\mathcal{M},d)\in\mathbb{M}^{dec}(\mathcal{S},\mathcal{L})$ and suppose that $Q_0,Q_1\in\mathcal{M}$ and $\varepsilon\in\dom(\mathcal{M})$ are such that $\varepsilon<\varepsilon_{Q_0}=\varepsilon_{Q_1}$. If $x\in d(\varepsilon)\cap Q_0[\omega_1]$, then $\Psi_{Q_0[\omega_1],Q_1[\omega_1]}(x)\in d(\varepsilon)$.
	\end{proposition}
	
	\subsection{Amalgamation lemmas}\label{subsection-deco-amal-lemmas} We will show that $\mathbb{M}^{dec}(\mathcal{S},\mathcal{L})$ satisfies the exact same amalgamation lemmas as $\mathbb{M}(\mathcal{S},\mathcal{L})$. The proofs will be built on the results that were obtained for $\mathbb{M}(\mathcal{S},\mathcal{L})$, and we will only have to focus on defining the right decoration for each $(\mathcal{S},\mathcal{L})$-symmetric system that was constructed in Subsection \ref{subsection-amalgamation-lemmas}. These results will be used in the subsequent subsection to prove that $\mathbb{M}^{dec}(\mathcal{S},\mathcal{L})$ preserves all cardinals and $2^{\aleph_1}=\aleph_2$.
	
	\begin{remark}
		Let $(\mathcal{M},d)\in\mathbb{M}^{dec}(\mathcal{S},\mathcal{L})$ and $Q\in\mathcal{S}\cup\mathcal{L}$. Note that if $(\mathcal{M},d)\subseteq Q$, then $\mathcal{M},d\in Q$. Therefore, $\mathcal{M},d\subseteq Q$ and, in particular, $d(\varepsilon)\in Q$ and $d(\varepsilon)\subseteq Q$, for all $\varepsilon\in\dom(\mathcal{M})\subseteq Q$.
	\end{remark}
	
	\begin{lemma}\label{deco-ontop}
		Let $(\mathcal{M},d)\in\mathbb{M}^{dec}(\mathcal{S},\mathcal{L})$ and let $Q\in\mathcal{S}\cup\mathcal{L}$ be such that $(\mathcal{M},d)\subseteq Q$. Then there is a condition $(\mathcal{M}_Q,d_Q)\in\mathbb{M}^{dec}(\mathcal{S},\mathcal{L})$ such that $Q\in\mathcal{M}_Q$ and $(\mathcal{M}_Q,d_Q)\leq(\mathcal{M},d)$.
	\end{lemma}
	\begin{proof}
		Let $\mathcal{M}_Q$ be the $(\mathcal{S},\mathcal{L})$-symmetric system given by Lemma \ref{amal-ontop}. If $Q\in\mathcal{L}$, we let $d_Q$ be defined exactly like $d$, except that the domain of $d_Q$ is now $\dom(\mathcal{M}_Q)=\dom(\mathcal{M})\cup\{\varepsilon_Q\}$, and we let $d_Q(\varepsilon_Q)=\emptyset$. Checking that $d_Q$ is a decoration on $\mathcal{M}_Q$ is an easy exercise.
		
		If $Q\in\mathcal{S}$, recall that $\mathcal{M}_Q$ was defined as the set
		\[
		\mathcal{M}\cup\{Q\}\cup\{X\cap Q:X\in\mathcal{M}\cap\mathcal{L}\}.
		\]
		Therefore, $\dom(\mathcal{M}_Q)$ equals
		\[
		\dom(\mathcal{M})\cup\{\varepsilon_Q\}\cup\{\varepsilon_{X\cap Q}:X\in\mathcal{M}\cap\mathcal{L}\}.
		\]
		Moreover, since $\mathcal{M}\subseteq Q$, it is not too hard to see that if $X\in\mathcal{M}\cap\mathcal{L}$ and $P\in\mathcal{M}\cap X$, then $P\in X\cap Q$. Hence, $\varepsilon_{X\cap Q}$ is the immediate predecessor of $\varepsilon_X$ in $\dom(\mathcal{M}_Q)$. 
		
		Let $d_Q$ be the function on $\dom(\mathcal{M}_Q)$ defined as follows:
		\begin{itemize}
			\item $d_Q(\varepsilon)=d(\varepsilon)$, if $\varepsilon\in\dom(\mathcal{M})$.
			
			\item $d_Q(\varepsilon)=\emptyset$, if $\varepsilon\in\dom(\mathcal{M}_Q)\setminus\dom(\mathcal{M})$.
		\end{itemize} 
		Let us show that $d_Q$ is a decoration on $\mathcal{M}_Q$. Suppose that $\varepsilon_0$ is the immediate predecessor of $\varepsilon_1$ in $\dom(\mathcal{M}_Q)$. For every $x\in d_Q(\varepsilon_0)$, we need to find a model $P\in\mathcal{M}_Q(\varepsilon_1)$ such that $x\in P$ to verify clause (2.a) of Definition \ref{def-deco}. If $\varepsilon_0,\varepsilon_1\in\dom(\mathcal{M})$, the conclusion follows from the fact that $d$ is a decoration on $\mathcal{M}$. Suppose that $\varepsilon_0\in\dom(\mathcal{M})$ and $\varepsilon_1\in\dom(\mathcal{M}_Q)\setminus\dom(\mathcal{M})$. Then, either $\varepsilon_1=\varepsilon_{X\cap Q}$ for some $X\in\mathcal{M}\cap\mathcal{L}$, or $\varepsilon_1=\varepsilon_Q$. Hence, in both cases the conclusion follows from $(\mathcal{M},d)\subseteq Q$ and the fact that $d$ is a decoration on $\mathcal{M}$. Lastly, suppose that $\varepsilon_0\in\dom(\mathcal{M}_Q)\setminus\mathcal{M}$. Then the conclusion follows vacuously because $d_Q(\varepsilon_0)=\emptyset$.
		
		Condition (2.b) from the definition of $\mathbb{M}^{dec}(\mathcal{S},\mathcal{L})$ for $d_Q$ follows from $(\mathcal{M},d)\subseteq Q$ and the fact that, by Proposition \ref{prop9}, for every $X_0,X_1\in\mathcal{M}\cap\mathcal{L}$ such that $\varepsilon_{X_0}=\varepsilon_{X_1}$, the $\omega_1$-isomorphism $(X_0\cap Q)[\omega_1]\cong(X_1\cap Q)[\omega_1]$ is witnessed by $\Psi_{X_0,X_1}$.
		
		Checking that $(\mathcal{M}_Q,d_Q)$ extends $(\mathcal{M},d)$ is straightforward.
	\end{proof}
	
	\begin{definition}\label{deco-def-restr}
		Let $(\mathcal{M},d)\in\mathbb{M}^{dec}(\mathcal{S},\mathcal{L})$ and $Q\in\mathcal{M}\cup\mathcal{M}[\omega_1]$. Define $(\mathcal{M},d)\restr Q$ as the pair $(\mathcal{M}\cap Q,d\restr Q)$, where $d\restr Q$ is the function on $\dom(\mathcal{M}\cap Q)$ defined by $d\restr Q(\varepsilon)=d(\varepsilon)\cap Q$, for every $\varepsilon\in\dom(\mathcal{M}\cap Q)$.
	\end{definition}
	
	\begin{lemma}\label{deco-restr}
		Let $(\mathcal{M},d)\in\mathbb{M}^{dec}(\mathcal{S},\mathcal{L})$, and let $Q$ be a model in $\mathcal{M}\cup\mathcal{M}[\omega_1]$. Then $(\mathcal{M},d)\restr Q\in\mathbb{M}^{dec}(\mathcal{S},\mathcal{L})\cap Q$ and $(\mathcal{M},d)\leq(\mathcal{M},d)\restr Q$.
	\end{lemma}
	\begin{proof}
		We know from Lemma \ref{amal-rest} that $\mathcal{M}\cap Q$ is an $(\mathcal{S},\mathcal{L})$-symmetric system, is a member of $Q$, and $\mathcal{M}\supseteq\mathcal{M}\cap Q$. So, it is clear that $(\mathcal{M},d)\restr Q\in Q$. Checking that $(\mathcal{M},d)\leq(\mathcal{M},d)\restr Q$ is straightforward. Hence, we only need to show that $d\restr Q$ is a decoration on $\mathcal{M}\cap Q$.
		
		Let us check condition (2.a) first. Suppose that $\varepsilon_0$ is the immediate predecessor of $\varepsilon_1$ in $\dom(\mathcal{M}\cap Q)$, and let $x\in d\restr Q(\varepsilon_0)=d(\varepsilon_0)\cap Q$. We need to find some $P\in\mathcal{M}(\varepsilon_1)\cap Q$ such that $x\in P$. We divide the proof in two cases.
		
		\textbf{Case 1.} Suppose first that $\varepsilon_1$ is the immediate successor of $\varepsilon_0$ in $\dom(\mathcal{M})$. Since $d$ is a decoration on $\mathcal{M}$, there is some $R\in\mathcal{M}(\varepsilon_1)$ such that $x\in R$. Find, by the shoulder axiom, some $Q'\in\mathcal{M}(\varepsilon_Q)$ such that $R\in Q'[\omega_1]$, and let $P:=\Psi_{Q'[\omega_1],Q[\omega_1]}(R)$. It follows easily from the properties of $(\mathcal{S},\mathcal{L})$-symmetric systems that $x\in P\in\mathcal{M}(\varepsilon_1)\cap Q[\omega_1]$. If $Q$ is a large model, $Q\in\mathcal{M}[\omega_1]$, or $Q$ is small and there are no large models $X\in\mathcal{M}\restr(\varepsilon_1,\varepsilon_Q)$, then $P\in Q$ (in the third case by Proposition \ref{prop18}), and we are done. Hence, suppose that $Q$ is small and there are large models $X\in\mathcal{M}\restr(\varepsilon_1,\varepsilon_Q)$. Denote $P$ by $P_0$.
		
		\begin{remark}
			By Proposition \ref{prop13}, there is a large model $X_0\in \mathcal{M}\cap Q[\omega_1]$ of maximal height and such that $P_0\in X_0$. By the maximality of the height of $X_0$ and Proposition \ref{prop18}, $X_0$ is in fact a member of $Q$, and hence $X_0\cap Q\in\mathcal{M}$. Since $\varepsilon_{P_0}=\varepsilon_1\in\dom(\mathcal{M}\cap Q)\cap X_0$, we must have $\varepsilon_1<\varepsilon_{X_0\cap Q}$. By another application of Proposition \ref{prop13}, there is $N_0\in\mathcal{M}(\varepsilon_{X_0\cap Q})\cap X_0$ such that $P_0\in N_0[\omega_1]$. Let $P_1:=\Psi_{N_0[\omega_1],(X_0\cap Q)[\omega_1]}(P_0)$. Since $x\in X_0\cap Q$, it follows from the properties of $(\mathcal{S},\mathcal{L})$-symmetric systems that $x\in P_1$ and $P_1\in \mathcal{M}(\varepsilon_1)\cap (X_0\cap Q)[\omega_1]$.
		\end{remark}
		
		Recursively apply the above remark to obtain $X_0,\dots,X_n$, $N_0,\dots,N_n$, $P_0,\dots,P_n,P_{n+1}$, all in $\mathcal{M}$, such that for all $i\leq n<\omega$:
		\begin{itemize}
			\item $N_i\in\mathcal{S}$ and $X_i\in\mathcal{L}\cap Q$.
			
			\item If $0<i$, then $X_{i}\in X_{i-1}$.
			
			\item $\varepsilon_1=\varepsilon_{P_i}<\varepsilon_{X_i\cap Q}=\varepsilon_{N_i}$ and $x\in P_i\in X_i$.
			
			\item $N_i\in X_i$, $P_i\in N_i[\omega_1]$ and $P_{i+1}=\Psi_{N_i[\omega_1],(X_i\cap Q)[\omega_1]}(P_i)$.
			
			\item $\mathcal{L}\cap(\mathcal{M}\restr(\varepsilon_1,\varepsilon_{X_n\cap Q}))=\emptyset$.
		\end{itemize}
		Thus $\varepsilon_1=\varepsilon_{P_n}=\varepsilon_{P_{n+1}}$ and $x\in P_{n+1}\in X_n\cap Q$ by Proposition \ref{prop18}, and hence $P_{n+1}$ witnesses the conclusion of (2.a).

			
			
			
			

		
		\textbf{Case 2.} Suppose now that $\varepsilon_2<\varepsilon_1$ is the immediate successor of $\varepsilon_0$ in $\dom(\mathcal{M})$. Note that in this case $Q$ must be a small model, and by Proposition \ref{prop-res-gaps}, we must have $\varepsilon_0<\varepsilon_2=\varepsilon_{X\cap Q}<\varepsilon_X=\varepsilon_1$, for some large model $X\in\mathcal{M}\cap Q$. Since $x\in d\restr Q(\varepsilon_0)=d(\varepsilon_0)\cap Q$ and $d$ is a decoration on $\mathcal{M}$, there must be some small model $N\in\mathcal{M}(\varepsilon_2)$ such that $x\in N$. Let $R\in\mathcal{M}(\varepsilon_1)$ and $Q'\in\mathcal{M}(\varepsilon_Q)$ be such that $N\in R\in Q'[\omega_1]$, given by two applications of the shoulder axiom. Since $R$ is a large model, $x\in N\subseteq R$. Let $P:=\Psi_{Q'[\omega_1],Q[\omega_1]}(R)$. It follows from the properties of $(\mathcal{S},\mathcal{L})$-symmetric systems that $x\in P\in\mathcal{M}(\varepsilon_1)\cap Q[\omega_1]$. The rest of the proof is exactly the same as in the first case.

		Lastly, we check condition (2.b). Suppose that $\varepsilon\in\dom(\mathcal{M}\cap Q)$, and let $x\in d\restr Q(\varepsilon)$ and $P_0,P_1\in\mathcal{M}\cap Q$ such that $x\in P_0$ and $\varepsilon_{P_0}=\varepsilon_{P_1}$. We need to show that $\Psi_{P_0[\omega_1],P_1[\omega_1]}(x)\in d\restr Q(\varepsilon)=d(\varepsilon)\cap Q$. On the one hand, note that $\Psi_{P_0[\omega_1],P_1[\omega_1]}(x)\in d(\varepsilon)$, because $d$ is a decoration on $\mathcal{M}$. On the other hand, note that since $x,P_0,P_1\in Q$, then $\Psi_{P_0[\omega_1],P_1[\omega_1]}(x)\in Q$. Hence, $\Psi_{P_0[\omega_1],P_1[\omega_1]}(x)$ is a member of $d\restr Q(\varepsilon)$ as we wanted.
	\end{proof}
	
	\begin{notation}
		If $(\mathcal{M}_0,d_0),(\mathcal{M}_1,d_1)\in\mathbb{M}^{dec}(\mathcal{S},\mathcal{L})$, we will write $$(\mathcal{M}_0,d_0)\cong(\mathcal{M}_1,d_1)$$ to denote that there is a map $\Psi$ witnessing $\mathcal{M}_0\cong\mathcal{M}_1$, such that for every $\varepsilon\in\dom(\mathcal{M}_0)=\dom(\mathcal{M}_1)$ with $\varepsilon<\max\dom(\mathcal{M}_0)$, $d_1(\varepsilon)=\Psi(d_0(\varepsilon))$.
	\end{notation}
	
	\begin{lemma}\label{deco-amal-iso}
		Let $n<\omega$ and suppose that $(\mathcal{M}_0,d_0),\dots,(\mathcal{M}_n,d_n)$ are members of $\mathbb{M}^{dec}(\mathcal{S},\mathcal{L})$ such that $(\mathcal{M}_i,d_i)\cong(\mathcal{M}_j,d_j)$, for all $i,j\leq n$. Then there is $(\mathcal{M},d)\in\mathbb{M}^{dec}(\mathcal{S},\mathcal{L})$ such that $(\mathcal{M},d)\leq(\mathcal{M}_i,d_i)$ for every $i\leq n$.
	\end{lemma}
	\begin{proof}
		For any two $i,j\leq n$, we will denote the isomorphism witnessing $(\mathcal{M}_i,d_i)\cong(\mathcal{M}_j,d_j)$ by $\Psi_{i,j}$, and the intersection $\bigcup\mathcal{M}_i[\omega_1]\cap\bigcup\mathcal{M}_j[\omega_1]$ by $\Delta_{i,j}$. We know from Lemma \ref{pureamalgamation2} that $\bigcup_{i\leq n}\mathcal{M}_i$ is an $(\mathcal{S},\mathcal{L})$-symmetric system. Denote $\bigcup_{i\leq n}\mathcal{M}_i$ by $\mathcal{M}$, and note that $\dom(\mathcal{M})=\dom(\mathcal{M}_i)$ for every $i\leq n$.
		
		Define the function $d$ on $\dom(\mathcal{M})$ by letting $d(\varepsilon)=\bigcup_{i\leq n}d_i(\varepsilon)$ for every $\varepsilon\in\dom(\mathcal{M})$. Let us check that $d$ is a decoration on $\mathcal{M}$. To see that condition (2.a) is satisfied, let $\varepsilon_0$ be the immediate predecessor of $\varepsilon_1$ in $\dom(\mathcal{M})$ and let $x\in d(\varepsilon_0)$. Then, there must be some $i\leq n$ such that $x\in d_i(\varepsilon_0)$, and hence there needs to exist a model $Q\in\mathcal{M}_i(\varepsilon_1)$ such that $x\in Q$. 
		
		Let us now check condition (2.b). Suppose that $\varepsilon\in\dom(\mathcal{M})$, and let $x\in d(\varepsilon)$ and $Q,P\in\mathcal{M}$ be such that $x\in Q$ and $\varepsilon_Q=\varepsilon_P$. We need to check that $\Psi_{Q[\omega_1],P[\omega_1]}(x)\in d(\varepsilon)$. Let $i,j,k\leq n$ such that $x\in d_i(\varepsilon)$, $Q\in\mathcal{M}_j$, and $P\in\mathcal{M}_k$. Since $(\mathcal{M}_i,d_i)\cong(\mathcal{M}_j,d_j)$, we have $d_j(\varepsilon)=\Psi_{i,j}(d_i(\varepsilon))$. Hence, $\Psi_{i,j}(x)\in d_j(\varepsilon)$. However, note that since $x\in Q\in\mathcal{M}_j$, $x$ must be a member of $\Delta_{i,j}$, so $x$ is fixed by the isomorphism $\Psi_{i,j}$, and thus $x\in d_j(\varepsilon)$. Let us denote $Q':=\Psi_{j,k}(Q)$ and $x':=\Psi_{j,k}(x)$, and note that as $(\mathcal{M}_j,d_j)\cong(\mathcal{M}_k,d_k)$, we have $Q'\in\mathcal{M}_k$ and $x'\in d_k(\varepsilon)\cap Q'$. Now, since $d_k$ is a decoration on $\mathcal{M}_k$, $\Psi_{Q'[\omega_1],P[\omega_1]}(x')\in d_k(\varepsilon)$. But then  we are done, because 
		\begin{align*}
			\Psi_{Q[\omega_1],P[\omega_1]}(x)&=\Psi_{Q'[\omega_1],P[\omega_1]}(\Psi_{j,k}(x))\\
			&=\Psi_{Q'[\omega_1],P[\omega_1]}(x')\in d_k(\varepsilon)\subseteq d(\varepsilon).
		\end{align*}
		So, we can conclude that $d$ is a decoration on $\mathcal{M}$.	
		
		Lastly, it is obvious from the definition of $d$ that $(\mathcal{M},d)$ extends all conditions $(\mathcal{M}_i,d_i)$.
	\end{proof}
	
	Recall from Section \ref{section-pure} that if $\mathcal{M},\mathcal{W}$ are two $(\mathcal{S},\mathcal{L})$-symmetric systems, $\mathcal{W}\land\mathcal{M}$ denotes the set of all models of the form $X\cap M$, where $M\in\mathcal{M}\cap\mathcal{S}$ and $X\in\mathcal{W}\cap\mathcal{L}$ are such that $X\in M$.
	
	\begin{proposition}\label{prop-deco-amalgamation}
		Let $(\mathcal{M},d)\in\mathbb{M}^{dec}(\mathcal{S},\mathcal{L})$ and let $Q\in\mathcal{M}$. Let $(\mathcal{W},e)$ be a condition in $\mathbb{M}^{dec}(\mathcal{S},\mathcal{L})$ with the following properties:
		\begin{enumerate}
			\item $(\mathcal{W},e)\leq(\mathcal{M},d)\restr Q[\omega_1]$ and $(\mathcal{W},e)\subseteq Q[\omega_1]$,
			
			\item $\mathcal{W}(\max\dom(\mathcal{W}))\subseteq Q$,
			
			\item $e(\max\dom(\mathcal{W}))\subseteq Q$, and
			
			\item $\mathcal{W}\land\mathcal{M}\subseteq\mathcal{W}$.
		\end{enumerate}
		Then there is $(\mathcal{V},f)\in\mathbb{M}^{dec}(\mathcal{S},\mathcal{L})$ that extends both $(\mathcal{M},d)$ and $(\mathcal{W},e)$.
	\end{proposition}
	\begin{proof}
		Recall from Proposition \ref{prop-pure-amalgamation-3} that 
		\[
		\mathcal{V}=\mathcal{M}(\geq\varepsilon_Q)\cup\{\Psi_{Q[\omega_1],Q'[\omega_1]}(W):W\in\mathcal{W},Q'\in\mathcal{M}(\varepsilon_Q)\}
		\]
		is an $(\mathcal{S},\mathcal{L})$-symmetric system extending both $\mathcal{M}$ and $\mathcal{W}$.
		
		Define the function $f$ on $\dom(\mathcal{V})$ by letting 
		\begin{itemize}
			\item $f(\varepsilon)=\bigcup\{\Psi_{Q[\omega_1],Q'[\omega_1]}(e(\varepsilon)):Q'\in\mathcal{M}(\varepsilon_Q)\}$, for every $\varepsilon\in\dom(\mathcal{V})$ such that $\varepsilon<\varepsilon_Q$, and
			
			\item $f(\varepsilon)=d(\varepsilon)$, for every $\varepsilon\in\dom(\mathcal{V})$ such that $\varepsilon\geq\varepsilon_Q$.
		\end{itemize}
		
		Let us show that it is a decoration on $\mathcal{V}$. Let us start by showing condition (2.a) from Definition \ref{def-deco}. Let $\varepsilon_0$ be the immediate predecessor of $\varepsilon_1$ in $\dom(\mathcal{V})$ and let $x\in f(\varepsilon_0)$. We have to find a model $P\in\mathcal{V}(\varepsilon_1)$ such that $x\in P$. We divide the proof in two cases.
		
		\textbf{Case 1.} Suppose first that $\varepsilon_0<\varepsilon_Q$, and note that $\varepsilon_1\leq\varepsilon_Q$. By definition of $f$ there must be a model $Q'\in\mathcal{M}(\varepsilon_Q)$ such that $x\in\Psi_{Q[\omega_1],Q'[\omega_1]}(e(\varepsilon_0))$. Hence, there is some $y\in e(\varepsilon_0)$ such that $\Psi_{Q[\omega_1],Q'[\omega_1]}(y)=x$. If $\varepsilon_1<\varepsilon_Q$, as $\dom(\mathcal{V}(<\varepsilon_Q))=\dom(\mathcal{W})$, $\varepsilon_0$ is also the immediate predecessor of $\varepsilon_1$ in $\dom(\mathcal{W})$. Hence, as $e$ is a decoration on $\mathcal{W}$, there must be a model $W\in\mathcal{W}(\varepsilon_1)$ such that $y\in W$. Therefore, $x\in\Psi_{Q[\omega_1],Q'[\omega_1]}(W)\in\mathcal{V}(\varepsilon_1)$. If $\varepsilon_1=\varepsilon_Q$, as $e(\max\dom(\mathcal{W}))=e(\varepsilon_0)\subseteq Q$, we have $y\in e(\varepsilon_0)\subseteq Q$, and thus, $x\in Q'\in\mathcal{V}(\varepsilon_1)$. 
		
		\textbf{Case 2.} Now, suppose that $\varepsilon_0\geq\varepsilon_Q$. By definition of $f$, we have $f(\varepsilon_0)=d(\varepsilon_0)$, and as $\mathcal{V}(\geq\varepsilon_Q)=\mathcal{M}(\geq\varepsilon_Q)$ and $d$ is a decoration on $\mathcal{M}$, there must be some $P\in\mathcal{M}(\varepsilon_1)$ such that $x\in P$.
		
		Let us now check condition (2.b). Let $\varepsilon\in\dom(\mathcal{V})$ and $x\in f(\varepsilon)$, and suppose that $P_0,P_1\in\mathcal{V}$ are such that $x\in P_0$ and $\varepsilon<\varepsilon_{P_0}=\varepsilon_{P_1}$. We need to check that $\Psi_{P_0[\omega_1],P_1[\omega_1]}(x)\in f(\varepsilon)$. We divide the proof in four cases.
		
		\textbf{Case 1.} Suppose first that $\varepsilon<\varepsilon_{P_0}=\varepsilon_{P_1}<\varepsilon_Q$. By definition of $\mathcal{V}$ and $f$, there are $Q',Q_0,Q_1\in\mathcal{M}(\varepsilon_Q)$ such that
		\begin{enumerate}
			\item $x\in\Psi_{Q[\omega_1],Q'[\omega_1]}(e(\varepsilon))$,
			\item $P_0\in\Psi_{Q[\omega_1],Q_0[\omega_1]}(\mathcal{W})$, and
			\item $P_1\in\Psi_{Q[\omega_1],Q_1[\omega_1]}(\mathcal{W})$.
		\end{enumerate}
		Hence, there are $y\in e(\varepsilon)$ and $W_0,W_1\in\mathcal{W}$ such that
		\begin{enumerate}[label=(\arabic*')]
			\item $\Psi_{Q[\omega_1],Q'[\omega_1]}(y)=x$,
			\item $\Psi_{Q[\omega_1],Q_0[\omega_1]}(W_0)=P_0$, and
			\item $\Psi_{Q[\omega_1],Q_1[\omega_1]}(W_1)=P_1$.
		\end{enumerate}
		Note that as $x\in\Psi_{Q[\omega_1],Q'[\omega_1]}(e(\varepsilon))\cap P_0$, $x$ is also a member of $Q'[\omega_1]\cap Q_0[\omega_1]$, and hence
		\begin{align*}
			x=\Psi_{Q'[\omega_1],Q_0[\omega_1]}(x)\in&\hbox{ }\Psi_{Q'[\omega_1],Q_0[\omega_1]}(\Psi_{Q[\omega_1],Q'[\omega_1]}(e(\varepsilon)))\\
			&=\Psi_{Q[\omega_1],Q_0[\omega_1]}(e(\varepsilon)),
		\end{align*} 
		and 
		\begin{align*}
			x=\Psi_{Q'[\omega_1],Q_0[\omega_1]}(x)&=\Psi_{Q'[\omega_1],Q_0[\omega_1]}(\Psi_{Q[\omega_1],Q'[\omega_1]}(y))\\
			&=\Psi_{Q[\omega_1],Q_0[\omega_1]}(y).
		\end{align*}
		Therefore, $x\in \Psi_{Q[\omega_1],Q_0[\omega_1]}(e(\varepsilon))\cap P_0$, and by symmetry, $y\in e(\varepsilon)\cap W_0$. Hence, as $e$ is a decoration on $\mathcal{W}$, we have $\Psi_{W_0[\omega_1],W_1[\omega_1]}(y)\in e(\varepsilon)$, and we can conclude that  
		\[
		\Psi_{Q[\omega_1],Q_1[\omega_1]}(\Psi_{W_0[\omega_1],W_1[\omega_1]}(y))\in\Psi_{Q[\omega_1],Q_1[\omega_1]}(e(\varepsilon))\subseteq f(\varepsilon).
		\]
		But then we are done as $\Psi_{Q[\omega_1],Q_0[\omega_1]}(W_0)=P_0$ and $\Psi_{Q[\omega_1],Q_1[\omega_1]}(W_1)=P_1$ imply that
		\begin{align*}
			\Psi_{Q[\omega_1],Q_1[\omega_1]}(\Psi_{W_0[\omega_1],W_1[\omega_1]}(y))&=\Psi_{Q[\omega_1],Q_1[\omega_1]}(\Psi_{W_0[\omega_1],W_1[\omega_1]}(\Psi_{Q_0[\omega_1],Q[\omega_1]}(x)))\\
			&=\Psi_{P_0[\omega_1],P_1[\omega_1]}(x).
		\end{align*}
		
		\textbf{Case 2.} If $\varepsilon<\varepsilon_{P_0}=\varepsilon_{P_1}=\varepsilon_Q$, the argument is simpler. Suppose that $Q'\in\mathcal{M}(\varepsilon_Q)$ is such that $x\in \Psi_{Q[\omega_1],Q'[\omega_1]}(e(\varepsilon))$. 
		Then, as $x\in Q'[\omega_1]\cap P_0$, we have that $x=\Psi_{Q'[\omega_1],P_0[\omega_1]}(x)$ is a member of 
		\[
		\Psi_{Q'[\omega_1],P_0[\omega_1]}(\Psi_{Q[\omega_1],Q'[\omega_1]}(e(\varepsilon)))=\Psi_{Q[\omega_1],P_0[\omega_1]}(e(\varepsilon)).
		\]
		Therefore, 
		\begin{align*}
			\Psi_{P_0[\omega_1],P_1[\omega_1]}(x)\in &\hbox{ }\Psi_{P_0[\omega_1],P_1[\omega_1]}(\Psi_{Q[\omega_1],P_0[\omega_1]}(e(\varepsilon)))\\
			&=\Psi_{Q[\omega_1],P_1[\omega_1]}(e(\varepsilon))\subseteq f(\varepsilon).
		\end{align*}
		
		\textbf{Case 3.} Suppose now that $\varepsilon<\varepsilon_Q<\varepsilon_{P_0}=\varepsilon_{P_1}$. Let $Q'\in\mathcal{M}(\varepsilon_Q)$ such that $x\in \Psi_{Q[\omega_1],Q'[\omega_1]}(e(\varepsilon))$. Using the shoulder axiom for $\mathcal{M}$, find $P_2\in\mathcal{M}(\varepsilon_{P_0})$ such that $Q'\in P_2[\omega_1]$, and let $Q_0:=\Psi_{P_2[\omega_1],P_0[\omega_1]}(Q')$, which is an element of $\mathcal{M}$ by Proposition \ref{prop6}. Note that $x\in P_2[\omega_1]\cap P_0$, and hence $x=\Psi_{P_2[\omega_1],P_0[\omega_1]}(x)\in Q_0[\omega_1]$. So, if we let $Q_1:=\Psi_{P_0[\omega_1],P_1[\omega_1]}(Q_0)$, which is again an element of $\mathcal{M}$, we can argue as in the last case to show that
		\[
		\Psi_{P_0[\omega_1],P_1[\omega_1]}(x)=\Psi_{Q_0[\omega_1],Q_1[\omega_1]}(x)\in \Psi_{Q[\omega_1],Q_1[\omega_1]}(e(\varepsilon))\subseteq f(\varepsilon).
		\]
		
		\textbf{Case 4.} Lastly, if $\varepsilon_Q\leq\varepsilon<\varepsilon_{P_0}=\varepsilon_{P_1}$, the conclusion follows from the definition of $f$ and the fact that $d$ is a decoration on $\mathcal{M}$.
		
		In order to show that $(\mathcal{V},f)$ extends both $(\mathcal{M},d)$ and $(\mathcal{W},e)$, it is enough to note that, by definition, $(\mathcal{V},f)$ strengthens $(\mathcal{W},e)$, and that as $(\mathcal{W},e)$ extends $(\mathcal{M},d)\restr Q[\omega_1]$, then $(\mathcal{V}(<\varepsilon_Q),f\restr\varepsilon_Q)\leq(\mathcal{M}(<\varepsilon_Q),d\restr\varepsilon_Q)$.
	\end{proof}
	
	\begin{lemma}\label{deco-amal-L}
		Let $(\mathcal{M},d)\in\mathbb{M}^{dec}(\mathcal{S},\mathcal{L})$ and let $X\in\mathcal{M}\cap\mathcal{L}$. Let $(\mathcal{W},e)$ be a condition in $\mathbb{M}^{dec}(\mathcal{S},\mathcal{L})\cap X$ such that $(\mathcal{W},e)\leq(\mathcal{M},d)\restr X$. Then there is $(\mathcal{V},f)\in\mathbb{M}^{dec}(\mathcal{S},\mathcal{L})$ that extends both $(\mathcal{M},d)$ and $(\mathcal{W},e)$.
	\end{lemma}
	\begin{proof}
		By Proposition \ref{prop-pure-amalgamation-1}, we have $\mathcal{W}\land\mathcal{M}\subseteq\mathcal{W}$. Note also that, by $(\mathcal{W},e)\in X$, we have $\mathcal{W}(\max\dom(\mathcal{W}))\subseteq X$ and $e(\max\dom(\mathcal{W}))\subseteq X$. Hence, by Proposition \ref{prop-deco-amalgamation} there is a condition $(\mathcal{V},f)\in\mathbb{M}^{dec}(\mathcal{S},\mathcal{L})$ extending both $(\mathcal{M},d)$ and $(\mathcal{W},e)$.
	\end{proof}

	\begin{lemma}\label{deco-amal-S}
		Let $(\mathcal{M},d)\in\mathbb{M}^{dec}(\mathcal{S},\mathcal{L})$ and let $M\in\mathcal{M}\cap\mathcal{S}$. Let $(\mathcal{W},e)$ be a condition in $\mathbb{M}^{dec}(\mathcal{S},\mathcal{L})\cap M$ such that $(\mathcal{W},e)\leq(\mathcal{M},d)\restr M$. Then there is $(\mathcal{U},g)\in\mathbb{M}^{dec}(\mathcal{S},\mathcal{L})$ extending both $(\mathcal{M},d)$ and $(\mathcal{W},e)$.
	\end{lemma}
	\begin{proof}
		The argument will be built on the proof of Lemma \ref{lemma-amal-S}, which is the analogous amalgamation lemma for the pure side condition forcing, $\mathbb{M}(\mathcal{S},\mathcal{L})$. Let $X_n$ be a large model containing both $(\mathcal{M},d)$ and $(\mathcal{W},e)$, and fix a maximal $\in$-chain $X_0\in\dots\in X_{n-1}$ of large models in $\mathcal{M}\cap M$. For every $i\leq n$ denote
		\begin{itemize}
			\item $(\mathcal{M}_i,d_i)=(\mathcal{M},d)\restr X_i$,
			\item $(\mathcal{W}_i,e_i)=(\mathcal{W},e)\restr X_i$, and
			\item $M_i=X_i\cap M$.
		\end{itemize} 
		Observe that $(\mathcal{M}_n,d_n)=(\mathcal{M},d)$, $(\mathcal{W}_n,e_n)=(\mathcal{W},e)$ and $M_n=M$. Moreover, note that $(\mathcal{M}_i,d_i),(\mathcal{W}_i,e_i)\in\mathbb{M}^{dec}(\mathcal{S},\mathcal{L})$ by Lemma \ref{deco-restr}, and that $(\mathcal{W}_i,e_i)\leq(\mathcal{M}_i,d_i)\restr M_i$ and $(\mathcal{W}_i,e_i)\in M_i$. We will build, by induction on $i\leq n$, a decreasing sequence of $\mathbb{M}^{dec}(\mathcal{S},\mathcal{L})$-conditions $(\mathcal{U}_i,g_i)$ contained in $X_i$ and extending both $(\mathcal{M}_i,d_i)$ and $(\mathcal{W}_i,e_i)$. For all $i\leq n$, let $\mathcal{V}_i\subseteq\mathcal{V}_i^*\subseteq\mathcal{U}_i$ be the $(\mathcal{S},\mathcal{L})$-symmetric systems built in the proof of Lemma \ref{lemma-amal-S}. We will define decorations $f_i,f_i^*,g_i$ on $\mathcal{V}_i,\mathcal{V}_i^*,\mathcal{U}_i$, respectively, so that $(\mathcal{U}_i,g_i)\leq(\mathcal{V}_i^*,f_i^*)\leq(\mathcal{V}_i,f_i)$, and such that $(\mathcal{U}_i,g_i)$ witnesses the compatibility of $(\mathcal{M}_i,d_i)$ and $(\mathcal{W}_i,e_i)$ in $X_i$. 
		
		Let $i<n$ and let $\varepsilon_i^-$ denote either $\varepsilon_{X_{i-1}}$, in case $i>0$, or an arbitrary ordinal smaller than the minimum of $\dom(\mathcal{W})$, otherwise. Suppose that $(\mathcal{V}_i,f_i)$ is a condition in $\mathbb{M}^{dec}(\mathcal{S},\mathcal{L})$ with the following properties:
		\begin{enumerate}[label=(\alph*)]
			\item $(\mathcal{V}_i,f_i)\leq(\mathcal{M}_i,d_i)\restr M_i[\omega_1]$ and $(\mathcal{V}_i,f_i)\subseteq M_i[\omega_1]$.
			
			\item $(\mathcal{V}_i,f_i)\leq(\mathcal{W}_i,e_i)$.
			
			\item $\mathcal{V}_i(\geq\varepsilon_i^-)=\mathcal{W}_i(\geq\varepsilon_i^-)$, and $f_i(\varepsilon)=e_i(\varepsilon)$ for all $\varepsilon\in\dom(\mathcal{W}_i(\geq\varepsilon_i^-))$.
		\end{enumerate}
		If $i=0$, we simply let $(\mathcal{V}_0,f_0)=(\mathcal{W}_0,e_0)$.
		
		Let $\mathcal{V}_i^*$ be the $(\mathcal{S},\mathcal{L})$-symmetric system extending $\mathcal{V}_i$, defined exactly as in the proof of Lemma \ref{lemma-amal-S}. That is, $\mathcal{V}_i^*$ is the union of $\mathcal{V}_i$ and all the sets $F_W$, for $W\in\mathcal{W}_i(>\varepsilon_i^-)\cap\mathcal{L}$. Recall that $F_W$ is the set $\{W\cap N:N\in E_i\}$, where $E_i$ is an $\in$-chain of small models $N_0\in\dots\in N_k$ of maximal length of $\mathcal{M}_i\restr[\varepsilon_{M_i},\varepsilon_i^+)$ such that $N_0=M_i$, and $\varepsilon_i^+$ is the least $\omega_2$-height of any large model in $\mathcal{M}_i(>\varepsilon_{M_i})$. Recall that $\mathcal{V}_i^*$ coincides with the union of $\mathcal{V}_i$ and
		\[
		\mathcal{V}_i(>\varepsilon_i^-)\land\mathcal{M}_i\restr[\varepsilon_{M_i},\varepsilon_i^+).
		\]
		
		
		
		Let $f_i^*$ be the function on $\dom(\mathcal{V}_i^*)$ defined as follows:
		\begin{itemize}
			\item $f_i^*(\varepsilon)=f_i(\varepsilon)$, for every $\varepsilon\in\dom(\mathcal{V}_i)$.
			
			\item $f_i^*(\varepsilon)=\emptyset$, for every $\varepsilon\in\dom(\mathcal{V}_i^*)\setminus\dom(\mathcal{V}_i)$
		\end{itemize}
		Note that if $\varepsilon\in\dom(\mathcal{V}_i^*)\setminus\dom(\mathcal{V}_i)$, then there is $W\cap N\in F_W$, where $W\in\mathcal{W}_i(>\varepsilon_i^-)\cap\mathcal{L}$, such that $\varepsilon=\varepsilon_{W\cap N}$. It is easy to verify that $f_i^*$ is a decoration on $\mathcal{V}_i^*$ using the same argument from the proof of Lemma \ref{deco-ontop}. Therefore, $(\mathcal{V}_i^*,f_i^*)$ is a condition in $\mathbb{M}^{dec}(\mathcal{S},\mathcal{L})$, and it is easy to verify that it has the following properties:
		\begin{enumerate}
			\item $(\mathcal{V}_i^*,f_i^*)\leq(\mathcal{M}_i,d_i)\restr M_i[\omega_1]$ and $(\mathcal{V}_i^*,f_i^*)\subseteq M_i[\omega_1]$,
			
			\item $\mathcal{V}_i^*(\max\dom(\mathcal{V}_i^*))=\mathcal{V}_i^*(\max\dom(\mathcal{\mathcal{W}}_i))\subseteq M_i$,
			
			\item $f_i^*(\max\dom(\mathcal{V}_i^*))=f_i(\max\dom(\mathcal{\mathcal{W}}_i))\subseteq M_i$, and
			
			\item $\mathcal{V}_i^*\land\mathcal{M}_i\subseteq\mathcal{V}_i^*$.
		\end{enumerate}
		Hence, by Proposition \ref{prop-deco-amalgamation}, there is a condition $(\mathcal{U}_i,g_i)\in\mathbb{M}^{dec}(\mathcal{S},\mathcal{L})$ that extends both $(\mathcal{M}_i,d_i)$ and $(\mathcal{V}_i^*,f_i^*)$, and hence, in particular, also extends $(\mathcal{W}_i,e_i)$. Recall from the proof of Proposition \ref{prop-deco-amalgamation} that the $(\mathcal{S},\mathcal{L})$-symmetric system $\mathcal{U}_i$ is defined as 
		\[
		\mathcal{U}_i:=\mathcal{M}_i(\geq\varepsilon_{M_i})\cup\{\Psi_{M_i[\omega_1],N_i[\omega_1]}(V):V\in\mathcal{V}_i^*,N_i\in\mathcal{M}_i(\varepsilon_{M_i})\},
		\]
		and the decoration $g_i$ is defined for every $\varepsilon\in\dom(\mathcal{U}_i)$ as
		\begin{itemize}
			\item $g_i(\varepsilon)=\bigcup\{\Psi_{M_i[\omega_1],N_i[\omega_1]}(f_i^*(\varepsilon)):N_i\in\mathcal{M}_i(\varepsilon_{M_i})\}$, if $\varepsilon<\varepsilon_{M_i}$, and
			
			\item $g_i(\varepsilon)=d_i(\varepsilon)$, if $\varepsilon\geq\varepsilon_{M_i}$.
		\end{itemize}
		Therefore, $(\mathcal{U}_i,g_i)\subseteq X_i$.
		
		If $i=n$, we have obtained a condition $(\mathcal{U}_n,g_n)\in\mathbb{M}^{dec}(\mathcal{S},\mathcal{L})$ that extends $(\mathcal{M}_n,d_n)=(\mathcal{M},d)$ and $(\mathcal{W}_n,e_n)=(\mathcal{W},e)$, and we are done.
		
		If $i<n$, we have to define the condition $(\mathcal{V}_{i+1},f_{i+1})$ so that the induction can go through. Observe that $(\mathcal{U}_i,g_i)$ has the following properties:
		\begin{enumerate}
			\item $(\mathcal{U}_i,g_i)\leq(\mathcal{W}_{i+1},e_{i+1})\restr X_i=(\mathcal{W}_i,e_i)$ and $(\mathcal{U}_i,g_i)\subseteq X_i$.
			
			\item Since $(\mathcal{U}_i,g_i)\subseteq X_i$, in particular, $\mathcal{U}_i(\max\dom(\mathcal{U}_i))$ and $g_i(\max\dom(\mathcal{U}_i))$ are both subsets of $X_i$.
			
			\item $\mathcal{U}_i\land\mathcal{W}_{i+1}\subseteq\mathcal{U}_i$, by Proposition \ref{prop-pure-amalgamation-1}.
		\end{enumerate}
		Therefore, by Proposition \ref{prop-deco-amalgamation} (and Proposition \ref{prop-pure-amalgamation-3}) we can define $(\mathcal{V}_{i+1},f_{i+1})$ by letting $\mathcal{V}_{i+1}$ be the $(\mathcal{S},\mathcal{L})$-symmetric system
		\[
		\mathcal{V}_{i+1}:=\mathcal{W}_{i+1}(\geq\varepsilon_{X_i})\cup\{\Psi_{X_i,Y_i}(U):U\in\mathcal{U}_i,Y_i\in\mathcal{W}_{i+1}(\varepsilon_{X_i})\},
		\]
		and $f_{i+1}$ the decoration on $\mathcal{V}_{i+1}$ defined, for every $\varepsilon\in\dom(\mathcal{V}_{i+1})$, as
		\begin{itemize}
			\item $f_{i+1}(\varepsilon)=\bigcup\{\Psi_{X_i,Y_i}(g_i(\varepsilon)):Y_i\in\mathcal{W}_{i+1}(\varepsilon_{X_i})\}$, if $\varepsilon<\varepsilon_{X_i}$, and
			
			\item $f_{i+1}(\varepsilon)=e_{i+1}(\varepsilon)$, if $\varepsilon\geq\varepsilon_{X_i}$.
		\end{itemize}
		Let us verify that $(\mathcal{V}_{i+1},f_{i+1})$ has the properties (a)-(c) that we assumed for $(\mathcal{V}_i,f_i)$ at the beginning of the inductive construction of the conditions $(\mathcal{U}_i,g_i)$. More precisely, since the $(\mathcal{S},\mathcal{L})$-symmetric system $\mathcal{V}_{i+1}$ is defined exactly as in the proof of Lemma \ref{lemma-amal-S}, we only need to check that the following conditions are satisfied:
		\begin{enumerate}[label=(\alph*)]
			\item $f_{i+1}(\varepsilon)\supseteq (d_{i+1}\restr M_{i+1}[\omega_1])(\varepsilon)$, for every $\varepsilon\in\dom(\mathcal{M}_{i+1}\cap M_{i+1}[\omega_1])$, and $f_{i+1}\subseteq M_{i+1}[\omega_1]$.
			
			\item $f_{i+1}(\varepsilon)\supseteq e_{i+1}(\varepsilon)$, for every $\varepsilon\in\dom(\mathcal{W}_{i+1}(<\varepsilon_{X_i}))$.
			
			\item $f_{i+1}(\varepsilon)=e_{i+1}(\varepsilon)$, for every $\varepsilon\in\dom(\mathcal{W}_{i+1}(\geq\varepsilon_{X_i}))$.
		\end{enumerate}
		By the construction of $(\mathcal{V}_{i+1},f_{i+1})$ and Proposition \ref{prop-deco-amalgamation}, we already have (b) and (c). Hence, we only need to prove (a).
		
		Recall that $d_{i+1}\restr M_{i+1}[\omega_1]$ is the function on $\dom(\mathcal{M}_{i+1}\cap M_{i+1}[\omega_1])$ defined by letting $(d_{i+1}\restr M_{i+1}[\omega_1])(\varepsilon)=d_{i+1}(\varepsilon)\cap M_{i+1}[\omega_1]$, for every ordinal $\varepsilon\in\dom(\mathcal{M}_{i+1}\cap M_{i+1}[\omega_1])$. Hence, item (a) follows from the following two facts:
		\begin{itemize}
			\item Since $\mathcal{M}_{i+1}(\varepsilon_{X_i})\cap M_{i+1}[\omega_1]=\mathcal{M}_{i+1}(\varepsilon_{X_i})\cap M_{i+1}$ and $(\mathcal{U}_i,g_i)\leq(\mathcal{M}_i,d_i)$, we have
			\begin{align*}
				d_{i+1}(\varepsilon)\cap M_{i+1}[\omega_1]&=\bigcup\{\Psi_{X_i,Y_i}(d_i(\varepsilon)):Y_i\in\mathcal{M}_{i+1}(\varepsilon_{X_i})\cap M_{i+1}\}\\
				&\subseteq\bigcup\{\Psi_{X_i,Y_i}(g_i(\varepsilon)):Y_i\in\mathcal{W}_{i+1}(\varepsilon_{X_i})\} \\
				&=f_{i+1}(\varepsilon),
			\end{align*}
			for every $\varepsilon<\varepsilon_{X_i}$ in $\dom(\mathcal{M}_{i+1}\cap M_{i+1}[\omega_1])$.
			
			\item As $(\mathcal{W}_{i+1},e_{i+1})\leq(\mathcal{M}_{i+1},d_{i+1})\restr M_{i+1}$, on the one hand we have $\dom(\mathcal{M}_{i+1}(\geq\varepsilon_{X_i})\cap M_{i+1}[\omega_1])=\dom(\mathcal{M}_{i+1}(\geq\varepsilon_{X_i})\cap M_{i+1})$, and on the other hand,
			\begin{align*}
				d_{i+1}(\varepsilon)\cap M_{i+1}[\omega_1]&=d_{i+1}(\varepsilon)\cap M_{i+1}\\
				&\subseteq e_{i+1}(\varepsilon)=f_{i+1}(\varepsilon),
			\end{align*}
			for every $\varepsilon\geq\varepsilon_{X_i}$ in $\dom(\mathcal{M}_{i+1}\cap M_{i+1}[\omega_1])$. 
		\end{itemize}
		
		In summary, we have shown that the condition $(\mathcal{V}_{i+1},f_{i+1})$ has the following properties, as we wanted:
		\begin{enumerate}[label=(\alph*)]
			\item $(\mathcal{V}_{i+1},f_{i+1})\leq(\mathcal{M}_{i+1},d_{i+1})\restr M_{i+1}[\omega_1]$ and $(\mathcal{V}_{i+1},f_{i+1})\subseteq M_{i+1}[\omega_1]$.
			
			\item $(\mathcal{V}_{i+1},f_{i+1})\leq(\mathcal{W}_{i+1},e_{i+1})$.
			
			\item $\mathcal{V}_{i+1}(\geq\varepsilon_{X_i})=\mathcal{W}_{i+1}(\geq\varepsilon_{X_i})$, and $f_{i+1}(\varepsilon)=e_{i+1}(\varepsilon)$ for every ordinal $\varepsilon\in\dom(\mathcal{W}_{i+1}(\geq\varepsilon_{X_i}))$.
		\end{enumerate}
	\end{proof}
	
	\subsection{Preservation lemmas} The following results are proven exactly like the corresponding preservation theorems for $\mathbb{M}(\mathcal{S},\mathcal{L})$ from Section \ref{section-pure}, using the amalgamation lemmas from the last subsection.
	
	\begin{theorem}\label{deco-preservation}
		The forcing $\mathbb{M}^{dec}(\mathcal{S},\mathcal{L})$ is strongly $\mathcal{S}$-proper and strongly $\mathcal{L}$-proper. Moreover, if $2^{\aleph_1}=\aleph_2$, then $\mathbb{M}^{dec}(\mathcal{S},\mathcal{L})$ has the $\aleph_3$-Knaster condition.
	\end{theorem}
	
	By combining the last theorem with Lemmas \ref{preservacio-proper} and \ref{preservation-cc}, we get the following preservation theorem.
	
	\begin{corollary}\label{deco-preservation2}
		If $2^{\aleph_1}=\aleph_2$ holds and $\mathcal{L}$ is stationary in $H(\kappa)$, then $\mathbb{M}^{dec}(\mathcal{S},\mathcal{L})$ preserves all cardinals.
	\end{corollary}
	
	\begin{theorem}\label{deco-preservation3}
		If $\mathcal{L}$ is stationary in $H(\kappa)$, then $\mathbb{M}^{dec}(\mathcal{S},\mathcal{L})$ preserves $2^{\aleph_1}=\aleph_2$.
	\end{theorem}
	
	\subsection{First application: Forcing a club subset of $\omega_2$}\label{subsection-deco-club}
	
	Let $G$ be an $\mathbb{M}^{dec}(\mathcal{S},\mathcal{L})$-generic filter over $V$. Let
	\[
	\mathcal{M}_G^{dec}:=\{Q\in\mathcal{M}:(\mathcal{M},d)\in G\},
	\]
	and
	\[
	C_G^{dec}:=\{\varepsilon_Q:Q\in\mathcal{M}_G^{dec}\}=\dom(\mathcal{M}_G^{dec}).
	\]
	
	As in the case of the poset $\mathbb{M}(\mathcal{S},\mathcal{L})$, it is easily seen that $\mathcal{M}_G^{dec}$ satisfies clauses (A)-(E) from the definition of $(\mathcal{S},\mathcal{L})$-symmetric system. Moreover, as in Subsection \ref{subsection-pure-first-application}, $\mathcal{M}_G^{dec}\cap\mathcal{S}$ is always unbounded in $H(\kappa)^V$, whereas $\mathcal{M}_G^{dec}\cap\mathcal{L}$ is unbounded in $H(\kappa)^V$, only if $\mathcal{L}$ is unbounded in $H(\kappa)^V$.
	
	The argument in the proof of the following theorem is inspired by the proof of Lemma 5.2 of \cite{Neeman2014Forcingwithsequencesofmodelsoftwotypes}.
	
	\begin{theorem}
		If $\mathcal{L}$ is stationary in $H(\kappa)^V$, then $C_G^{dec}$ is a club subset of $\omega_2$.
	\end{theorem}
	\begin{proof}
		By Lemma \ref{deco-ontop} and the unboundedness of $\mathcal{S}\cup\mathcal{L}$, for every $\varepsilon<\omega_2$ the set of conditions $(\mathcal{M},d)\in\mathbb{M}^{dec}(\mathcal{S},\mathcal{L})$ such that $\varepsilon\in Q$, for some $Q\in\mathcal{M}$, is dense. Hence, by genericity there is some $\varepsilon'\in C_G^{dec}$ such that $\varepsilon'>\varepsilon$. It follows that $C_G^{dec}$ is unbounded in $\omega_2$.
		
		Let us now show that $C_G^{dec}$ is closed. Let $\varepsilon<\omega_2$ and suppose that $(\mathcal{M},d)\in\mathbb{M}^{dec}(\mathcal{S},\mathcal{L})$ forces $\check{\varepsilon}\notin\dot{C}_G^{dec}$, where $\dot{C}_G^{dec}$ is an $\mathbb{M}^{dec}(\mathcal{S},\mathcal{L})$-name for the club $C_G^{dec}$. We will show that $(\mathcal{M},d)$ can be extended to a condition forcing that $\check{\varepsilon}$ is not a limit point of $\dot{C}_G^{dec}$.
		
		By extending $(\mathcal{M},d)$ if necessary (using the unboundedness of $\mathcal{S}\cup\mathcal{L}$ and Lemma \ref{deco-ontop}), suppose that there is some $Q\in\mathcal{M}$ such that $\varepsilon_Q>\varepsilon$. Let $\varepsilon^+$ be the least ordinal in $\dom(\mathcal{M})$ such that $\varepsilon^+>\varepsilon$, and fix a model $Q^+\in\mathcal{M}(\varepsilon^+)$. We may assume that $C_G^{dec}\cap\varepsilon\neq\emptyset$. Let $\varepsilon^-$ be the immediate predecessor of $\varepsilon^+$ in $\dom(\mathcal{M})$, and note that $\varepsilon^-<\varepsilon$, otherwise we would contradict the assumption $(\mathcal{M},d)\Vdash``\check{\varepsilon}\notin\dot{C}_G^{dec}"$. Since $\varepsilon^+=\sup(Q^+\cap\omega_2)>\varepsilon$, there is some $\nu\in Q^+\cap\omega_2$ such that $\nu\geq\varepsilon$.
		
		Let $f$ be the function on $\dom(\mathcal{M})$ defined exactly like $d$, except that 
		\[
		f(\varepsilon^-)=d(\varepsilon^-)\cup\{\Psi_{Q^+[\omega_1],P[\omega_1]}(\nu):P\in\mathcal{M}(\varepsilon^+)\}.
		\]
		Note that since $\nu\in Q^+\cap\omega_2$, we have $\Psi_{Q^+[\omega_1],P[\omega_1]}(\nu)=\nu$. Thus $f(\varepsilon^-)=d(\varepsilon^-)\cup\{\nu\}$. It is not too hard to see that $f$ is a decoration on $\mathcal{M}$ and $\nu\in f(\varepsilon^-)$. Hence, $(\mathcal{M},f)\in\mathbb{M}^{dec}(\mathcal{S},\mathcal{L})$ and $(\mathcal{M},f)\leq(\mathcal{M},d)$. We claim that for every $(\mathcal{W},e)\in\mathbb{M}^{dec}(\mathcal{S},\mathcal{L})$ such that $(\mathcal{W},e)\leq(\mathcal{M},f)$, there is no $\eta\in\dom(\mathcal{W})$ such that $\varepsilon^-<\eta\leq\varepsilon$. Indeed, note that $\nu\in e(\varepsilon^-)$, so if $\varepsilon^*$ is the successor of $\varepsilon^-$ in $\dom(\mathcal{W})$, there must be some $W\in\mathcal{W}(\varepsilon^*)$ such that $\nu\in W$. Therefore, $\varepsilon_W=\varepsilon^*>\nu\geq\varepsilon$.
		
		It follows that $(\mathcal{M},f)$ forces that $C_G^{dec}$ has no elements between $\varepsilon^-$ and $\varepsilon$, and, in particular, that $\varepsilon$ is not a limit point of $C_G^{dec}$.
	\end{proof}
	
	\begin{lemma}
		If $B$ is an infinite subset of $\omega_2$ from the ground model $V$, then $B\not\subseteq C_G^{dec}$.
	\end{lemma}
	\begin{proof}
		It is immediate to see that the set $D_B$ of conditions $(\mathcal{M},d)\in\mathbb{M}^{dec}(\mathcal{S},\mathcal{L})$ for which there exist triples $(\varepsilon,\nu,\eta)$ such that $\varepsilon<\nu<\eta<\omega_2$, $\varepsilon\in\dom(\mathcal{M})$, $\nu\in B$, and $\eta\in d(\varepsilon)$, is a dense subset of $\mathbb{M}^{dec}(\mathcal{S},\mathcal{L})$.
	\end{proof}

	\subsection{Second application: Bounding canonical functions}
	
	The second application of the decorated forcing with $(\mathcal{S},\mathcal{L})$-symmetric systems involves adding a function on $\omega_2$ that bounds every canonical function below $\omega_3$ on a club.
	
	Recall that for every nonzero $\alpha<\omega_3$, if $\pi:\omega_2\to\alpha$ is a surjection, then the function $g_\alpha:\omega_2\to\omega_2$ defined by letting $g_\alpha(\nu)=\ot(\pi"\nu)$, for every $\nu<\omega_2$, is called a \emph{canonical function for $\alpha$}. Equivalently, $g'_\alpha$ is a canonical function for $\alpha$ if the poset $\mathcal{P}(\omega_2)/\mathrm{NS}_{\omega_2}$\footnote{$\mathrm{NS}_{\omega_2}$ denotes the non-stationary ideal on $\omega_2$.} forces $j(g_\alpha')(\omega_2^V)=\alpha$, where $j$ is the elementary embedding induced by the $\mathcal{P}(\omega_2)/\mathrm{NS}_{\omega_2}$-generic filter. In other words, the canonical function for $\alpha$ represents the ordinal $\alpha$ in the generic ultrapower of $V$ obtained after forcing with $\mathcal{P}(\omega_2)/\mathrm{NS}_{\omega_2}$.
	
	It is a well-known fact that canonical functions are uniquely determined modulo clubs, therefore given a nonzero ordinal $\alpha<\omega_3$, the choice of the surjection $\pi:\omega_2\to\alpha$ is irrelevant. More precisely, if $\pi_0,\pi_1:\omega_2\to\alpha$ are two surjections, then there is a club $C\subseteq\omega_2$ such that for every $\nu\in C$, $\pi_0"\nu=\pi_1"\nu$.
	
	
	Let $G$ be an $\mathbb{M}^{dec}(\mathcal{S},\mathcal{L})$-generic filter over $V$.
	
	
	\begin{theorem}
		If $2^{\aleph_1}=\aleph_2$ and $\mathcal{L}$ is stationary in $H(\kappa)$, the forcing $\mathbb{M}^{dec}(\mathcal{S},\mathcal{L})$ adds a function $f^{dec}_G:C^{dec}_G\to\omega_2$, such that for every nonzero $\alpha<\omega_3$, there is a canonical function $g_\alpha$ for $\alpha$ and some $\nu_0<\omega_2$, such that for every $\nu>\nu_0$ in $C^{dec}_G$, $g_\alpha(\nu)<f^{dec}_G(\nu)$.
	\end{theorem}
	\begin{proof}
		Define the function $f^{dec}_G: C^{dec}_G\to\omega_2$ as follows. For every $\varepsilon_Q\in C^{dec}_G$, let $f^{dec}_G(\varepsilon_Q)=\ot(Q[\omega_1]\cap\omega_3)$. Since for any $(\mathcal{S},\mathcal{L})$-symmetric system $\mathcal{M}$ and any $Q_0,Q_1\in\mathcal{M}$, if $\varepsilon_{Q_0}=\varepsilon_{Q_1}$, then $\ot(Q_0[\omega_1]\cap\omega_3)=\ot(Q_1[\omega_1]\cap\omega_3)$, the function $f^{dec}_G$ is well-defined.
		
		Let us show that for every nonzero $\alpha<\omega_3$, if $\pi_\alpha:\omega_2^V\to\alpha$ is a surjection, there is some $\nu_0<\omega_2$ such that for every $\nu>\nu_0$ in $C^{dec}_G$, $\ot(\pi_\alpha"\nu)<f^{dec}_G(\nu)$. Since $\mathcal{M}^{dec}_G$ covers $H(\kappa)$, we can find some $Q_0\in\mathcal{M}^{dec}_G$ such that $\alpha,\pi_\alpha\in Q_0$. Let $\nu_0:=\varepsilon_{Q_0}$ and fix some $\nu>\nu_0$ in $C_G^{dec}$. Since $\mathcal{M}_G^{dec}$ satisfies the shoulder axiom, there is $Q\in\mathcal{M}_G^{dec}$ such that $\varepsilon_Q=\nu$ and $Q_0\in Q[\omega_1]$. Note that $Q_0\subseteq Q[\omega_1]$, and thus $\alpha,\pi_\alpha\in Q[\omega_1]$. Therefore, $\pi_\alpha"\varepsilon_Q\subseteq Q[\omega_1]\cap\alpha$, and we can conclude that
		\[
		\ot(\pi_\alpha"\varepsilon_Q)\leq\ot(Q[\omega_1]\cap\alpha)<f^{dec}_G(\varepsilon_Q).
		\]
	\end{proof}

	\section{Simplified $(\omega_2,1)$-morasses}\label{section-morass}

	This final section will be devoted to defining a variant $\mathbb{M}^{mor}(\mathcal{S},\mathcal{L})$ of the forcing with decorated $(\mathcal{S},\mathcal{L})$-symmetric systems, so that all the nice preservation properties of the poset $\mathbb{M}^{dec}(\mathcal{S},\mathcal{L})$ still hold, but such that $\mathbb{M}^{mor}(\mathcal{S},\mathcal{L})$ adds a simplified $(\omega_2,1)$-morass.
	
	Morasses are complicated combinatorial objects that were introduced by Jensen in order to establish gap-2 two cardinal transfer principles in the constructible universe $L$.\footnote{Let $\mathscr{L}$ be a countable language with a unary predicate $U$. An $\mathscr{L}$-structure $A$ has \emph{type $(\mu,\eta)$} if $|A|=\mu$ and $|U^A|=\eta$. Jensen's gap-2 two cardinal principle $(\lambda^{++},\lambda)\rightarrow(\mu^{++},\mu)$ asserts that every theory with an $\mathscr{L}$-model of type $(\lambda^{++},\lambda)$ also has a model of type $(\mu^{++},\mu)$.} Similar to other principles living in $L$ such as $\diamondsuit$ and $\square$, morasses proved to be very useful in uncovering the combinatorial structure of $L$ (see \cite{Devlin1973Aspectsofconstructibility} and \cite{Kanamori1983Morassesincombinatorialsettheory}).
	
	Given an infinite cardinal $\kappa$, a $(\kappa,1)$-morass is a directed set of size $\kappa^+$, which allows one to construct objects of size $\kappa^+$ by $\kappa$-many approximations of size $<\kappa$. Morasses, in addition to existing in $L$, can also be added to models of $\mathrm{ZFC}$ by forcing. However, unlike the principles $\diamondsuit$ and $\square$, morasses are such complex objects that their use was restricted to a very small group of experts. Fortunately, starting in \cite{Velleman1982morassesdiamondandforcing}, Velleman undertook a very fruitful program to deeper understand morasses. Among many other results, Velleman introduced in \cite{Velleman1984simplifiedmorasses} a very striking simplification of $(\kappa,1)$-morasses, which he called \emph{simplified $(\kappa,1)$-morasses}, and whose existence was proven to be equivalent to the existence of $(\kappa,1)$-morasses (Corollary 4.2 in \cite{Velleman1984simplifiedmorasses}).
	
	Velleman's work opened the door of constructions along morasses to a broader audience of set theorists. In particular, Koszmider developed the technique of forcing with side conditions in morasses, inspired also by the work of Todor\v cevi\'c. A very important application of this technique, which appeared in \cite{Koszmider2000onstrongchainsofuncoutablefunctions}, was in proving the consistency of the existence of a strong chain of functions in $^{\omega_1}\omega_1$ of length $\omega_2$.\footnote{A sequence $\langle f_\alpha:\alpha<\omega_2\rangle$ of functions in $^{\omega_1}\omega_1$ is said to form a \emph{strong chain} if it is increasing coordinatewise, except for a finite set of coordinates.} But most importantly for our current purposes, Koszmider in that same paper pointed out that the technique of forcing with side conditions in morasses seemed to be equivalent to Todor\v cevi\'c's method of forcing with matrices (or symmetric systems) of models as side conditions. In fact, he went one step further and observed that a natural way of looking at $(\omega_1,1)$-morasses was to see them as families similar to $\{M\cap\omega_2:M\preceq H(\omega_3),|M|=\aleph_0\}$ with some extra properties. 
	
	This connection between symmetric systems of elementary submodels of a single type and Velleman's simplified morasses was made explicit by Miyamoto in \cite{Miyamoto2023Astronglysigmaclosedposetthatforcesasimplifiedmorass}, \cite{Miyamoto2014Matricesofisomorphicmodelsandmorasslikestructures}, \cite{Miyamoto2015Squaresbymatriceswithcoherentsequences}, \cite{Miyamoto2022Forcingaclubbyageneralizedfastfunction}, \cite{Miyamoto2023Negatingapartitionrelationbyafamilyofsimplified}. In this section, and building on some of the ideas of Miyamoto, we will strengthen these results by forcing a simplified $(\omega_2,1)$-morass with a variant of the forcing with decorated $(\mathcal{S},\mathcal{L})$-symmetric systems, that we will denote by $\mathbb{M}^{mor}(\mathcal{S},\mathcal{L})$. 
	
	By results of Velleman, the existence of a simplified $(\omega_2,1)$-morass implies the existence of an $(\omega_2,1)$-morass. Therefore, the generic extension obtained after forcing with $\mathbb{M}^{mor}(\mathcal{S},\mathcal{L})$ also satisfies all the consequences of the existence of an $(\omega_2,1)$-morass. Let us mention some of them:
	\begin{enumerate}
		\item The existence of an $\omega_2$-Kurepa tree (\cite{Velleman1982morassesdiamondandforcing}). 
		
		\item Weak $\square_{\omega_2}$ (\cite{ShelahStanley1982SforcingAblackboxtheoremformorasses}, \cite{Velleman1982morassesdiamondandforcing}).\footnote{The principle weak $\square_{\omega_2}$ is a weakening of $\square_{\omega_2}$ that asserts that there is a square sequence for a set of limit ordinals of $\omega_3$ containing all ordinals of cofinality $\omega_2$.}
		
		\item Jensen's gap-2 two cardinal principle $(\lambda^{++},\lambda)\rightarrow(\omega_3,\omega_1)$, for every $\lambda$ (\cite{Velleman1982morassesdiamondandforcing}).
		
		\item The existence of an $\omega_3$-Suslin tree, if $2^{\aleph_1}=\aleph_2$ (\cite{Velleman1982morassesdiamondandforcing}, \cite{ShelahStanley1982SforcingAblackboxtheoremformorasses}).\footnote{In fact, it implies the existence of an $\omega_3$-super-Suslin tree, which may fail to be $\omega_3$-Suslin, but implies the existence of an $\omega_3$-Suslin tree.}
		
		\item The existence of an $\omega_3$-Aronszajn tree (claimed in \cite{Koszmider2000onstrongchainsofuncoutablefunctions}).
		
		\item The existence of an $\omega_2$-thin-tall superatomic Boolean algebra and an $(\omega_2,\alpha,\omega_3)$-superatomic Boolean algebra, for every ordinal $\alpha$ such that $\omega_2\leq\alpha<\omega_3$ and $\cf(\alpha)=\omega_2$ (\cite{KoepkeMartinez1995superatomicbooleanalgebrasconstructedfrom}). 
		
	\end{enumerate}
	
	
	\subsection{Simplified morasses and elementary submodels}
	
	If $A,B\in[\omega_3]^{\leq\omega_1}$, let $A<B$ denote that $\alpha<\beta$ for every $\alpha\in A$ and every $\beta\in B$. Moreover, if $\mathcal{F}\subseteq[\omega_3]^{\leq\omega_1}$, we denote by $\mathcal{F}\restr A$ the set $\{B\in\mathcal{F}:B\subsetneq A\}$.
		
	\begin{definition}\label{def-simplified-morass}
		A \emph{simplified $(\omega_2,1)$-morass} is a family $\mathcal{F}\subseteq[\omega_3]^{\leq\omega_1}$ satisfying the following conditions:
		\begin{enumerate}
			\item $(\mathcal{F},\subseteq)$ is well-founded.
			
			\item $\mathcal{F}$ is \emph{locally small}, that is, $|\mathcal{F}\restr A|\leq\aleph_1$ for every $A\in\mathcal{F}$.
			
			\item $\mathcal{F}$ is \emph{homogeneous}, that is, if $A,B\in\mathcal{F}$ are such that $\rank(A)=\rank(B)$,\footnote{Recall that any well-founded relation (in this case $(\mathcal{F},\subseteq)$) has a naturally associated notion of rank.} then $A$ and $B$ have the same order-type and if $f_{A,B}$ is the unique order-preserving map from $A$ onto $B$, then $$\mathcal{F}\restr B=\{f_{A,B}"(C):C\in\mathcal{F}\restr A\}.$$
			
			\item $\mathcal{F}$ is \emph{directed}, that is, for any two $A,B\in\mathcal{F}$ there is $C\in\mathcal{F}$ such that $A,B\subseteq C$.
			
			\item $\mathcal{F}$ is \emph{locally almost directed}, that is, for every $A\in\mathcal{F}$ either
			\begin{enumerate}
				\item $\mathcal{F}\restr A$ is directed, or
				\item $\mathcal{F}\restr A$ \emph{splits}, that is, there are $A_0,A_1\in\mathcal{F}$ such that
				\begin{itemize}
					\item $\rank(A_0)=\rank(A_1)$,
					\item $A_0\cap A_1$ is a proper initial segment of both $A_0$ and $A_1$,
					\item $A_0\cap A_1<A_0\setminus A_1<A_1\setminus A_0$, and
					\item $\mathcal{F}\restr A=\mathcal{F}\restr A_0\cup\mathcal{F}\restr A_1\cup\{A_0,A_1\}$.
				\end{itemize}    
				In this case, we call the pair $\{A_0,A_1\}$ \emph{a split of} $\mathcal{F}\restr A$.
			\end{enumerate}
			
			\item $\mathcal{F}$ \emph{covers} $\omega_3$, that is, $\bigcup\mathcal{F}=\omega_3$.
		\end{enumerate}
	\end{definition}
	
	
	For the rest of the section, we make the convention that, since there will be no possible confusion, every time that we refer to a directed set, what we actually mean is that the set is directed with respect to inclusion.
		
	The simplified $(\omega_2,1)$-morass that we intend to add generically with the forcing $\mathbb{M}^{mor}(\mathcal{S},\mathcal{L})$, will be given by the set
	\[
	\mathcal{F}_G:=\{Q[\omega_1]\cap\omega_3:Q\in\mathcal{M}_G^{mor}\},
	\]
	where $G$ is a generic filter for $\mathbb{M}^{mor}(\mathcal{S},\mathcal{L})$ over $V$, and $\mathcal{M}_G^{mor}$ is the generic $(\mathcal{S},\mathcal{L})$-symmetric system $\mathcal{M}_G^{mor}=\{Q\in\mathcal{M}:(\mathcal{M},d)\in G\}$. Lemma \ref{agreement} ensures that, if $Q_0,Q_1\in\mathcal{S}\cup\mathcal{L}$ are $\omega_1$-isomorphic, then $Q_0[\omega_1]\cap Q_1[\omega_1]\cap\omega_3$ is an initial segment of both $Q_0[\omega_1]\cap\omega_3$ and $Q_1[\omega_1]\cap\omega_3$. However, it is not in general true that either 
	\[
	(Q_0[\omega_1]\setminus Q_1[\omega_1])\cap\omega_3<(Q_1[\omega_1]\setminus Q_0[\omega_1])\cap\omega_3,
	\]
	or 
	\[(Q_1[\omega_1]\setminus Q_0[\omega_1])\cap\omega_3<(Q_0[\omega_1]\setminus Q_1[\omega_1])\cap\omega_3.
	\]
	We introduce the following strengthening of the notion of $\omega_1$-isomorphism to account for this.

	
	
	\begin{definition}
		Two models $Q_0,Q_1\in\mathcal{S}\cup\mathcal{L}$ are \emph{strongly $\omega_1$-isomorphic} if the following conditions hold:
		\begin{enumerate}
			\item $Q_0$ and $Q_1$ are $\omega_1$-isomorphic.
			
			\item $Q_0[\omega_1]\cap Q_1[\omega_1]\cap\omega_3$ is a proper initial segment of both $Q_0[\omega_1]\cap\omega_3$ and $Q_1[\omega_1]\cap\omega_3$.
			
			\item Either $$(Q_0[\omega_1]\setminus Q_1[\omega_1])\cap\omega_3<(Q_1[\omega_1]\setminus Q_0[\omega_1])\cap\omega_3,$$ or $$(Q_1[\omega_1]\setminus Q_0[\omega_1])\cap\omega_3<(Q_0[\omega_1]\setminus Q_1[\omega_1])\cap\omega_3.$$
		\end{enumerate}
	
		In this case, we will say that the $\omega_1$-isomorphism $\Psi_{Q_0[\omega_1],Q_1[\omega_1]}$ is a \emph{strong $\omega_1$-isomorphism}.
	\end{definition}

	The next proposition is easy to verify.

	\begin{proposition}\label{morass-copy-strongly-iso}
		Let $Q_0,Q_1,P_0,P_1\in\mathcal{S}\cup\mathcal{L}$. If $Q_0$ and $Q_1$ are strongly $\omega_1$-isomorphic, $P_0$ and $P_1$ are $\omega_1$-isomorphic, and $Q_0,Q_1\in P_0[\omega_1]$, then $\Psi_{P_0[\omega_1],P_1[\omega_1]}(Q_0)$ and $\Psi_{P_0[\omega_1],P_1[\omega_1]}(Q_1)$ are strongly $\omega_1$-isomorphic.
	\end{proposition}

	\begin{proposition}\label{morass-strongly-iso}
		Let $Q_0,Q_1\in\mathcal{S}\cup\mathcal{L}$ be two $\omega_1$-isomorphic models.  If there are $\alpha<\beta<\omega_3$ such that $\alpha\in Q_0[\omega_1]$, $\beta\in Q_1[\omega_1]$, $\Psi_{Q_0[\omega_1],Q_1[\omega_1]}(\alpha)=\beta$, $Q_0[\omega_1]\cap\omega_3\subseteq\beta$ and $Q_0[\omega_1]\cap\alpha=Q_1[\omega_1]\cap\beta$, then $Q_0$ and $Q_1$ are strongly $\omega_1$-isomorphic. 
	\end{proposition}
	\begin{proof}
		Let $\xi_0\in (Q_0[\omega_1]\setminus Q_1[\omega_1])\cap\omega_3$ and $\xi_1\in (Q_1[\omega_1]\setminus Q_0[\omega_1])\cap\omega_3$. We claim that $\xi_0<\xi_1$. Suppose on the contrary that $\xi_0\geq\xi_1$. Then, as $\xi_0\in Q_0[\omega_1]\cap\omega_3\subseteq\beta$, we have $\xi_1<\beta$, and hence
		\[
		\xi_1\in Q_1[\omega_1]\cap\beta=Q_0[\omega_1]\cap\alpha,
		\]
		which contradicts the assumption $\xi_1\in (Q_1[\omega_1]\setminus Q_0[\omega_1])\cap\omega_3$.
	\end{proof}

	\begin{notation}
		If $\mathcal{M}\subseteq\mathcal{S}\cup\mathcal{L}$, we will denote by $\mathcal{F}_\mathcal{M}$ the set
		\[
		\mathcal{F}_\mathcal{M}=\{Q[\omega_1]\cap\omega_3:Q\in\mathcal{M}\}.
		\]
	\end{notation}	
		
	\begin{definition}
		Let $\mathcal{M}$ be a finite subset of $\mathcal{S}\cup\mathcal{L}$. We say that $\mathcal{M}$ is \emph{almost directed} if either
		\begin{enumerate}
			\item $\mathcal{F}_\mathcal{M}$ is directed with respect to inclusion, or
			\item $\mathcal{F}_\mathcal{M}$ \emph{splits}, that is, there are two strongly $\omega_1$-isomorphic models $P_0,P_1\in\mathcal{M}$ such that
			\begin{align*}
				\mathcal{F}_\mathcal{M}=\{P_0[\omega_1]\cap\omega_3,P_1[\omega_1]\cap\omega_3\}&\cup\mathcal{F}_\mathcal{M}\restr (P_0[\omega_1]\cap\omega_3)\\
				&\cup\mathcal{F}_\mathcal{M}\restr(P_1[\omega_1]\cap\omega_3).
			\end{align*}	 
			In this case, we call the pair $\{P_0[\omega_1]\cap\omega_3,P_1[\omega_1]\cap\omega_3\}$ \emph{a split of} $\mathcal{F}_\mathcal{M}$.
		\end{enumerate}
	\end{definition}
	
	\begin{definition}
		Let $\mathcal{M}$ be a finite subset of $\mathcal{S}\cup\mathcal{L}$. We say that $\mathcal{M}$ is \emph{locally almost directed} if for every $Q\in\mathcal{M}$, the set $\mathcal{M}\cap Q[\omega_1]$ is almost directed.
	\end{definition}
	
	
	
	
	
	The next two results follow immediately from the definitions.
	
	\begin{proposition}
		Let $\mathcal{M}$ be a finite subset of $\mathcal{S}\cup\mathcal{L}$ and let $Q,P\in\mathcal{M}$ be such that $P[\omega_1]\cap\omega_3\subseteq Q[\omega_1]\cap\omega_3$. Then 
		\[
		(\mathcal{F}_{\mathcal{M}}\restr(Q[\omega_1]\cap\omega_3))\restr(P[\omega_1]\cap\omega_3)=\mathcal{F}_{\mathcal{M}}\restr(P[\omega_1]\cap\omega_3).
		\] 
	\end{proposition}
	
	\begin{proposition}\label{morass-prop}
		Let $\mathcal{M}$ be a nonempty finite subset of $\mathcal{S}\cup\mathcal{L}$. Then the following hold:
		\begin{enumerate}
			\item $\mathcal{F}_\mathcal{M}$ is directed with respect to inclusion if and only if there is $Q\in\mathcal{M}$ such that $\varepsilon_Q$ is maximal in $\dom(\mathcal{M})$ and
			\[
			\mathcal{F}_\mathcal{M}=\{Q[\omega_1]\cap\omega_3\}\cup\mathcal{F}_\mathcal{M}\restr (Q[\omega_1]\cap\omega_3).
			\]
			
			\item If $\{P_0[\omega_1]\cap\omega_3,P_1[\omega_1]\cap\omega_3\}$ is a split of $\mathcal{F}_\mathcal{M}$, then $\varepsilon_{P_0}=\varepsilon_{P_1}$ are maximal in $\dom(\mathcal{M})$.
		\end{enumerate}
	\end{proposition}
	\begin{proof}
		Let us denote $\varepsilon_{max}=\max(\dom(\mathcal{M}))$ and 
		\[
		\mathcal{F}_{max}=\{P[\omega_1]\cap\omega_3:P\in\mathcal{M}(\varepsilon_{max})\}.
		\]
		
		For (1), let $A$ be $\subseteq$-maximal in the finite set $\mathcal{F}_{max}$. Hence $A\in\mathcal{F}_{max}$ and for any $P[\omega_1]\cap\omega_3\in\mathcal{F}_{max}$, $A$ is not a proper subset of $P[\omega_1]\cap\omega_3$. It follows easily that $\mathcal{F}_{\mathcal{M}}$ is directed if and only if $\mathcal{F}_{\mathcal{M}}=\{A\}\cup\mathcal{F}_{\mathcal{M}}\restr A$.
		
		For (2), let $\{A_0,A_1\}$ be a split of $\mathcal{F}_{\mathcal{M}}$, and suppose that $Q_0,Q_1\in\mathcal{M}$ are such that
		\begin{itemize}
			\item $Q_0$ and $Q_1$ are strongly $\omega_1$-isomorphic,
			\item $A_0=Q_0[\omega_1]\cap\omega_3$ and $A_1=Q_1[\omega_1]\cap\omega_3$, and
			\item $\mathcal{F}_{\mathcal{M}}=\{A_0,A_1\}\cup\mathcal{F}_{\mathcal{M}}\restr A_0\cup\mathcal{F}_{\mathcal{M}}\restr A_1$.
		\end{itemize}
		Take any two $P_0,P_1\in\mathcal{M}$ such that $A_0=P_0[\omega_1]\cap\omega_3$ and $A_1=P_1[\omega_1]\cap\omega_3$. We claim that $\varepsilon_{P_0}=\varepsilon_{P_1}=\varepsilon_{max}$. For $i=0,1$,
		\begin{align*}
		\varepsilon_{P_i}=P_i[\omega_1]\cap\omega_2&=(P_i[\omega_1]\cap\omega_3)\cap\omega_2\\
		&=(Q_i[\omega_1]\cap\omega_3)\cap\omega_2=Q_i[\omega_1]\cap\omega_2=\varepsilon_{Q_i}. 
		\end{align*}
		Since $\varepsilon_{Q_0}=\varepsilon_{Q_1}$, we also have $\varepsilon_{P_0}=\varepsilon_{P_1}$. Next, for any $Q\in\mathcal{M}$, either $Q[\omega_1]\cap\omega_3\subseteq A_0$ or $Q[\omega_1]\cap\omega_3\subseteq A_1$. Hence,
		\[
		\varepsilon_Q=Q[\omega_1]\cap\omega_2=(Q[\omega_1]\cap\omega_3)\cap\omega_2\subseteq A_0\cap\omega_2=A_1\cap\omega_2,
		\]
		and we can conclude that $\varepsilon_{P_0}=\varepsilon_{P_1}=\varepsilon_{max}$.
	\end{proof}
	
	\begin{proposition}\label{morass-rank}
		Let $\mathcal{M}$ be an $(\mathcal{S},\mathcal{L})$-symmetric system and let $Q,P\in\mathcal{M}$. If $P[\omega_1]\cap\omega_3\in\mathcal{F}_{\mathcal{M}}\restr(Q[\omega_1]\cap\omega_3)$, then $\varepsilon_P<\varepsilon_Q$.
	\end{proposition}
	\begin{proof}
		If $P[\omega_1]\cap\omega_3\in\mathcal{F}_{\mathcal{M}}\restr(Q[\omega_1]\cap\omega_3)$, then $P[\omega_1]\cap\omega_3\subsetneq Q[\omega_1]\cap\omega_3$. Therefore, by Corollary \ref{prelim-height-hull},
		\[
		\varepsilon_P=P[\omega_1]\cap\omega_2\leq Q[\omega_1]\cap\omega_2=\varepsilon_Q.
		\]
		Note that if $\varepsilon_P=\varepsilon_Q$, then $P$ and $Q$ would be $\omega_1$-isomorphic. Hence $P[\omega_1]\cap\omega_3$ would be fixed by the isomorphism $\Psi_{P[\omega_1],Q[\omega_1]}$, and thus
		\[
		P[\omega_1]\cap\omega_3=\Psi_{P[\omega_1],Q[\omega_1]}(P[\omega_1]\cap\omega_3)=Q[\omega_1]\cap\omega_3,
		\]
		which is impossible because $P[\omega_1]\cap \omega_3$ is a proper subset of $Q[\omega_1]\cap\omega_3$ by assumption. 
	\end{proof}
	
	\begin{proposition}\label{morass-restr}
		Let $\mathcal{M}$ be an $(\mathcal{S},\mathcal{L})$-symmetric system. For all $Q\in\mathcal{M}$,
		\[
		\mathcal{F}_\mathcal{M}\cap Q[\omega_1]=\mathcal{F}_\mathcal{M}\restr(Q[\omega_1]\cap\omega_3)=\mathcal{F}_{\mathcal{M}\cap Q[\omega_1]}.
		\]
	\end{proposition}
	\begin{proof}
		Let us start with the equality $\mathcal{F}_\mathcal{M}\cap Q[\omega_1]=\mathcal{F}_\mathcal{M}\restr(Q[\omega_1]\cap\omega_3)$. If $P[\omega_1]\cap\omega_3\in\mathcal{F}_\mathcal{M}\cap Q[\omega_1]$, then clearly $P[\omega_1]\cap\omega_3\subsetneq Q[\omega_1]\cap\omega_3$. Let us show the other inclusion. Suppose that $R[\omega_1]\cap\omega_3\in\mathcal{F}_\mathcal{M}\restr(Q[\omega_1]\cap\omega_3)$. Then by Proposition \ref{morass-rank}, $\varepsilon_R<\varepsilon_Q$. Hence, using the shoulder axiom for $\mathcal{M}$, we can find $Q'\in\mathcal{M}(\varepsilon_Q)$ such that $R\in Q'[\omega_1]$ and, in particular, $R[\omega_1]\cap\omega_3\in Q'[\omega_1]$. Since $R[\omega_1]\cap\omega_3\subseteq Q'[\omega_1]\cap Q[\omega_1]$, the isomorphism $\Psi_{Q'[\omega_1],Q[\omega_1]}$ fixes $R[\omega_1]\cap\omega_3$, and hence,
		\[
		R[\omega_1]\cap\omega_3=\Psi_{Q'[\omega_1],Q[\omega_1]}(R[\omega_1]\cap\omega_3)\in \mathcal{F}_\mathcal{M}\cap Q[\omega_1].
		\]
		
		Let us now prove the equality $\mathcal{F}_\mathcal{M}\cap Q[\omega_1]=\mathcal{F}_{\mathcal{M}\cap Q[\omega_1]}$. First, we show the inclusion $\mathcal{F}_\mathcal{M}\cap Q[\omega_1]\subseteq\mathcal{F}_{\mathcal{M}\cap Q[\omega_1]}$. Suppose that $P[\omega_1]\cap\omega_3\in\mathcal{F}_\mathcal{M}\cap Q[\omega_1]$. By Proposition \ref{morass-rank} and the equality 
		\[
		\mathcal{F}_\mathcal{M}\cap Q[\omega_1]=\mathcal{F}_\mathcal{M}\restr(Q[\omega_1]\cap\omega_3),
		\]
		proven in the last paragraph, we have $\varepsilon_P<\varepsilon_Q$, and hence, by the shoulder axiom for $\mathcal{M}$, there is some $Q'\in\mathcal{M}$ such that $\varepsilon_{Q'}=\varepsilon_Q$ and $P\in Q'[\omega_1]$. Note that, in particular, $P[\omega_1]\cap\omega_3\in Q'[\omega_1]$. Therefore, if we let $R:=\Psi_{Q'[\omega_1],Q[\omega_1]}(P)$, which is an element of $\mathcal{M}\cap Q[\omega_1]$ by the symmetry of $\mathcal{M}$, we can conclude that, since $P[\omega_1]\cap\omega_3\in Q[\omega_1]\cap Q'[\omega_3]$,
		\begin{align*}
			P[\omega_1]\cap\omega_3&=\Psi_{Q'[\omega_1],Q[\omega_1]}(P[\omega_1]\cap\omega_3)\\
			&=R[\omega_1]\cap\omega_3\in\mathcal{F}_{\mathcal{M}\cap Q[\omega_1]}.
		\end{align*}
		Let us now show the other inclusion. Suppose that $P[\omega_1]\cap\omega_3\in\mathcal{F}_{\mathcal{M}\cap Q[\omega_1]}$. Then there is some $R\in\mathcal{M}\cap Q[\omega_1]$ such that $P[\omega_1]\cap\omega_3=R[\omega_1]\cap\omega_3$. However, note that, as $R\in Q[\omega_1]$, we have $R[\omega_1]\cap\omega_3\in Q[\omega_1]$. Therefore, $P[\omega_1]\cap\omega_3$ is a member of $\mathcal{F}_\mathcal{M}\cap Q[\omega_1]$. 
	\end{proof}	
	
	\subsection{Morass-like symmetric systems}
	
	\begin{definition}\label{def-morass-like-SSM2T}		
		An $(\mathcal{S},\mathcal{L})$-symmetric system $\mathcal{M}$ will be called \emph{morass-like} if $\mathcal{M}$ is almost directed and locally almost directed.
	\end{definition}
	
	\begin{definition}
		Let $\mathbb{M}^{mor}(\mathcal{S},\mathcal{L})$ be the poset whose conditions are pairs $(\mathcal{M},d)$ such that $(\mathcal{M},d)\in\mathbb{M}^{dec}(\mathcal{S},\mathcal{L})$ and $\mathcal{M}$ is morass-like. The order on $\mathbb{M}^{mor}(\mathcal{S},\mathcal{L})$ is induced by the order on $\mathbb{M}^{dec}(\mathcal{S},\mathcal{L})$. That is, if $(\mathcal{M},d),(\mathcal{N},f)\in\mathbb{M}^{mor}(\mathcal{S},\mathcal{L})$, we set $(\mathcal{N},f)\leq(\mathcal{M},d)$ if and only if
		\begin{enumerate}
			\item $\mathcal{N}\supseteq\mathcal{M}$, and
			\item $f(\varepsilon)\supseteq d(\varepsilon)$, for every $\varepsilon\in\dom(\mathcal{M})$.
		\end{enumerate}
	\end{definition}
	
	\begin{definition}
		Let $\mathcal{M}$ be an $(\mathcal{S},\mathcal{L})$-symmetric system. Given $Q,P\in\mathcal{M}$, we will say that $P$ is an \emph{$\mathcal{M}$-successor of $Q$}, or equivalently, that $Q$ is an \emph{$\mathcal{M}$-predecessor of $P$}, if $Q\in P$ and $\varepsilon_P$ is the immediate successor of $\varepsilon_Q$ in $\dom(\mathcal{M})$. 
	\end{definition}
	
	\begin{remark}\label{morass-remark}
		Suppose that $\mathcal{M}$ is a morass-like $(\mathcal{S},\mathcal{L})$-symmetric system. Given $Q\in\mathcal{M}$, one of the consequences of Proposition \ref{morass-prop} is that the almost directedness of $\mathcal{M}\cap Q[\omega_1]$ is always witnessed by $\mathcal{M}$-predecessors of $Q$. 
		
		For the rest of the paper, we freely use the equation in Proposition \ref{morass-restr} for relevant $(\mathcal{S},\mathcal{L})$-symmetric systems. In particular, when we check the locally almost directedness of $\mathcal{M}$ with respect to the subsystem $\mathcal{M}\cap Q[\omega_1]$, where $Q\in\mathcal{M}$, we indistinctly denote the projection of $\mathcal{M}\cap Q[\omega_1]$ by $\mathcal{F}_{\mathcal{M}\cap Q[\omega_1]}$, $\mathcal{F}_{\mathcal{M}}\cap Q[\omega_1]$, or $\mathcal{F}_{\mathcal{M}}\restr(Q[\omega_1]\cap\omega_3)$. Hence, as split of $\mathcal{F}_{\mathcal{M}}\cap Q[\omega_1]$ means the split of the projection $\mathcal{F}_{\mathcal{M}\cap Q[\omega_1]}$.
		
		We caution, however, that the notation may suggest a procedure in which $\mathcal{F}_{\mathcal{M}}$ is first computed using the full system $\mathcal{M}$ and the resulting value is then used to compute $\mathcal{F}_\mathcal{M}\cap Q[\omega_1]$ or $\mathcal{F}_\mathcal{M}\restr(Q[\omega_1]\cap\omega_3)$. Although $\mathcal{F}_{\mathcal{M}\cap Q[\omega_1]}$, $\mathcal{F}_{\mathcal{M}}\cap Q[\omega_1]$, and $\mathcal{F}_{\mathcal{M}}\restr(Q[\omega_1]\cap\omega_3)$ are equal, this may obscure the conceptual distinction between computing $\mathcal{F}_\mathcal{M}\cap Q[\omega_1]$ or $\mathcal{F}_\mathcal{M}\restr(Q[\omega_1]\cap\omega_3)$ via $\mathcal{F}_\mathcal{M}$ and computing $\mathcal{F}_{\mathcal{M}\cap Q[\omega_1]}$ directly from the subsystem $\mathcal{M}\cap Q[\omega_1]$.
		
		The same applies to other projections as well.

	\end{remark}
	
	
	The proof of the next result, which is analogous to Proposition \ref{isocopy}, uses Proposition \ref{morass-copy-strongly-iso} and is an easy exercise.
		
	\begin{proposition}\label{morass-isocopy}
		Let $\mathcal{M}$ be a morass-like $(\mathcal{S},\mathcal{L})$-symmetric system. Let $Q_0,Q_1\in\mathcal{S}\cup\mathcal{L}$ be $\omega_1$-isomorphic. If $\mathcal{M}\subseteq Q_0[\omega_1]$, then $\Psi_{Q_0[\omega_1],Q_1[\omega_1]}(\mathcal{M})$ is a morass-like $(\mathcal{S},\mathcal{L})$-symmetric system.
	\end{proposition}

	\subsection{Amalgamation lemmas}
	
	As in previous sections, we will show that $\mathbb{M}^{mor}(\mathcal{S},\mathcal{L})$ satisfies the standard amalgamation lemmas for $(\mathcal{S},\mathcal{L})$-symmetric systems. We will show that if we start from morass-like $(\mathcal{S},\mathcal{L})$-symmetric systems, then the $(\mathcal{S},\mathcal{L})$-symmetric systems witnessing the conclusions of the lemmas from Subsection \ref{subsection-amalgamation-lemmas} are also morass-like. From this, it is easy to derive amalgamation lemmas for $\mathbb{M}^{mor}(\mathcal{S},\mathcal{L})$, analogous to the ones for the decorated forcing $\mathbb{M}^{dec}(\mathcal{S},\mathcal{L})$ from Subsection \ref{subsection-deco-amal-lemmas}.
	
	Indeed, recall that the first components of the conditions witnessing the amalgamation lemmas for the decorated forcing $\mathbb{M}^{dec}(\mathcal{S},\mathcal{L})$ are the $(\mathcal{S},\mathcal{L})$-symmetric systems obtained in Subsection \ref{subsection-amalgamation-lemmas}. Therefore, once we show that these $(\mathcal{S},\mathcal{L})$-symmetric systems can be made morass-like, the corresponding decorations can be obtained by essentially repeating the constructions from Subsection \ref{subsection-deco-amal-lemmas}. 
	

	\begin{lemma}\label{morass-on-top}
		Let $\mathcal{M}$ be a morass-like $(\mathcal{S},\mathcal{L})$-symmetric system and let $Q\in\mathcal{S}\cup\mathcal{L}$ such that $\mathcal{M}\subseteq Q$. Then there is a morass-like $(\mathcal{S},\mathcal{L})$-symmetric system $\mathcal{M}_Q$ such that $\mathcal{M}\cup\{Q\}\subseteq\mathcal{M}_Q$.
	\end{lemma}
	\begin{proof}
		If $Q\in\mathcal{L}$, it is easy to see that $\mathcal{M}_Q=\mathcal{M}\cup\{Q\}$ is a morass-like $(\mathcal{S},\mathcal{L})$-symmetric system, where $Q[\omega_1]\cap\omega_3$ witnesses that $\mathcal{F}_{\mathcal{M}_Q}$ is directed with respect to inclusion.
		
		Suppose that $Q\in\mathcal{S}$, and recall from the proof of Lemma \ref{amal-ontop} that 
		\[
		\mathcal{M}_Q=\mathcal{M}\cup\{Q\}\cup\{X\cap Q:X\in\mathcal{M}\cap\mathcal{L}\}
		\]
		is an $(\mathcal{S},\mathcal{L})$-symmetric system. We claim that it is also morass-like. 
		
		Note that $\mathcal{F}_{\mathcal{M}_Q}$ is directed, as $P[\omega_1]\cap\omega_3\in Q[\omega_1]$ for every $P\in\mathcal{M}_Q\setminus\{Q\}$. Hence, we only need to check that $\mathcal{M}_Q$ is locally almost directed.
		
		Let $X\in\mathcal{M}\cap\mathcal{L}$. Note that, as $\mathcal{M}\subseteq Q$, for every $P\in\mathcal{M}\cap X$ we have $P\in X\cap Q$. Hence, in particular, $P[\omega_1]\cap\omega_3\subseteq (X\cap Q)[\omega_1]\cap\omega_3$. Moreover, note that $(X\cap Q)[\omega_1]\cap\omega_3\subseteq X\cap\omega_3$. From this, we can deduce that the following conditions hold for every $P\in\mathcal{M}_Q$, which in turn imply the local almost directedness of $\mathcal{M}_Q$:
		\begin{enumerate}
			\item Suppose that $P\in\mathcal{M}\cap\mathcal{S}$.
			\begin{itemize}
				\item If $\mathcal{F}_\mathcal{M}\cap P[\omega_1]$ is directed, then $\mathcal{F}_{\mathcal{M}_Q}\cap P[\omega_1]$ is directed, and this is witnessed by the same model witnessing the directedness of $\mathcal{F}_\mathcal{M}\cap P[\omega_1]$.
				\item If $\mathcal{F}_\mathcal{M}\cap P[\omega_1]$ splits, then $\mathcal{F}_{\mathcal{M}_Q}\cap P[\omega_1]$ splits, and this is witnessed by the same pair of models witnessing that $\mathcal{F}_\mathcal{M}\cap P[\omega_1]$ splits.
			\end{itemize}
			
			\item Suppose that $P\in\mathcal{M}\cap\mathcal{L}$. Then $\mathcal{F}_{\mathcal{M}_Q}\cap P$ is directed, as witnessed by $(P\cap Q)[\omega_1]\cap\omega_3$.
			
			\item Suppose that $P=X\cap Q$ for some $X\in\mathcal{M}\cap\mathcal{L}$. 
			\begin{itemize}
				\item If $\mathcal{F}_\mathcal{M}\cap X$ is directed, then $\mathcal{F}_{\mathcal{M}_Q}\cap P[\omega_1]$ is directed, and this is witnessed by the same model witnessing the directedness of $\mathcal{F}_\mathcal{M}\cap X$.  
				
				\item If $\mathcal{F}_\mathcal{M}\cap X$ splits, then $\mathcal{F}_{\mathcal{M}_Q}\cap P[\omega_1]$ splits, and this is witnessed by the same pair of models witnessing that $\mathcal{F}_\mathcal{M}\cap X$ splits. 
			\end{itemize}
			
			\item Suppose that $P=Q$.
			\begin{itemize}
				\item If $\mathcal{F}_\mathcal{M}$ is directed, then $\mathcal{F}_{\mathcal{M}_Q}\cap Q[\omega_1]$ is directed, and this is witnessed by the same model witnessing the directedness of $\mathcal{F}_\mathcal{M}$.
				
				\item If $\mathcal{F}_\mathcal{M}$ splits, then $\mathcal{F}_{\mathcal{M}_Q}\cap Q[\omega_1]$ splits, and this is witnessed by the same pair of models witnessing that $\mathcal{F}_\mathcal{M}$ splits.
			\end{itemize}
		\end{enumerate}
	\end{proof}
	
	\begin{corollary}\label{morass-on-top-corollary}
		Let $(\mathcal{M},d)\in\mathbb{M}^{mor}(\mathcal{S},\mathcal{L})$ and let $Q\in\mathcal{S}\cup\mathcal{L}$ be such that $(\mathcal{M},d)\subseteq Q$. Then there is a condition $(\mathcal{M}_Q,d_Q)\in\mathbb{M}^{mor}(\mathcal{S},\mathcal{L})$ such that $Q\in\mathcal{M}_Q$ and $(\mathcal{M}_Q,d_Q)\leq(\mathcal{M},d)$.
	\end{corollary}
	\begin{proof}
		The proof is exactly the same as the one for Lemma \ref{deco-ontop}, except that we use Lemma \ref{morass-on-top} instead of Lemma \ref{amal-ontop}.
	\end{proof}
	
	\begin{lemma}\label{morass-refl}
		Let $\mathcal{M}$ be a morass-like $(\mathcal{S},\mathcal{L})$-symmetric system and let $Q\in\mathcal{M}\cup\mathcal{M}[\omega_1]$. Then $\mathcal{M}\cap Q$ is a morass-like $(\mathcal{S},\mathcal{L})$-symmetric system.
	\end{lemma}
	\begin{proof}
		We know from Lemma \ref{amal-rest} that $\mathcal{M}\cap Q$ is an $(\mathcal{S},\mathcal{L})$-symmetric system.
		
		If $Q\in\mathcal{L}$ or $Q\in\mathcal{M}[\omega_1]$, by Proposition \ref{morass-restr} we have  
		\[
		\mathcal{F}_{\mathcal{M}\cap Q}=\mathcal{F}_{\mathcal{M}\cap Q[\omega_1]}=\mathcal{F}_\mathcal{M}\cap Q[\omega_1]=\mathcal{F}_\mathcal{M}\cap Q.\]
		Hence, as $\mathcal{M}$ is locally almost directed, then $\mathcal{M}\cap Q$ is almost directed and locally almost directed.
		
		Suppose now that $Q\in\mathcal{S}$. In order to see that $\mathcal{M}\cap Q$ is almost directed we only need to observe that the $\mathcal{M}$-predecessors of $Q$ are also $(\mathcal{M}\cap Q)$-predecessors. Therefore, if $\mathcal{F}_\mathcal{M}\cap Q[\omega_1]$ is directed and, by Remark \ref{morass-remark}, this is witnessed by an $\mathcal{M}$-predecessor of $Q$, then $\mathcal{F}_{\mathcal{M}\cap Q}$ is directed, and this is witnessed by the same model witnessing the directedness of $\mathcal{F}_\mathcal{M}\cap Q[\omega_1]$. Similarly, if $\mathcal{F}_\mathcal{M}\cap Q[\omega_1]$ splits, then $\mathcal{F}_{\mathcal{M}\cap Q}$ splits, and this is witnessed by the same pair of $\mathcal{M}$-predecessors of $Q$ witnessing that $\mathcal{F}_\mathcal{M}\cap Q[\omega_1]$ splits.
		
		Let us now show that $\mathcal{M}\cap Q$ is locally almost directed. It is important to keep in mind the picture of $\mathcal{M}\cap Q$ given by Proposition \ref{prop-res-gaps} and Lemma \ref{amal-rest}. If $X\in\mathcal{M}\cap Q\cap\mathcal{L}$, note that the $(\mathcal{M}\cap Q)$-predecessors of $X$ are exactly the $\mathcal{M}$-predecessors of $X\cap Q$. Therefore, using Remark \ref{morass-remark},
		\begin{itemize}
			\item if $\mathcal{F}_\mathcal{M}\cap (X\cap Q)[\omega_1]$ is directed, the $\mathcal{M}$-predecessor of $X\cap Q$ witnessing this, also witnesses that $\mathcal{F}_{\mathcal{M}\cap Q}\cap X$ is directed, and
			
			\item if $\mathcal{F}_\mathcal{M}\cap (X\cap Q)[\omega_1]$ splits, the pair of $\mathcal{M}$-predecessors of $X\cap Q$ witnessing this, also witnesses that $\mathcal{F}_{\mathcal{M}\cap Q}\cap X$ splits.
		\end{itemize}
		If $N\in\mathcal{M}\cap Q\cap\mathcal{S}$ and $P$ is an $\mathcal{M}$-predecessor of $N$, then $P\in\mathcal{M}\cap Q$. Hence, every $\mathcal{M}$-predecessor of $N$ is also an $(\mathcal{M}\cap Q)$-predecessor. Therefore, 
		\begin{itemize}
			\item if an $\mathcal{M}$-predecessor of $N$ witnesses that $\mathcal{F}_\mathcal{M}\cap N[\omega_1]$ is directed, then it also witnesses that $\mathcal{F}_{\mathcal{M}\cap Q}\cap N[\omega_1]$ is directed, and
			
			\item if a pair of $\mathcal{M}$-predecessors of $N$ is a split of $\mathcal{F}_\mathcal{M}\cap N[\omega_1]$, then it is also a split of $\mathcal{F}_{\mathcal{M}\cap Q}\cap N[\omega_1]$.
		\end{itemize}
		So, we can conclude that $\mathcal{M}\cap Q$ is locally almost directed.
	\end{proof}
	
	Recall that if $(\mathcal{M},d)\in\mathbb{M}^{dec}(\mathcal{S},\mathcal{L})$ and $Q\in\mathcal{M}\cup\mathcal{M}[\omega_1]$, then $(\mathcal{M},d)\restr Q$ was defined as the pair $(\mathcal{M}\cap Q,d\restr Q)$, where $d\restr Q$ is the function on $\dom(\mathcal{M}\cap Q)$ defined by $d\restr Q(\varepsilon)=d(\varepsilon)\cap Q$, for every $\varepsilon\in\dom(\mathcal{M}\cap Q)$.
	
	\begin{corollary}\label{morass-refl-corollary}
		Let $(\mathcal{M},d)\in\mathbb{M}^{mor}(\mathcal{S},\mathcal{L})$, and let $Q\in\mathcal{M}\cup\mathcal{M}[\omega_1]$. Then $(\mathcal{M},d)\restr Q\in\mathbb{M}^{mor}(\mathcal{S},\mathcal{L})\cap Q$ and $(\mathcal{M},d)\leq(\mathcal{M},d)\restr Q$.
	\end{corollary}
	\begin{proof}
		This is proven exactly like Lemma \ref{deco-restr}, except that we use Lemma \ref{morass-refl} instead of Lemma \ref{amal-rest}.
	\end{proof}
	
	\begin{lemma}\label{morass-amal-iso}
		Let $\mathcal{M}_0$ and $\mathcal{M}_1$ be two morass-like $(\mathcal{S},\mathcal{L})$-symmetric systems such that $\mathcal{M}_0\cong\mathcal{M}_1$. Suppose that $R_0\in\mathcal{M}_0$ and $R_1\in\mathcal{M}_1$ are two strongly $\omega_1$-isomorphic models such that $R_0[\omega_1]\cap\omega_3$ and $R_1[\omega_1]\cap\omega_3$ witness that $\mathcal{F}_{\mathcal{M}_0}$ and $\mathcal{F}_{\mathcal{M}_1}$ are directed, respectively. Then $\mathcal{M}_0\cup\mathcal{M}_1$ is a morass-like $(\mathcal{S},\mathcal{L})$-symmetric system.
	\end{lemma}
	\begin{proof}
		By Lemma \ref{pureamalgamation2} we know that $\mathcal{M}_0\cup\mathcal{M}_1$ is an $(\mathcal{S},\mathcal{L})$-symmetric system. Moreover, $\{R_0[\omega_1]\cap\omega_3,R_1[\omega_1]\cap\omega_3\}$ is a split of $\mathcal{F}_{\mathcal{M}_0\cup\mathcal{M}_1}$, and thus $\mathcal{M}_0\cup\mathcal{M}_1$ is almost directed. Hence, we only need to check that $\mathcal{M}_0\cup\mathcal{M}_1$ is locally almost directed. Let us denote by $\Psi_{0,1}$ the isomorphism witnessing $\mathcal{M}_0\cong\mathcal{M}_1$.
		
		Let $i\in\{0,1\}$ and $j\in\{0,1\}\setminus\{i\}$, and suppose that $P_i\in\mathcal{M}_i$. We will show that $(\mathcal{M}_0\cup\mathcal{M}_1)\cap P_i[\omega_1]$ is almost directed. Let $Q_j\in\mathcal{M}_j$ such that $Q_j[\omega_1]\cap\omega_3\in\mathcal{F}_{\mathcal{M}_0\cup\mathcal{M}_1}\cap P_i[\omega_1]$. On the one hand, note that since $Q_j[\omega_1]\cap\omega_3\in\mathcal{F}_{\mathcal{M}_j}$, we have $\Psi_{j,i}(Q_j[\omega_1]\cap\omega_3)\in\mathcal{F}_{\mathcal{M}_i}$. On the other hand, since $Q_j[\omega_1]\cap\omega_3\in P_i[\omega_1]$ and $\Psi_{j,i}$ is the identity on $\bigcup\mathcal{M}_0[\omega_1]\cap\bigcup\mathcal{M}_1[\omega_1]$, we have $\Psi_{j,i}(Q_j[\omega_1]\cap\omega_3)=Q_j[\omega_1]\cap\omega_3$. Hence, $Q_j[\omega_1]\cap\omega_3\in\mathcal{F}_{\mathcal{M}_i}$. This proves that 
		\[
		\mathcal{F}_{\mathcal{M}_0\cup\mathcal{M}_1}\cap P_i[\omega_1]=\mathcal{F}_{\mathcal{M}_i}\cap P_i[\omega_1],
		\]
		and as $\mathcal{M}_i$ is a morass-like $(\mathcal{S},\mathcal{L})$-symmetric system, we can conclude that $(\mathcal{M}_0\cup\mathcal{M}_1)\cap P_i[\omega_1]$ is almost directed.
	\end{proof}
	
	Recall from Section \ref{section-decorations}, that if $(\mathcal{M}_0,d_0),(\mathcal{M}_1,d_1)\in\mathbb{M}^{dec}(\mathcal{S},\mathcal{L})$, then $(\mathcal{M}_0,d_0)\cong(\mathcal{M}_1,d_1)$ denotes that there exists a map $\Psi$ witnessing $\mathcal{M}_0\cong\mathcal{M}_1$, such that for every $\varepsilon\in\dom(\mathcal{M}_0)=\dom(\mathcal{M}_1)$ with $\varepsilon<\max\dom(\mathcal{M}_0)$, $d_1(\varepsilon)=\Psi(d_0(\varepsilon))$.
	
	\begin{corollary}\label{morass-amal-iso-corollary}
		Let $(\mathcal{M}_0,d_0),(\mathcal{M}_1,d_1)\in\mathbb{M}^{mor}(\mathcal{S},\mathcal{L})$ such that $(\mathcal{M}_0,d_0)\cong(\mathcal{M}_1,d_1)$. Suppose that $R_0\in\mathcal{M}_0$ and $R_1\in\mathcal{M}_1$ are two strongly $\omega_1$-isomorphic models such that $R_0[\omega_1]\cap\omega_3$ and $R_1[\omega_1]\cap\omega_3$ witness that $\mathcal{F}_{\mathcal{M}_0}$ and $\mathcal{F}_{\mathcal{M}_1}$ are directed, respectively. Then there is $(\mathcal{M},d)\in\mathbb{M}^{mor}(\mathcal{S},\mathcal{L})$ extending both $(\mathcal{M}_0,d_0)$ and $(\mathcal{M}_1,d_1)$.
	\end{corollary}
	\begin{proof}
		This is proven exactly like Lemma \ref{deco-amal-iso}, but using Lemma \ref{morass-amal-iso} rather than Lemma \ref{pureamalgamation2}. 
	\end{proof}
	
	\begin{proposition}\label{morass-prop-technical}
		Let $\mathcal{M}$ be a morass-like $(\mathcal{S},\mathcal{L})$-symmetric system and let $Q\in\mathcal{M}$. Let $\mathcal{W}$ be another morass-like $(\mathcal{S},\mathcal{L})$-symmetric system such that
		\begin{enumerate}
			\item $\mathcal{M}\cap Q[\omega_1]\subseteq\mathcal{W}\subseteq Q[\omega_1]$,
			
			\item $\mathcal{W}(\max\dom(\mathcal{W}))\subseteq Q$, and
			
			\item $\mathcal{W}\land\mathcal{M}\subseteq\mathcal{W}$.
		\end{enumerate}   
		Then,
		\[
		\mathcal{V}:=\mathcal{M}(\geq\varepsilon_Q)\cup\{\Psi_{Q[\omega_1],Q'[\omega_1]}(W):W\in\mathcal{W},Q'\in\mathcal{M}(\varepsilon_Q)\}
		\]
		is a morass-like $(\mathcal{S},\mathcal{L})$-symmetric system such that $\mathcal{M}\cup\mathcal{W}\subseteq\mathcal{V}$.
	\end{proposition}
	\begin{proof}
		We know from Proposition \ref{prop-pure-amalgamation-3} that $\mathcal{V}$ is an $(\mathcal{S},\mathcal{L})$-symmetric system. Let us check that it is morass-like.
		
		Since every model of $\mathcal{V}$ of maximal $\omega_2$-height is an element of $\mathcal{M}$ and $\mathcal{M}$ is morass-like, $\mathcal{V}$ is almost directed, and this is witnessed by the same models that witness the almost directedness of $\mathcal{M}$.
		
		Let us prove that $\mathcal{V}$ is locally almost directed. Let $P\in\mathcal{V}$. We need to show that $\mathcal{V}\cap P[\omega_1]$ is almost directed.
		
		\textbf{Case 1.} Suppose that $\varepsilon_P<\varepsilon_Q$. By the definition of $\mathcal{V}$, there are models $Q'\in\mathcal{M}(\varepsilon_Q)$ and $W\in\mathcal{W}$ such that $P=\Psi_{Q[\omega_1],Q'[\omega_1]}(W)$. We claim that \begin{align*}
			\Psi_{Q[\omega_1],Q'[\omega_1]}(\mathcal{F}_{\mathcal{W}}\cap W[\omega_1])&=\mathcal{F}_{\Psi_{Q[\omega_1],Q'[\omega_1]}(\mathcal{W})}\cap P[\omega_1]\\&=\mathcal{F}_{\mathcal{V}}\cap P[\omega_1].
		\end{align*}
		The first equality is clear. Let us show the second one. It is clear that $\mathcal{F}_{\Psi_{Q[\omega_1],Q'[\omega_1]}(\mathcal{W})}\cap P[\omega_1]\subseteq\mathcal{F}_{\mathcal{V}}\cap P[\omega_1]$. Hence, let $R[\omega_1]\cap\omega_3\in\mathcal{F}_\mathcal{V}\cap P[\omega_1]$. Since $\varepsilon_R<\varepsilon_P<\varepsilon_Q$, there must be some $Q''\in\mathcal{M}(\varepsilon_Q)$ and $V\in\mathcal{W}$ such that $R=\Psi_{Q[\omega_1],Q''[\omega_1]}(V)$. Since $R\in Q''[\omega_1]$ and $R[\omega_1]\cap\omega_3\in P[\omega_1]\subseteq Q'[\omega_1]$, the isomorphism $\Psi_{Q''[\omega_1],Q'[\omega_1]}$ fixes $R[\omega_1]\cap\omega_3$. Therefore,
		\begin{align*}
			R[\omega_1]\cap\omega_3&=\Psi_{Q''[\omega_1],Q'[\omega_1]}(R[\omega_1]\cap\omega_3)\\
			&=\Psi_{Q''[\omega_1],Q'[\omega_1]}(\Psi_{Q[\omega_1],Q''[\omega_1]}(V[\omega_1]\cap\omega_3))\\
			&=\Psi_{Q[\omega_1],Q'[\omega_1]}(V[\omega_1]\cap\omega_3)\in\mathcal{F}_{\Psi_{Q[\omega_1],Q'[\omega_1]}(\mathcal{W})}\cap P[\omega_1].
		\end{align*}
		Hence, the almost directedness of $\mathcal{V}\cap P[\omega_1]$ follows from the almost directedness of $\mathcal{W}\cap W[\omega_1]$ and Proposition \ref{morass-isocopy}. 
		
		\textbf{Case 2.} Suppose that $\varepsilon_P=\varepsilon_Q$. Using an argument similar to the one from the last case, it is not too hard to see that $\mathcal{F}_\mathcal{V}\cap Q[\omega_1]=\mathcal{F}_\mathcal{W}$. Therefore,
		\begin{align*}
			\mathcal{F}_\mathcal{V}\cap P[\omega_1]&=\Psi_{Q[\omega_1],P[\omega_1]}(\mathcal{F}_\mathcal{V}\cap Q[\omega_1])\\
			&=\Psi_{Q[\omega_1],P[\omega_1]}(\mathcal{F}_\mathcal{W})=\mathcal{F}_{\Psi_{Q[\omega_1],P[\omega_1]}(\mathcal{W})},
		\end{align*}
		and hence, the almost directedness of $\mathcal{V}\cap P[\omega_1]$ follows from the almost directedness of $\mathcal{W}$ and Proposition \ref{morass-isocopy}. 
		
		\textbf{Case 3.} Suppose that $\varepsilon_P>\varepsilon_Q$. By the definition of $\mathcal{V}$, note that $P\in\mathcal{M}$ and that all $\mathcal{V}$-predecessors of $P$ are also members of $\mathcal{M}$. Therefore, by Remark \ref{morass-remark}, the almost directedness of $\mathcal{V}\cap P[\omega_1]$ is witnessed by the same models that witness the almost directedness of $\mathcal{M}\cap P[\omega_1]$.
	\end{proof}

	\begin{corollary}\label{morass-prop-technical-corollary}
		Let $(\mathcal{M},d)\in\mathbb{M}^{mor}(\mathcal{S},\mathcal{L})$ and let $Q\in\mathcal{M}$. Let $(\mathcal{W},e)$ be a condition in $\mathbb{M}^{mor}(\mathcal{S},\mathcal{L})$ with the following properties:
		\begin{enumerate}
			\item $(\mathcal{W},e)\leq(\mathcal{M},d)\restr Q[\omega_1]$ and $(\mathcal{W},e)\subseteq Q[\omega_1]$,
			
			\item $\mathcal{W}(\max\dom(\mathcal{W}))\subseteq Q$,
			
			\item $e(\max\dom(\mathcal{W}))\subseteq Q$, and
			
			\item $\mathcal{W}\land\mathcal{M}\subseteq\mathcal{W}$.
		\end{enumerate}
		Then there is $(\mathcal{V},f)\in\mathbb{M}^{mor}(\mathcal{S},\mathcal{L})$ that extends both $(\mathcal{M},d)$ and $(\mathcal{W},e)$.
	\end{corollary}
	\begin{proof}
		This is proven exactly like Proposition \ref{prop-deco-amalgamation}, except that we use Proposition \ref{morass-prop-technical} instead of Proposition \ref{prop-pure-amalgamation-3}.
	\end{proof}
	
	\begin{lemma}\label{morass-amal-L}
		Let $\mathcal{M}$ be a morass-like $(\mathcal{S},\mathcal{L})$-symmetric system and let $X\in\mathcal{M}\cap\mathcal{L}$. Let $\mathcal{W}$ be another morass-like $(\mathcal{S},\mathcal{L})$-symmetric system such that $\mathcal{M}\cap X\subseteq\mathcal{W}\subseteq X$. Then there is a morass-like $(\mathcal{S},\mathcal{L})$-symmetric system $\mathcal{V}$ such that $\mathcal{M}\cup\mathcal{W}\subseteq\mathcal{V}$.
	\end{lemma}
	\begin{proof}
		By Proposition \ref{prop-pure-amalgamation-1}, we have $\mathcal{W}\land\mathcal{M}\subseteq\mathcal{W}$. Hence, from Proposition \ref{morass-prop-technical} it follows that 
		\[
		\mathcal{V}=\mathcal{M}(\geq\varepsilon_X)\cup\{\Psi_{X,Y}(W):W\in\mathcal{W},Y\in\mathcal{M}(\varepsilon_X)\}
		\]
		is a morass-like $(\mathcal{S},\mathcal{L})$-symmetric system witnessing the compatibilty of $\mathcal{M}$ and $\mathcal{W}$.
	\end{proof}	
	
	\begin{corollary}\label{morass-amal-L-corollary}
		Let $(\mathcal{M},d)\in\mathbb{M}^{mor}(\mathcal{S},\mathcal{L})$ and let $X\in\mathcal{M}\cap\mathcal{L}$. Let $(\mathcal{W},e)$ be a condition in $\mathbb{M}^{mor}(\mathcal{S},\mathcal{L})\cap X$ such that $(\mathcal{W},e)\leq(\mathcal{M},d)\restr X$. Then there is $(\mathcal{V},f)\in\mathbb{M}^{mor}(\mathcal{S},\mathcal{L})$ that extends both $(\mathcal{M},d)$ and $(\mathcal{W},e)$.
	\end{corollary}
	\begin{proof}
		This is proven exactly like Lemma \ref{deco-amal-L}, combining Corollary \ref{morass-prop-technical-corollary} and Lemma \ref{morass-amal-L}.
	\end{proof}

	\begin{lemma}\label{morass-amal-S}
		Let $\mathcal{M}$ be a morass-like $(\mathcal{S},\mathcal{L})$-symmetric system and let $M\in\mathcal{M}\cap\mathcal{S}$. Let $\mathcal{W}$ be another morass-like $(\mathcal{S},\mathcal{L})$-symmetric system such that $\mathcal{M}\cap M\subseteq\mathcal{W}\subseteq M$. Then there is a morass-like $(\mathcal{S},\mathcal{L})$-symmetric system $\mathcal{U}$ such that $\mathcal{M}\cup\mathcal{W}\subseteq\mathcal{U}$.
	\end{lemma}
	\begin{proof}
		The argument will be built upon the proof of Lemma \ref{lemma-amal-S}. Let us recall some of the notation. Let $X_n$ be a large model containing both $\mathcal{M}$ and $\mathcal{W}$, and fix a maximal $\in$-chain $X_0\in\dots\in X_{n-1}$ of large models in $\mathcal{M}\cap M$. For every $i\leq n$ denote 
		\begin{itemize}
			\item $\mathcal{M}_i=\mathcal{M}\cap X_i$,
			\item $\mathcal{W}_i=\mathcal{W}\cap X_i$, and
			\item $M_i=X_i\cap M$.
		\end{itemize}
		Observe that $\mathcal{M}_n=\mathcal{M}$, $\mathcal{W}_n=\mathcal{W}$ and $M_n=M$. Moreover, note that $\mathcal{M}_i$ and $\mathcal{W}_i$ are morass-like $(\mathcal{S},\mathcal{L})$-symmetric systems by Lemma \ref{morass-refl}, and that $\mathcal{M}_i\cap M_i\subseteq \mathcal{W}_i\subseteq M_i$. Recall from the proof of Lemma \ref{lemma-amal-S}, that we built, by induction on $i\leq n$, an $\subseteq$-increasing sequence of $(\mathcal{S},\mathcal{L})$-symmetric systems contained in $X_i$ and extending $\mathcal{M}_i\cup\mathcal{W}_i$. More precisely, for every $i\leq n$, we obtained three $(\mathcal{S},\mathcal{L})$-symmetric systems $\mathcal{V}_i\subseteq\mathcal{V}_i^*\subseteq\mathcal{U}_i$, such that $\mathcal{M}_i\cup\mathcal{W}_i\subseteq\mathcal{U}_i\subseteq X_i$. We will show that, for every $i\leq n$, the $(\mathcal{S},\mathcal{L})$-symmetric systems $\mathcal{V}_i,\mathcal{V}_i^*,\mathcal{U}_i$ are morass-like.
		
		Let $i<n$ and let $\varepsilon_i^-$ denote either $\varepsilon_{X_{i-1}}$, in case $i>0$, or an arbitrary ordinal smaller than the minimum of $\dom(\mathcal{W})$, otherwise. Recall that $\mathcal{V}_i$ was assumed to be an $(\mathcal{S},\mathcal{L})$-symmetric system with the following properties:
		\begin{enumerate}[label=(\alph*)]
			\item $\mathcal{M}_i\cap M_i[\omega_1]\subseteq\mathcal{V}_i\subseteq M_i[\omega_1]$.
			
			\item $\mathcal{W}_i\subseteq\mathcal{V}_i$.
			
			\item $\mathcal{V}_i(\geq\varepsilon_i^-)=\mathcal{W}_i(\geq\varepsilon_i^-)$.
		\end{enumerate}
		If $i=0$, we simply let $\mathcal{V}_0=\mathcal{W}_0$. In this case we need to add the extra assumption that $\mathcal{V}_i$ is morass-like.
		
		Recall that, for every $W\in\mathcal{W}_i(>\varepsilon_i^-)\cap\mathcal{L}$, the set $F_W$ was defined as 
		\[
		F_W=\{W\cap N:N\in E_i\},
		\]
		where $E_i$ is an $\in$-chain of small models $N_0\in\dots\in N_k$ of $\mathcal{M}_i\restr[\varepsilon_{M_i},\varepsilon_i^+)$ of maximal length such that $N_0=M_i$, and $\varepsilon_i^+$ is the least $\omega_2$-height of any large model in $\mathcal{M}_i(>\varepsilon_{M_i})$. Recall that $\mathcal{V}_i^*$ was the $(\mathcal{S},\mathcal{L})$-symmetric system that resulted from taking the union of $\mathcal{V}_i$ and all the sets $F_W$, where $W\in\mathcal{W}_i(>\varepsilon_i^-)\cap\mathcal{L}$. Moreover, recall that $\mathcal{V}_i^*$ coincides with the union of $\mathcal{V}_i$ and
		\[
		\mathcal{V}_i(>\varepsilon_i^-)\land\mathcal{M}_i\restr[\varepsilon_{M_i},\varepsilon_i^+).
		\]
		We need to check that $\mathcal{V}_i^*$ is morass-like.

		Since every model of $\mathcal{V}_i^*$ of maximal $\omega_2$-height is an element of $\mathcal{V}_i$, and $\mathcal{V}_i$ is morass-like, $\mathcal{V}_i^*$ is almost directed, and this is witnessed by the same models that witness the almost directedness of $\mathcal{V}_i$.
		
		Let us show that $\mathcal{V}_i^*$ is locally almost directed now. So, fix some $P\in\mathcal{V}_i^*$. We need to check that $\mathcal{V}_i^*\cap P[\omega_1]$ is almost directed. Let $W_0$ be a large model of minimal $\omega_2$-height among the models in $\mathcal{W}_i(>\varepsilon_i^-)\cap\mathcal{L}$. If $\varepsilon_P$ is strictly smaller than $\varepsilon_{W_0\cap M_i}$, then $P\in\mathcal{V}_i$ and $\mathcal{V}_i^*\cap P[\omega_1]=\mathcal{V}_i\cap P[\omega_1]$, which is almost directed by the induction hypothesis. So, let us assume that $\varepsilon_P\geq\varepsilon_{W_0\cap M_i}$. We divide the proof in two cases.
		
		\textbf{Case 1:} $P\in F_W$ for some $W\in \mathcal{W}_i(>\varepsilon_i^-)\cap\mathcal{L}$. If $P=W\cap M_i$, then $\mathcal{V}_i^*\cap P[\omega_1]$ is almost directed, and this is witnessed by the same models witnessing that $\mathcal{V}_i\cap W$ is almost directed. Indeed, note that all $\mathcal{V}_i$-predecessors of $W$ are members of $\mathcal{W}_i$, as $\mathcal{V}_i(\geq\varepsilon_i^-)=\mathcal{W}_i(\geq\varepsilon_i^-)$ by (c) of the induction hypothesis. Hence, all $\mathcal{V}_i$-predecessors of $W$ are also members of $M_i$. Therefore, since any $\mathcal{V}_i^*$-predecessor of $W\cap M_i=P$ is a $\mathcal{V}_i$-predecessor of $W$, it follows that $\mathcal{V}_i^*\cap P[\omega_1]$ is almost directed. Lastly, since $F_W$ forms an $\in$-chain, if $P\in F_W\setminus\{W\cap M_i\}$, and $W\cap N\in F_W$ is the predecessor of $P$ in $F_W$, then $(W\cap N)[\omega_1]\cap\omega_3$ witnesses that $\mathcal{F}_{\mathcal{V}_i^*}\cap P[\omega_1]$ is directed.
		
		\textbf{Case 2:} $P\in\mathcal{V}_i$. First, note that since $\mathcal{V}_i(\geq\varepsilon_i^-)=\mathcal{W}_i(\geq\varepsilon_i^-)$ by (c) of the induction hypothesis, $P\in\mathcal{W}_i$. If $P$ is a large model, by the exact same reasons as in the last case, we can argue that $\mathcal{F}_{\mathcal{V}_i^*}\cap P$ is directed, and this is witnessed by the predecessor of $P$ in the $\in$-chain $F_P$. Suppose now that $P$ is a small model. Then, by the construction of $\mathcal{V}_i^*$, the $\mathcal{V}_i^*$-predecessors of $P$ have to be members of $\mathcal{V}_i$ (and thus, also members of $\mathcal{W}_i$). Therefore, the same models that witness the almost directedness of $\mathcal{V}_i\cap P[\omega_1]$, also witness the almost directedness of $\mathcal{V}_i^*\cap P[\omega_1]$.
		
		This finishes the proof of the local almost directedness of $\mathcal{V}_i^*$, and hence $\mathcal{V}_i^*$ is a morass-like $(\mathcal{S},\mathcal{L})$-symmetric system.
		
		Recall (from the proof of Lemma \ref{lemma-amal-S}) that $\mathcal{V}_i^*\land\mathcal{M}_i\subseteq\mathcal{V}_i^*$ by Proposition \ref{prop-pure-amalgamation-2}. Therefore, by Proposition \ref{morass-prop-technical},
		\[
		\mathcal{U}_i:=\mathcal{M}_i(\geq\varepsilon_{M_i})\cup\{\Psi_{M_i[\omega_1],N_i[\omega_1]}(V):V\in\mathcal{V}_i^*,N_i\in\mathcal{M}_i(\varepsilon_{M_i})\}
		\] 
		is a morass-like $(\mathcal{S},\mathcal{L})$-symmetric system such that $\mathcal{V}_i^*\cup\mathcal{M}_i\subseteq\mathcal{U}_i\subseteq X_i$.
		
		If $i=n$, we have obtained a morass-like $(\mathcal{S},\mathcal{L})$-symmetric system $\mathcal{U}_n$ which extends $\mathcal{M}_n=\mathcal{M}$ and $\mathcal{W}_n=\mathcal{W}$, and we are done.
		
		If $i<n$, we need to make sure that the induction can go through by extending $\mathcal{U}_i$ to a morass-like $(\mathcal{S},\mathcal{L})$-symmetric system $\mathcal{V}_{i+1}$, so that $\mathcal{V}_{i+1}$ satisfies the conditions (a)-(c) that we assumed for $\mathcal{V}_i$ at the beginning of the inductive construction of the morass-like $(\mathcal{S},\mathcal{L})$-symmetric systems $\mathcal{U}_i$. Since $\mathcal{W}_{i+1}\cap X_i=\mathcal{W}_i\subseteq\mathcal{U}_i\subseteq X_i$, we have $\mathcal{U}_i\land\mathcal{W}_{i+1}\subseteq\mathcal{U}_i$ by Proposition \ref{prop-pure-amalgamation-1}. Therefore,  
		\[
		\mathcal{V}_{i+1}:=\mathcal{W}_{i+1}(\geq\varepsilon_{X_i})\cup\{\Psi_{X_i,Y_i}(U):U\in\mathcal{U}_i,Y_i\in\mathcal{W}_{i+1}(\varepsilon_{X_i})\}
		\]
		is a morass-like $(\mathcal{S},\mathcal{L})$-symmetric system by Proposition \ref{morass-prop-technical}, and as we saw in the proof of Lemma \ref{lemma-amal-S}, it has the following properties:
		\begin{enumerate}[label=(\alph*)]
			\item $\mathcal{M}_{i+1}\cap M_{i+1}[\omega_1]\subseteq\mathcal{V}_{i+1}\subseteq M_{i+1}[\omega_1]$.
			
			\item $\mathcal{W}_{i+1}\subseteq\mathcal{V}_{i+1}$.
			
			\item $\mathcal{V}_{i+1}(\geq\varepsilon_{X_i})=\mathcal{W}_{i+1}(\geq\varepsilon_{X_i})$.
		\end{enumerate}
		So, the induction can go through as we wanted.
	\end{proof}
	
	\begin{corollary}\label{morass-amal-S-corollary}
		Let $(\mathcal{M},d)\in\mathbb{M}^{mor}(\mathcal{S},\mathcal{L})$ and let $M\in\mathcal{M}\cap\mathcal{S}$. Let $(\mathcal{W},e)$ be a condition in $\mathbb{M}^{mor}(\mathcal{S},\mathcal{L})\cap M$ such that $(\mathcal{W},e)\leq(\mathcal{M},d)\restr M$. Then there is $(\mathcal{U},g)\in\mathbb{M}^{mor}(\mathcal{S},\mathcal{L})$ extending both $(\mathcal{M},d)$ and $(\mathcal{W},e)$.
	\end{corollary}
	\begin{proof}
		The proof is exactly the same as the one for Lemma \ref{deco-amal-S}, except that we need to appeal to the preceding amalgamation lemmas rather than the ones obtained in Subsection \ref{subsection-deco-amal-lemmas}.
	\end{proof}

	\subsection{Preservation lemmas}
			
	The following result is proven exactly like theorems \ref{pureproperS} and \ref{pureproperL} (and Theorem \ref{deco-preservation}), using the amalgamation lemmas from the last subsection.
	
	\begin{theorem}
		The forcing $\mathbb{M}^{mor}(\mathcal{S},\mathcal{L})$ is strongly $\mathcal{S}$-proper and strongly $\mathcal{L}$-proper.
	\end{theorem}
	
	Moreover, as in Subsection \ref{subsection-pure-preservation}, the following result can be isolated from the proof of the last theorem.
	
	\begin{lemma}\label{morass-genericity}
		Every $(\mathcal{M},d)\in\mathbb{M}^{mor}(\mathcal{S},\mathcal{L})$ is strongly $(Q,\mathbb{M}^{mor}(\mathcal{S},\mathcal{L}))$-generic for every $Q\in\mathcal{M}$.
	\end{lemma} 
	
	\begin{theorem}\label{morass-cc}
		If $2^{\aleph_1}=\aleph_2$, then $\mathbb{M}^{mor}(\mathcal{S},\mathcal{L})$ has the $\aleph_3$-Knaster condition. 
	\end{theorem}
	\begin{proof}
		For every $\alpha<\omega_3$, let $(\mathcal{M}_\alpha,d_\alpha)\in\mathbb{M}^{mor}(\mathcal{S},\mathcal{L})$ and $M_\alpha\in\mathcal{S}$ such that $\alpha,(\mathcal{M}_\alpha,d_\alpha)\in M_\alpha$. By Fodor's Lemma, there are a stationary set $I_0\subseteq S_{\omega_2}^{\omega_3}$ and an ordinal $\gamma<\omega_3$, such that $\sup(M_\alpha[\omega_1]\cap\alpha)=\gamma$, for every $\alpha\in I_0$. Moreover, as $2^{\aleph_1}=\aleph_2$, we have $$\aleph_2^{\aleph_1}=(2^{\aleph_1})^{\aleph_1}=2^{\aleph_1\cdot\aleph_1}=\aleph_2<\aleph_3.$$ Thus, since $\gamma<\omega_3$, there is a subset $I_1\subseteq I_0$ of size $\aleph_3$ such that for any two $\alpha,\beta\in I_1$, $M_\alpha[\omega_1]\cap\alpha=M_\beta[\omega_1]\cap\beta$. Note that for every $\alpha\in I_1$, there are $\aleph_3$-many $\beta\in I_1$ greater than $\alpha$ such that $M_\alpha[\omega_1]\cap\omega_3\subseteq\beta$. Hence, there is a subset $I_2\subseteq I_1$ of size $\aleph_3$ such that for any two $\alpha<\beta$ in $I_2$, $M_\alpha[\omega_1]\cap\omega_3\subseteq\beta$. Furthermore, again by $2^{\aleph_1}=\aleph_2$, there is a subset $I_3\subseteq I_2$ of size $\aleph_3$ such that $\{M_\alpha[\omega_1]:\alpha\in I_3\}$ forms a $\Delta$-system with root $\mathcal{X}$. Since there are only $\aleph_2$-many isomorphism types between structures of the form $(M_\alpha[\omega_1];\in,M_\alpha,\mathcal{M}_\alpha,d_\alpha,\alpha,\mathcal{X})$, there must be a subset $I_4\subseteq I_3$ of size $\aleph_3$, such that for any two $\alpha,\beta\in I_4$, the structures
		\[
		(M_\alpha[\omega_1];\in,M_\alpha,\mathcal{M}_\alpha,d_\alpha,\alpha,\mathcal{X},\lhd)
		\]
		and
		\[
		(M_\beta[\omega_1];\in,M_\beta,\mathcal{M}_\beta,d_\beta,\beta,\mathcal{X},\lhd)
		\]
		are isomorphic via an isomorphism $\Psi$ which is the identity on $\mathcal{X}$, and is such that $\Psi(\alpha)=\beta$. Therefore, by Proposition \ref{morass-strongly-iso}, $\{M_\alpha:\alpha\in I_4\}$ is a set of pairwise strongly $\omega_1$-isomorphic models.
				
		For every $\alpha\in I_4$, let $(\mathcal{M}_{M_\alpha},d_{M_\alpha})$ denote the condition that extends $(\mathcal{M}_\alpha,d_\alpha)$ such that $M_\alpha\in\mathcal{M}_{M_\alpha}$, given by Corollary \ref{morass-on-top-corollary}. Note that for every $\alpha\in I_4$, $M_\alpha[\omega_1]\cap\omega_3$ witnesses that $\mathcal{F}_{\mathcal{M}_{M_\alpha}}$ is directed. Therefore, since for any two ordinals $\alpha<\beta$ in $I_4$, we have
		\[
		(\mathcal{M}_{M_\alpha},d_{M_\alpha})\cong(\mathcal{M}_{M_\beta},d_{M_\beta})
		\]
		via the strong $\omega_1$-isomorphism $\Psi_{M_\alpha[\omega_1],M_\beta[\omega_1]}$, Corollary \ref{morass-amal-iso-corollary} ensures that the conditions $(\mathcal{M}_{M_\alpha},d_{M_\alpha})$ and $(\mathcal{M}_{M_\beta},d_{M_\beta})$ are compatible, and this in turn implies the compatibility of the conditions $(\mathcal{M}_\alpha,d_\alpha)$ and $(\mathcal{M}_\beta,d_\beta)$.
	\end{proof}
		
	Once again, in light of Lemmas \ref{preservacio-proper} and \ref{preservation-cc}, from the last theorems we can deduce the following preservation theorem.
	
	\begin{corollary}
		If $2^{\aleph_1}=\aleph_2$ holds and $\mathcal{L}$ is stationary in $H(\kappa)$, then $\mathbb{M}^{mor}(\mathcal{S},\mathcal{L})$ preserves all cardinals.	
	\end{corollary}
	
	\begin{theorem}
		If $\mathcal{L}$ is stationary in $H(\kappa)$, then $\mathbb{M}^{mor}(\mathcal{S},\mathcal{L})$ preserves $2^{\aleph_1}=\aleph_2$.
	\end{theorem}
	\begin{proof}
		The argument is exactly the same as the one from the proof of Theorem \ref{purepreservationCH}. The only difference is that we need to ensure that the $(\mathcal{S},\mathcal{L})$-symmetric system $\mathcal{M}_{\alpha,\beta}$ is morass-like. We recall the first half of the argument for the sake of completeness.
		
		Assume that $2^{\aleph_1}=\aleph_2$ holds and let $\langle\tau_\alpha:\alpha<\omega_3\rangle$ be a sequence of $\mathbb{M}^{mor}(\mathcal{S},\mathcal{L})$-names for subsets of $\omega_1$. Suppose that $(\mathcal{M},d)\in\mathbb{M}^{mor}(\mathcal{S},\mathcal{L})$ is a condition forcing $\langle\tau_\alpha:\alpha<\omega_3\rangle$ to be a sequence of pairwise different subsets of $\omega_1$. For each $\alpha<\omega_3$, let $X_\alpha^*$ be an elementary submodel of a large enough $H(\theta)$ such that $\mathbb{M}^{mor}(\mathcal{S},\mathcal{L}),\alpha,\tau_\alpha,\mathcal{M},d\in X_\alpha^*$ and $X_\alpha:=X_\alpha^*\cap H(\kappa)\in\mathcal{L}$. Arguing as in the proof of Theorem \ref{morass-cc}, we may assume that there are different $\alpha,\beta<\omega_3$ for which the structures $(X_\alpha^*;\in,\tau_\alpha,\alpha,\lhd)$ and $(X_\beta^*;\in,\tau_\beta,\beta,\lhd)$ are isomorphic via an isomorphism $\Psi$ that sends $\tau_\alpha$ to $\tau_\beta$ and has the following properties:
		\begin{itemize}
			\item $\Psi$ is the identity on $X_\alpha^*\cap X_\beta^*$.
			
			\item $X_\alpha^*\cap X_\beta^*\cap\omega_3$ is a proper initial segment of both $X_\alpha^*\cap\omega_3$ and $X_\beta^*\cap\omega_3$.
			
			\item Either
			\[
			(X_\alpha^*\setminus X_\beta^*)\cap\omega_3<(X_\beta^*\setminus X_\alpha^*)\cap\omega_3,
			\]
			or
			\[
			(X_\beta^*\setminus X_\alpha^*)\cap\omega_3<(X_\alpha^*\setminus X_\beta^*)\cap\omega_3.
			\]
		\end{itemize}
		Let $\mathcal{M}_{\alpha,\beta}:=\mathcal{M}\cup\{X_\alpha,X_\beta\}$, which is easily seen to be an $(\mathcal{S},\mathcal{L})$-symmetric system. In order to show that it is morass-like, it suffices to observe that the pair $\{X_\alpha\cap\omega_3,X_\beta\cap\omega_3\}$ witnesses the almost directedness of $\mathcal{M}_{\alpha,\beta}$. The local almost directedness follows from the fact that $\mathcal{M}$ is morass-like. Moreover, define a map $d_{\alpha,\beta}$ on $\dom(\mathcal{M}_{\alpha,\beta})$ such that $d_{\alpha,\beta}(\varepsilon)=d(\varepsilon)$, for each $\varepsilon\in\dom(\mathcal{M})$, and $d_{\alpha,\beta}(\varepsilon_{X_\alpha})=\emptyset$ (or any finite subset of $H(\kappa)$). Hence, $(\mathcal{M}_{\alpha,\beta},d_{\alpha,\beta})\in\mathbb{M}^{mor}(\mathcal{S},\mathcal{L})$. The rest of the proof is exactly the same as the second half of the proof of Theorem \ref{purepreservationCH}.
	\end{proof}
		
	\subsection{Application: Forcing a simplified $(\omega_2,1)$-morass}
	
	Assume for the rest of the section that $\mathcal{L}$ is stationary in $H(\kappa)$ and $2^{\aleph_1}=\aleph_2$ in $V$. Let $G$ be an $\mathbb{M}^{mor}(\mathcal{S},\mathcal{L})$-generic filter over $V$. Let $\mathcal{M}_G^{mor}$ denote the set $$\{Q\in\mathcal{M}:(\mathcal{M},d)\in G\},$$ and let $\mathcal{F}_G$ denote $\mathcal{F}_{\mathcal{M}_G^{mor}}$. That is,
	\[
	\mathcal{F}_G=\{Q[\omega_1]\cap\omega_3:Q\in\mathcal{M}_G^{mor}\}.
	\]

	\begin{proposition}\label{morass-well-founded}
		$(\mathcal{F}_G,\subseteq)$ is well-founded.
	\end{proposition}
	\begin{proof}
		Follows from Proposition \ref{morass-rank} and the fact that there cannot be infinite descending sequences of ordinals.
	\end{proof}

	\begin{proposition}\label{morass-prop-rank}
		For every $Q,P\in\mathcal{M}_G^{mor}$, the following hold:
		\begin{itemize}
			\item $\rank(Q[\omega_1]\cap\omega_3)<\rank(P[\omega_1]\cap\omega_3)$ iff $\varepsilon_Q<\varepsilon_P$.
			\item $\rank(Q[\omega_1]\cap\omega_3)=\rank(P[\omega_1]\cap\omega_3)$ iff $\varepsilon_Q=\varepsilon_P$.
		\end{itemize}
	\end{proposition}
	\begin{proof}
		First, if $\varepsilon_Q=\varepsilon_P$, then $\rank(Q[\omega_1]\cap\omega_3)=\rank(P[\omega_1]\cap\omega_3)$. To observe this implication, suppose $\varepsilon_Q=\varepsilon_P$. Then $(Q[\omega_1];\in,Q,\lhd)$ and $(P[\omega_1];\in,Q,\lhd)$ are isomorphic via $\Psi_{Q[\omega_1],P[\omega_1]}$. By the symmetry of $\mathcal{M}_G^{mor}$, it follows that $(\mathcal{M}_G^{mor}\cap Q[\omega_1];\in)$ and $(\mathcal{M}_G^{mor}\cap P[\omega_1];\in)$ are isomorphic. Similarly, $(\mathcal{F}_{\mathcal{M}_G^{mor}\cap Q[\omega_1]};\subseteq)$ and $(\mathcal{F}_{\mathcal{M}_G^{mor}\cap P[\omega_1]};\subseteq)$ are isomorphic. Hence, $\rank(Q[\omega_1]\cap\omega_3)=\rank(P[\omega_1]\cap\omega_3)$.
		
		Next, if $\varepsilon_Q<\varepsilon_P$, then $\rank(Q[\omega_1]\cap\omega_3)<\rank(P[\omega_1]\cap\omega_3)$. To observe this implication, suppose $\varepsilon_Q<\varepsilon_P$. Then by the shoulder axiom for $\mathcal{M}_G^{mor}$, there exists $P'\in\mathcal{M}_G^{mor}$ such that $\varepsilon_{P'}=\varepsilon_P$ and $Q\in P'[\omega_1]$. Therefore,
		\[
		\rank(Q[\omega_1]\cap\omega_3)<\rank(P'[\omega_1]\cap\omega_3)=\rank(P[\omega_1]\cap\omega_3).
		\]
	\end{proof}

	\begin{theorem}
		$\mathcal{F}_G$ is a simplified $(\omega_2,1)$-morass.
	\end{theorem}
	\begin{proof}
		
		The well-foundedness of $(\mathcal{F}_G,\subseteq)$ follows from Proposition \ref{morass-well-founded}.
		
		
		To see that $\mathcal{F}_G$ is locally small, we only need to note that, for every $Q\in\mathcal{M}_G^{mor}$, the size of $Q[\omega_1]$ is $\aleph_1$, and that by Proposition \ref{morass-restr},
		\[
		\mathcal{F}_G\cap Q[\omega_1]=\mathcal{F}_G\restr(Q[\omega_1]\cap\omega_3).
		\]

		The next step is to prove that $\mathcal{F}_G$ is homogeneous. Fix $Q,P\in\mathcal{M}_G^{mor}$ such that $\rank(Q[\omega_1]\cap\omega_3)=\rank(P[\omega_1]\cap\omega_3)$. Then, by Proposition \ref{morass-prop-rank}, $\varepsilon_Q=\varepsilon_P$. Hence, in light of Proposition \ref{morass-restr},
		\begin{align*}
			\mathcal{F}_G\restr(P[\omega_1]\cap\omega_3)&=\mathcal{F}_G\cap P[\omega_1]\\
			&=\Psi_{Q[\omega_1],P[\omega_1]}"(\mathcal{F}_G\cap Q[\omega_1])\\
			&=\{\Psi_{Q[\omega_1],P[\omega_1]}(R[\omega_1]\cap\omega_3):R[\omega_1]\cap\omega_3\in\mathcal{F}_G\cap Q[\omega_1]\}\\
			&=\{\Psi_{Q[\omega_1],P[\omega_1]}(R[\omega_1]\cap\omega_3):R[\omega_1]\cap\omega_3\in\mathcal{F}_G\restr(Q[\omega_1]\cap\omega_3)\}.
		\end{align*}
		
		The directedness of $\mathcal{F}_G$ follows from the fact that for every $x\in H(\kappa)^V$, the set
		\[
		D_x:=\{(\mathcal{M},d)\in\mathbb{M}^{mor}(\mathcal{S},\mathcal{L}):\exists Q\in\mathcal{M}(x\in Q)\}
		\]
		is dense in $\mathbb{M}^{mor}(\mathcal{S},\mathcal{L})$. Therefore, for any pair of models $\{P_0,P_1\}\subseteq\mathcal{M}_G^{mor}$, there is $Q\in\mathcal{M}_G^{mor}$ such that $\{P_0,P_1\}\in Q$, which in particular implies
		\[	P_0[\omega_1]\cap\omega_3,P_1[\omega_1]\cap\omega_3\subsetneq Q[\omega_1]\cap\omega_3.
		\]
		
		
		Let us now show that $\mathcal{F}_G$ is locally almost directed. Fix $Q\in\mathcal{M}_G^{mor}$. We need to check that $\mathcal{F}_G\restr(Q[\omega_1]\cap\omega_3)$ is almost directed. If $\varepsilon_Q$ is the minimum of $\dom(\mathcal{M}_G^{mor})$, then the almost directedness follows trivially. So, suppose that $\varepsilon_Q$ is not the minimum of $\dom(\mathcal{M}_G^{mor})$. We divide the proof into two cases.
				
		\textbf{Case 1.} Assume first that $\varepsilon_Q$ is a limit point of $\dom(\mathcal{M}_G^{mor})$. We claim that $\mathcal{F}_G\restr (Q[\omega_1]\cap\omega_3)$ is directed. Let $P_0,P_1\in\mathcal{M}_G^{mor}\cap Q[\omega_1]$ and $(\mathcal{M},d)\in G$ be such that $Q,P_0,P_1\in\mathcal{M}$. We claim that there is $P\in\mathcal{M}_G^{mor}\cap Q[\omega_1]$ such that $P_0,P_1\in P[\omega_1]$, which in turn implies that
		\[
		P_0[\omega_1]\cap\omega_3,P_1[\omega_1]\cap\omega_3\subseteq P[\omega_1]\cap\omega_3,
		\]
		and by appealing to Proposition \ref{morass-restr}, this is enough to get the conclusion.
		
		Let $\varepsilon'$ be the immediate predecessor of $\varepsilon_Q$ in $\dom(\mathcal{M})$. For $i=0,1$, if $\varepsilon_{P_i}<\varepsilon'$, use Proposition \ref{prop13} to find $P_i'\in\mathcal{M}(\varepsilon')$ such that $P_i\in P_i'[\omega_1]$ and $P_i'\in Q$; otherwise, let $P_i'=P_i$.

				
		\begin{claim}
			For every $x\in Q$, the set
			\[
			E_x=\{(\mathcal{N},f)\in\mathbb{M}^{mor}(\mathcal{S},\mathcal{L}):\exists R\in\mathcal{N}\big(\varepsilon_R<\varepsilon_Q\land x\in f(\varepsilon_R)\big)\}
			\]
			is open and dense below $(\mathcal{M},d)$.
		\end{claim}
		\begin{proof}
			Suppose that $(\mathcal{N}_0,f_0)\in\mathbb{M}^{mor}(\mathcal{S},\mathcal{L})$ extends $(\mathcal{M},d)$. Let $R\in\mathcal{N}_0$ such that $\varepsilon_R$ is the predecessor of $\varepsilon_Q$ in $\dom(\mathcal{N}_0)$. Let $\mathcal{N}=\mathcal{N}_0$, and let $f$ be the function on $\dom(\mathcal{N})$ defined exactly like $f_0$, except that
			\[
			f(\varepsilon_R)=f_0(\varepsilon_R)\cup\{\Psi_{Q[\omega_1],Q'[\omega_1]}(x):Q'\in\mathcal{N}_0(\varepsilon_Q)\}.
			\]
			Since $x\in Q$, it is not too hard to see that $f$ is a decoration on $\mathcal{N}$ and $x\in f(\varepsilon_R)$. Therefore, $(\mathcal{N},f)$ is a condition in $E_x$ extending $(\mathcal{N}_0,f_0)$.
		\end{proof}

		Since $E_{\{P_0',P_1'\}}$ is open and dense below $(\mathcal{M},d)$ by the last claim, and $(\mathcal{M},d)\in G$, there is $(\mathcal{N},f)\in G\cap E_{\{P_0',P_1'\}}$ such that $(\mathcal{N},f)\leq(\mathcal{M},d)$. Therefore, there is $R\in\mathcal{N}$ such that $\varepsilon_R<\varepsilon_Q$ and $\{P_0',P_1'\}\in f(\varepsilon_R)$.
		
		Since $\varepsilon_Q$ is a limit point of $\dom(\mathcal{M}_G^{mor})$, there is $(\mathcal{W},e)\leq(\mathcal{N},f)$ in $G$ with $\varepsilon\in\dom(\mathcal{W})$ such that $\varepsilon_R<\varepsilon<\varepsilon_Q$. Hence, since  $\{P_0',P_1'\}\in e(\varepsilon_R)$, there must be a model $S\in\mathcal{W}$ such that $\varepsilon_S$ is the immediate successor of $\varepsilon_R$ in $\dom(\mathcal{W})$ and $\{P_0',P_1'\}\in S$. Use the shoulder axiom to find $Q'\in\mathcal{W}$ such that $\varepsilon_{Q'}=\varepsilon_Q$ and $S\in Q'[\omega_1]$, and let $P:=\Psi_{Q'[\omega_1],Q[\omega_1]}(S)$, which is a member of $\mathcal{W}$ by symmetry. Since $\{P_0',P_1'\}\in S\subseteq Q'[\omega_1]$ and $\{P_0',P_1'\}\in Q$, the isomorphism $\Psi_{Q'[\omega_1],Q[\omega_1]}$ fixes the pair $\{P_0',P_1'\}$, and hence 
		\[
		\{P_0',P_1'\}\in P\in Q[\omega_1].
		\]
		Therefore, we can conclude that $P_0,P_1\in P[\omega_1]$, as we claimed.\footnote{This case is the only point in the whole section where the decorations seem to be necessary for the argument to go through.}
		

		
		\textbf{Case 2.} Suppose now that $\varepsilon_Q$ is a successor in $\dom(\mathcal{M}_G^{mor})$ and let $\varepsilon_0$ be its immediate predecessor. Let $(\mathcal{M},d)\in G$ be such that $Q\in\mathcal{M}$ and $\varepsilon_0\in\dom(\mathcal{M})$.
		
		\textbf{Subcase 2.1.} Assume first that $\mathcal{F}_\mathcal{M}\restr(Q[\omega_1]\cap\omega_3)$ splits. So, by Remark \ref{morass-remark}, we can find models $P_0,P_1\in\mathcal{M}\cap Q$ such that $\{P_0[\omega_1]\cap\omega_3,P_1[\omega_1]\cap\omega_3\}$ is a split of $\mathcal{F}_\mathcal{M}\restr(Q[\omega_1]\cap\omega_3)$ and $\varepsilon_0=\varepsilon_{P_0}=\varepsilon_{P_1}$. We claim that the pair $\{P_0[\omega_1]\cap\omega_3,P_1[\omega_1]\cap\omega_3\}$ is a split of $\mathcal{F}_G\restr(Q[\omega_1]\cap\omega_3)$. To do so, we will prove that, for every $(\mathcal{W},e)\leq (\mathcal{M},d)$ in $G$, the pair $\{P_0[\omega_1]\cap\omega_3,P_1[\omega_1]\cap\omega_3\}$ remains a split of $\mathcal{F}_\mathcal{W}\restr(Q[\omega_1]\cap\omega_3)$.
		
		Let us first argue that this does in fact show that $\{P_0[\omega_1]\cap\omega_3,P_1[\omega_1]\cap\omega_3\}$ is a split of $\mathcal{F}_G\restr(Q[\omega_1]\cap\omega_3)$. That is, we need to check the following three items (see clause (5b) in Definition \ref{def-simplified-morass}):
		\begin{itemize}
			\item $\rank(P_0[\omega_1]\cap\omega_3)=\rank(P_1[\omega_1]\cap\omega_3)$.
			
			\item $P_0[\omega_1]\cap P_1[\omega_1]\cap\omega_3$ is a proper initial segment of $P_0[\omega_1]\cap\omega_3$ and $P_1[\omega_1]\cap\omega_3$.		
			
			\item Either
			 \begin{align*}
			 P_0[\omega_1]\cap P_1[\omega_1]\cap\omega_3&<(P_0[\omega_1]\setminus P_1[\omega_1])\cap\omega_3\\
			 &<(P_1[\omega_1]\setminus P_0[\omega_1])\cap\omega_3,
			 \end{align*}
			 or 
			 \begin{align*}
				P_0[\omega_1]\cap P_1[\omega_1]\cap\omega_3&<(P_1[\omega_1]\setminus P_0[\omega_1])\cap\omega_3\\
				&<(P_0[\omega_1]\setminus P_1[\omega_1])\cap\omega_3.
			\end{align*}
			 
			\item $\mathcal{F}_G\restr(Q[\omega_1]\cap\omega_3)$ equals the set
			\begin{align*}
				\{P_0[\omega_1]\cap\omega_3,P_1[\omega_1]\cap\omega_3\}&\cup\mathcal{F}_G\restr (P_0[\omega_1]\cap\omega_3)\\
				&\cup\mathcal{F}_G\restr(P_1[\omega_1]\cap\omega_3).
			\end{align*}
		\end{itemize}
		The first item follows from Proposition \ref{morass-prop-rank} and the fact that $\varepsilon_{P_0}=\varepsilon_{P_1}$. The second and third items follow from $\{P_0[\omega_1]\cap\omega_3,P_1[\omega_1]\cap\omega_3\}$ being a split of $\mathcal{F}_{\mathcal{M}}\restr(Q[\omega_1]\cap\omega_3)$. The fourth item follows easily.		
		
		Let us now return to the proof that for every $(\mathcal{W},e)\in G$ extending $(\mathcal{M},d)$, the pair $\{P_0[\omega_1]\cap\omega_3,P_1[\omega_1]\cap\omega_3\}$ is a split of $\mathcal{F}_\mathcal{W}\restr(Q[\omega_1]\cap\omega_3)$. Fix a condition $(\mathcal{W},e)$ in $G$ such that $(\mathcal{W},e)\leq(\mathcal{M},d)$. Since $P_0,P_1\in\mathcal{W}$ and $\varepsilon_0=\varepsilon_{P_0}=\varepsilon_{P_1}$, it is clear that $\mathcal{F}_\mathcal{W}\restr(Q[\omega_1]\cap\omega_3)$ cannot be directed with respect to inclusion. Hence, there must be some $R_0,R_1\in\mathcal{W}$ such that $\varepsilon_{R_0}=\varepsilon_{R_1}=\varepsilon_0$ and $\{R_0[\omega_1]\cap\omega_3,R_1[\omega_1]\cap\omega_3\}$ is a split of $\mathcal{F}_\mathcal{W}\restr(Q[\omega_1]\cap\omega_3)$. So, by definition,
		\begin{align*}
			\mathcal{F}_\mathcal{W}\restr(Q[\omega_1]\cap\omega_3)=&\mathcal{F}_\mathcal{W}\restr(R_0[\omega_1]\cap\omega_3)\cup\mathcal{F}_\mathcal{W}\restr(R_1[\omega_1]\cap\omega_3)\\
			&\cup\{R_0[\omega_1]\cap\omega_3,R_1[\omega_1]\cap\omega_3\}.
		\end{align*}
		Note that as $\varepsilon_0=\varepsilon_{P_i}=\varepsilon_{R_i}$, for $i=0,1$, we must have  
		\[
		P_0[\omega_1]\cap\omega_3,P_1[\omega_1]\cap\omega_3\in\{R_0[\omega_1]\cap\omega_3,R_1[\omega_1]\cap\omega_3\},
		\]
		which implies that $\{P_0[\omega_1]\cap\omega_3,P_1[\omega_1]\cap\omega_3\}$ is a split of $\mathcal{F}_\mathcal{W}\restr(Q[\omega_1]\cap\omega_3)$, as we wanted.
				
		\textbf{Subcase 2.2.} Assume now that $\mathcal{F}_\mathcal{M}\restr(Q[\omega_1]\cap\omega_3)$ is directed with respect to inclusion. If there is an extension $(\mathcal{N},f)$ of $(\mathcal{M},d)$ in $G$ such that $\mathcal{F}_{\mathcal{N}}\restr(Q[\omega_1]\cap\omega_3)$ splits, then by the last subcase, for any further extension $(\mathcal{W},e)\in G$ of $(\mathcal{N},f)$, $\mathcal{F}_{\mathcal{W}}\restr(Q[\omega_1]\cap\omega_3)$ will remain splitting, and hence $\mathcal{F}_G\restr(Q[\omega_1]\cap\omega_3)$ will also split. So, let us assume that there is no extension $(\mathcal{N},f)$ of $(\mathcal{M},d)$ in $G$ such that $\mathcal{F}_{\mathcal{N}}\restr(Q[\omega_1]\cap\omega_3)$ splits. Or, in other words, assume that for every extension $(\mathcal{N},f)$ of $(\mathcal{M},d)$ in $G$, $\mathcal{F}_{\mathcal{N}}\restr(Q[\omega_1]\cap\omega_3)$ is directed with respect to inclusion. By Remark \ref{morass-remark} and Proposition \ref{morass-prop}, there is some $P\in\mathcal{M}\cap Q$ such that $\varepsilon_P=\varepsilon_0$ and
		\[
		\mathcal{F}_\mathcal{M}\restr(Q[\omega_1]\cap\omega_3)=\{P[\omega_1]\cap\omega_3\}\cup\mathcal{F}_\mathcal{M}\restr(P[\omega_1]\cap\omega_3).
		\]
		We will show that $\mathcal{F}_G\restr(Q[\omega_1]\cap\omega_3)$ is directed by showing that
		\[
		\mathcal{F}_G\restr(Q[\omega_1]\cap\omega_3)=\{P[\omega_1]\cap\omega_3\}\cup\mathcal{F}_G\restr(P[\omega_1]\cap\omega_3).
		\]
		Fix an arbitrary condition $(\mathcal{W},e)\in G$ such that $(\mathcal{W},e)\leq(\mathcal{M},d)$. By assumption, $\mathcal{F}_{\mathcal{W}}\restr(Q[\omega_1]\cap\omega_3)$ must be directed with respect to inclusion. So, there must be some $R\in\mathcal{W}\cap Q$ such that $\varepsilon_R=\varepsilon_0$ and 
		\[
		\mathcal{F}_\mathcal{W}\restr(Q[\omega_1]\cap\omega_3)=\{R[\omega_1]\cap\omega_3\}\cup\mathcal{F}_\mathcal{W}\restr(R[\omega_1]\cap\omega_3).
		\]
		Note that $P[\omega_1]\cap\omega_3\in\mathcal{F}_\mathcal{W}\restr(Q[\omega_1]\cap\omega_3)$. However, since $\varepsilon_P=\varepsilon_0=\varepsilon_R$, we must have $P[\omega_1]\cap\omega_3=R[\omega_1]\cap\omega_3$. Therefore, we can conclude that
		\[
		\mathcal{F}_\mathcal{W}\restr(Q[\omega_1]\cap\omega_3)=\{P[\omega_1]\cap\omega_3\}\cup\mathcal{F}_\mathcal{W}\restr(P[\omega_1]\cap\omega_3),
		\]
		which in turn implies, since $\mathcal{W}$ was chosen arbitrarily, that 
		\[
		\mathcal{F}_G\restr(Q[\omega_1]\cap\omega_3)=\{P[\omega_1]\cap\omega_3\}\cup\mathcal{F}_G\restr(P[\omega_1]\cap\omega_3),
		\]
		as we wanted. This finishes the proof of the local almost directedness of $\mathcal{F}_G$.
		
		Lastly, we need to check that $\mathcal{F}_G$ covers $\omega_3$. But this follows immediately from the fact that $\mathbb{M}^{mor}(\mathcal{S},\mathcal{L})$ preserves all cardinals and that $\mathcal{M}_G^{mor}$ covers $H(\kappa)^V$, by Corollary  \ref{morass-on-top-corollary}.
	\end{proof}

	\bibliographystyle{plain} 
	\bibliography{References}

@incollection {Abraham2010Properforcing,
	AUTHOR = {Abraham, Uri},
	TITLE = {Proper forcing},
	BOOKTITLE = {Handbook of set theory. {V}ols. 1, 2, 3},
	PAGES = {333--394},
	PUBLISHER = {Springer, Dordrecht},
	YEAR = {2010},
	ISBN = {978-1-4020-4843-2},
	MRCLASS = {03E40 (03E35)},
	MRNUMBER = {2768684},
	MRREVIEWER = {Martin\ Goldstern and Jakob\ Kellner},
	DOI = {10.1007/978-1-4020-5764-9\_6},
	URL = {https://doi.org/10.1007/978-1-4020-5764-9_6},
}

@article{AbrahamCummingsSmyth2007SomeresultsinpolychromaticsRamseytheory,
	title={Some results in polychromatic Ramsey theory},
	author={Abraham, Uri and Cummings, James and Smyth, Clifford},
	journal={The Journal of Symbolic Logic},
	volume={72},
	number={3},
	pages={865--896},
	year={2007},
	publisher={Cambridge University Press}
}

@article{AbrahamCummings2012moreresultsinpolychromatic,
	title={More results in polychromatic {Ramsey} theory},
	author={Abraham, Uri and Cummings, James},
	journal={Open Mathematics},
	volume={10},
	number={3},
	pages={1004--1016},
	year={2012},
	publisher={De Gruyter Open Access}
}

@article {Aspero:Mota2015a,
	AUTHOR = {Asper\'{o}, David and Mota, Miguel Angel},
	TITLE = {Forcing consequences of {PFA} together with the continuum
	large},
	JOURNAL = {Trans. Amer. Math. Soc.},
	FJOURNAL = {Transactions of the American Mathematical Society},
	VOLUME = {367},
	YEAR = {2015},
	NUMBER = {9},
	PAGES = {6103--6129},
	ISSN = {0002-9947,1088-6850},
	MRCLASS = {03E35 (03E05 03E57)},
	MRNUMBER = {3356931},
	MRREVIEWER = {Peter\ Holy},
	DOI = {10.1090/S0002-9947-2015-06205-9},
	URL = {https://doi.org/10.1090/S0002-9947-2015-06205-9},
}

@article {Aspero:Mota2015b,
	AUTHOR = {Asper\'{o}, David and Mota, Miguel Angel},
	TITLE = {A generalization of {M}artin's {A}xiom},
	JOURNAL = {Israel J. Math.},
	FJOURNAL = {Israel Journal of Mathematics},
	VOLUME = {210},
	YEAR = {2015},
	NUMBER = {1},
	PAGES = {193--231},
	ISSN = {0021-2172,1565-8511},
	MRCLASS = {03E50 (03E05 03E35 03E57)},
	MRNUMBER = {3430273},
	MRREVIEWER = {Matteo\ Viale},
	DOI = {10.1007/s11856-015-1250-0},
	URL = {https://doi.org/10.1007/s11856-015-1250-0},
}

@book {Barwise1975Admissiblesetsandstructures,
	AUTHOR = {Barwise, Jon},
	TITLE = {Admissible sets and structures},
	SERIES = {Perspectives in Mathematical Logic},
	NOTE = {An approach to definability theory},
	PUBLISHER = {Springer-Verlag, Berlin-New York},
	YEAR = {1975},
	PAGES = {xiii+394},
	MRCLASS = {02H99 (02B25 02F27)},
	MRNUMBER = {424560},
	MRREVIEWER = {Nigel\ J.\ Cutland},
}

@incollection {Burke1998Forcingaxioms,
	AUTHOR = {Burke, Maxim R.},
	TITLE = {Forcing axioms},
	BOOKTITLE = {Set theory ({C}ura\c{c}ao, 1995; {B}arcelona, 1996)},
	PAGES = {1--21},
	PUBLISHER = {Kluwer Acad. Publ., Dordrecht},
	YEAR = {1998},
	ISBN = {0-7923-4905-9},
	MRCLASS = {03E65 (03E05 03E35 03E40 03E50)},
	MRNUMBER = {1601960},
	MRREVIEWER = {Miroslav\ Repick\'{y}},
}

@book {Devlin1973Aspectsofconstructibility,
	AUTHOR = {Devlin, Keith J.},
	TITLE = {Aspects of constructibility},
	SERIES = {Lecture Notes in Mathematics},
	VOLUME = {Vol. 354},
	PUBLISHER = {Springer-Verlag, Berlin-New York},
	YEAR = {1973},
	PAGES = {xii+240},
	MRCLASS = {02K99 (02H05 04-XX)},
	MRNUMBER = {376351},
	MRREVIEWER = {F.\ R.\ Drake},
}

@article{DolinarDzamonja2013,
	title={Forcing $\square_{\omega_1}$ with finite conditions},
	author={Dolinar, Gregor and D{\v{z}}amonja, Mirna},
	journal={Annals of Pure and Applied Logic},
	volume={164},
	number={1},
	pages={49--64},
	year={2013},
	publisher={Elsevier}
}

@article {Dow1998Anintroductiontoapplicationsofelementarysubmodels,
	AUTHOR = {Dow, Alan},
	TITLE = {An introduction to applications of elementary submodels to
	topology},
	JOURNAL = {Topology Proc.},
	FJOURNAL = {Topology Proceedings},
	VOLUME = {13},
	YEAR = {1988},
	NUMBER = {1},
	PAGES = {17--72},
	ISSN = {0146-4124,2331-1290},
	MRCLASS = {54A35 (03C62 03E35 03E55 03E75 54A25)},
	MRNUMBER = {1031969},
	MRREVIEWER = {S.\ R.\ Kogalovski\u{\i}},
}

@incollection {Friedman2006Forcingwithfiniteconditions,
	AUTHOR = {Friedman, Sy-David},
	TITLE = {Forcing with finite conditions},
	BOOKTITLE = {Set theory},
	SERIES = {Trends Math.},
	PAGES = {285--295},
	PUBLISHER = {Birkh\"auser, Basel},
	YEAR = {2006},
	ISBN = {978-3-7643-7691-8; 3-7643-7691-0},
	MRCLASS = {03E40 (03E05)},
	MRNUMBER = {2267153},
	MRREVIEWER = {Juris\ Stepr\=ans},
	DOI = {10.1007/3-7643-7692-9\_10},
	URL = {https://doi.org/10.1007/3-7643-7692-9_10},
}

@phdthesis{Gallart2024thesis,
	author={Gallart, Curial},
	title={Side conditions of models of two types and high forcing axioms},
	school={School of mathematics, University of East Anglia},
	year={2024},
	type={PhD thesis}
}

@article {GitikMagidor2016SPFAbyfiniteconditions,
	AUTHOR = {Gitik, Moti and Magidor, Menachem},
	TITLE = {S{PFA} by finite conditions},
	JOURNAL = {Arch. Math. Logic},
	FJOURNAL = {Archive for Mathematical Logic},
	VOLUME = {55},
	YEAR = {2016},
	NUMBER = {5-6},
	PAGES = {649--661},
	ISSN = {0933-5846,1432-0665},
	MRCLASS = {03E35 (03E40 03E50 03E55 03E57)},
	MRNUMBER = {3523647},
	MRREVIEWER = {Andrzej\ Ros\l anowski},
	DOI = {10.1007/s00153-016-0485-8},
	URL = {https://doi.org/10.1007/s00153-016-0485-8},
}

@book {Jech2003Settheory,
AUTHOR = {Jech, Thomas},
TITLE = {Set theory},
SERIES = {Springer Monographs in Mathematics},
EDITION = {millennium},
PUBLISHER = {Springer-Verlag, Berlin},
YEAR = {2003},
PAGES = {xiv+769},
ISBN = {3-540-44085-2},
MRCLASS = {03Exx (03-01 03-02)},
MRNUMBER = {1940513},
MRREVIEWER = {Eva\ Coplakova},
}

@book {JustWeese1997DiscoveringmodernsetthoryII,
AUTHOR = {Just, Winfried and Weese, Martin},
TITLE = {Discovering modern set theory. {II}},
SERIES = {Graduate Studies in Mathematics},
VOLUME = {18},
NOTE = {Set-theoretic tools for every mathematician},
PUBLISHER = {American Mathematical Society, Providence, RI},
YEAR = {1997},
PAGES = {xiv+224},
ISBN = {0-8218-0528-2},
MRCLASS = {03-01 (03E05 04-01 04A20)},
MRNUMBER = {1474727},
MRREVIEWER = {J.\ M.\ Henle},
DOI = {10.1090/gsm/018},
URL = {https://doi.org/10.1090/gsm/018},
}

@incollection {Kanamori1983Morassesincombinatorialsettheory,
	AUTHOR = {Kanamori, Akihiro},
	TITLE = {Morasses in combinatorial set theory},
	BOOKTITLE = {Surveys in set theory},
	SERIES = {London Math. Soc. Lecture Note Ser.},
	VOLUME = {87},
	PAGES = {167--196},
	PUBLISHER = {Cambridge Univ. Press, Cambridge},
	YEAR = {1983},
	ISBN = {0-521-27733-7},
	MRCLASS = {03E05 (03-02 03E45)},
	MRNUMBER = {823780},
	MRREVIEWER = {F.\ R.\ Drake},
	DOI = {10.1017/CBO9780511758867.007},
	URL = {https://doi.org/10.1017/CBO9780511758867.007},
}

@article{KoepkeMartinez1995superatomicbooleanalgebrasconstructedfrom,
	title={Superatomic {B}oolean algebras constructed from morasses},
	author={Koepke, Peter and Martinez, Juan Carlos},
	journal={The Journal of Symbolic Logic},
	volume={60},
	number={3},
	pages={940--951},
	year={1995},
	publisher={Cambridge University Press}
}

@article {Koszmider2000onstrongchainsofuncoutablefunctions,
	AUTHOR = {Koszmider, Piotr},
	TITLE = {On strong chains of uncountable functions},
	JOURNAL = {Israel J. Math.},
	FJOURNAL = {Israel Journal of Mathematics},
	VOLUME = {118},
	YEAR = {2000},
	PAGES = {289--315},
	ISSN = {0021-2172,1565-8511},
	MRCLASS = {03E35 (03E05 03E40 03E45)},
	MRNUMBER = {1776085},
	MRREVIEWER = {Peter\ Vojt\'a\v s},
	DOI = {10.1007/BF02803525},
	URL = {https://doi.org/10.1007/BF02803525},
}

@book {Kunen2009Thefoudnationsofmathematics,
	AUTHOR = {Kunen, Kenneth},
	TITLE = {The foundations of mathematics},
	SERIES = {Studies in Logic (London)},
	VOLUME = {19},
	NOTE = {Mathematical Logic and Foundations},
	PUBLISHER = {College Publications, London},
	YEAR = {2009},
	PAGES = {viii+251},
	ISBN = {978-1-904987-14-7},
	MRCLASS = {03-01},
	MRNUMBER = {2667506},
}

@book {Kunen2011Settheory,
AUTHOR = {Kunen, Kenneth},
TITLE = {Set theory},
SERIES = {Studies in Logic (London)},
VOLUME = {34},
PUBLISHER = {College Publications, London},
YEAR = {2011},
PAGES = {viii+401},
ISBN = {978-1-84890-050-9},
MRCLASS = {03-01 (03E05 03E10 03E35 03E40 03E45 03E50)},
MRNUMBER = {2905394},
MRREVIEWER = {Klaas\ Pieter\ Hart},
}

@article {KuzeljevicTodorcevic2017Forcingwithmatrices,
AUTHOR = {Kuzeljevi{\'c}, Bori{\v{s}}a and Todor{\v{c}}evi{\'c}, Stevo},
TITLE = {Forcing with matrices of countable elementary submodels},
JOURNAL = {Proc. Amer. Math. Soc.},
FJOURNAL = {Proceedings of the American Mathematical Society},
VOLUME = {145},
YEAR = {2017},
NUMBER = {5},
PAGES = {2211--2222},
ISSN = {0002-9939,1088-6826},
MRCLASS = {03E40 (03E05 03E35)},
MRNUMBER = {3611332},
MRREVIEWER = {Andrzej\ Ros\l anowski},
DOI = {10.1090/proc/13133},
URL = {https://doi.org/10.1090/proc/13133},
}

@article {Mitchell2005Addingclubsofomega2withfiniteforcing,
AUTHOR = {Mitchell, William J.},
TITLE = {Adding closed unbounded subsets of {$\omega_2$} with finite
forcing},
JOURNAL = {Notre Dame J. Formal Logic},
FJOURNAL = {Notre Dame Journal of Formal Logic},
VOLUME = {46},
YEAR = {2005},
NUMBER = {3},
PAGES = {357--371},
ISSN = {0029-4527,1939-0726},
MRCLASS = {03E35 (03E04)},
MRNUMBER = {2162106},
DOI = {10.1305/ndjfl/1125409334},
URL = {https://doi.org/10.1305/ndjfl/1125409334},
}

@article{Mitchell2009I[omega2]canbethenonstationaryidealoncofomega1,
	title={${I}[\omega_2]$ can be the nonstationary ideal on $\operatorname{Cof}(\omega_1)$},
	author={Mitchell, William J.},
	journal={Trans. Amer. Math. Soc},
	volume={361},
	number={2},
	pages={561--601},
	year={2009}
}

@misc{Miyamoto2012AnoteonageneralizedMAfortheNorwichposets,
	title={A note on a generalized Martin’s Axiom for the Norwich posets},
	author={Miyamoto, Tadatoshi},
	date={2019},
	howpublished={Unpublished note}
}

@article{Miyamoto2014Matricesofisomorphicmodelsandmorasslikestructures,
title={Matrices of isomorphic models and morass-like structures (Reflection principles and set theory of large cardinals)},
author={Miyamoto, Tadatoshi},
journal={RIMS Kôkyûroku},
volume={1895},
pages={79--102},
year={2014},
publisher={Research Institute for Mathematical Sciences, Kyoto University}
}

@article{Miyamoto2015Squaresbymatriceswithcoherentsequences,
	title={Squares by Matrices with Coherent Sequences},
	author={Miyamoto, Tadatoshi},
	journal={RIMS Kôkyûroku},
	volume={1949},
	pages={62--72},
	year={2015},
	publisher={Research Institute for Mathematical Sciences, Kyoto University}
}

@misc{Miyamoto2019AnoteonPffcandWTCGS12,
	title={A note on ${P}_{ff}^c$ and {WTCG}$({S}_1^2)$},
	author={Miyamoto, Tadatoshi},
	date={2019},
	howpublished={Unpublished note}
}

@article{Miyamoto2022Forcingaclubbyageneralizedfastfunction,
	title={Forcing a Club by a Generalized Fast Function},
	author={Miyamoto, Tadatoshi},
	journal={Journal of the Nanzan Academic Society Humanities and Natural Sciences},
	number={23},
	pages={183--191},
	year={2022},
	publisher={Nanzan University}
}

@misc{Miyamoto2023Astronglysigmaclosedposetthatforcesasimplifiedmorass,
	title={A strongly $\sigma$-closed poset that forces a simplified $(\omega_1,1)$-morass},
	author={Miyamoto, Tadatoshi},
	date={2023},
	howpublished={Unpublished note}
}

@article{Miyamoto2023Negatingapartitionrelationbyafamilyofsimplified,
	title={Negating a partition relation by a family of simplified morasses},
	author={Miyamoto, Tadatoshi},
	journal={RIMS Kôkyûroku},
	volume={2261},
	pages={78--87},
	year={2023},
	publisher={Research Institute for Mathematical Sciences, Kyoto University}
}

@article {Mohammadpour2023SpecialisingtreeswithsmallapproximationsI,
	AUTHOR = {Mohammadpour, Rahman},
	TITLE = {Specialising trees with small approximations {I}},
	JOURNAL = {J. Symb. Log.},
	FJOURNAL = {The Journal of Symbolic Logic},
	VOLUME = {88},
	YEAR = {2023},
	NUMBER = {2},
	PAGES = {640--663},
	ISSN = {0022-4812,1943-5886},
	MRCLASS = {03E05 (03E35 03E57)},
	MRNUMBER = {4594294},
	MRREVIEWER = {Yair\ Hayut},
	DOI = {10.1017/jsl.2022.24},
	URL = {https://doi.org/10.1017/jsl.2022.24},
}

@article{Neeman2014Forcingwithsequencesofmodelsoftwotypes,
AUTHOR = {Neeman, Itay},
TITLE = {Forcing with sequences of models of two types},
JOURNAL = {Notre Dame J. Form. Log.},
FJOURNAL = {Notre Dame Journal of Formal Logic},
VOLUME = {55},
YEAR = {2014},
NUMBER = {2},
PAGES = {265--298},
ISSN = {0029-4527,1939-0726},
MRCLASS = {03E35},
MRNUMBER = {3201836},
MRREVIEWER = {John\ Krueger},
DOI = {10.1215/00294527-2420666},
URL = {https://doi.org/10.1215/00294527-2420666},
}

@article{Neeman2017Twoapplicationsoffinitesideconditionsatomega2,
AUTHOR = {Neeman, Itay},
TITLE = {Two applications of finite side conditions at $\omega_2$},
JOURNAL = {Arch. Math. Logic},
FJOURNAL = {Archive for Mathematical Logic},
VOLUME = {56},
YEAR = {2017},
NUMBER = {7-8},
PAGES = {983--1036},
ISSN = {0933-5846,1432-0665},
MRCLASS = {03E35 (03E05)},
MRNUMBER = {3696074},
MRREVIEWER = {Andrzej\ Ros\l anowski},
DOI = {10.1007/s00153-017-0550-y},
URL = {https://doi.org/10.1007/s00153-017-0550-y},
}

@unpublished{NeemanSlides1,
	title        = "Higher analog of the {PFA}",
	author       = "{Neeman, Itay}",
	note = "Slides presented at the Fields Institute, https://www.math.ucla.edu/~ineeman/",
	year         = 2012
}

@unpublished{NeemanSlides2,
	title        = "Higher analogues of properness",
	author       = "{Neeman, Itay}",
	note = "Slides presented at BIRS, https://www.math.ucla.edu/~ineeman/",
	year         = 2013
}

@book {Shelah2017Properandimproperforcing,
AUTHOR = {Shelah, Saharon},
TITLE = {Proper and improper forcing},
SERIES = {Perspectives in Mathematical Logic},
EDITION = {Second},
PUBLISHER = {Springer-Verlag, Berlin},
YEAR = {1998},
PAGES = {xlviii+1020},
ISBN = {3-540-51700-6},
MRCLASS = {03-02 (03E05 03E35 03E40 03E45 03E50)},
MRNUMBER = {1623206},
MRREVIEWER = {A.\ Kanamori},
DOI = {10.1007/978-3-662-12831-2},
URL = {https://doi.org/10.1007/978-3-662-12831-2},
}

@article{ShelahStanley1982SforcingAblackboxtheoremformorasses,
	title={S-forcing, I. A “black-box” theorem for morasses, with applications to super-Souslin trees},
	author={Shelah, Saharon and Stanley, Lee},
	journal={Israel Journal of Mathematics},
	volume={43},
	number={3},
	pages={185--224},
	year={1982},
	publisher={Springer}
}

@inbook {Todorcevic1984AnoteonthePFA,
AUTHOR = {Todor\v{c}evi\'c, Stevo},
TITLE = {A note on the proper forcing axiom},
BOOKTITLE = {Axiomatic set theory ({B}oulder, {C}olo., 1983)},
SERIES = {Contemp. Math.},
VOLUME = {31},
PAGES = {209--218},
PUBLISHER = {Amer. Math. Soc., Providence, RI},
YEAR = {1984},
ISBN = {0-8218-5026-1},
MRCLASS = {03E50 (03E35 03E40)},
MRNUMBER = {763902},
MRREVIEWER = {James\ Baumgartner},
DOI = {10.1090/conm/031/763902},
URL = {https://doi.org/10.1090/conm/031/763902},
}

@article{Todorcevic1985directedsets,
	title={Directed sets and cofinal types},
	author={Todor{\v{c}}evi{\'c}, Stevo},
	journal={Transactions of the American Mathematical Society},
	volume={290},
	number={2},
	pages={711--723},
	year={1985}
}

@article{Todorcevic2014notesonforcingaxioms,
	title={Notes on forcing axioms, volume 26 of Lecture Notes Series},
	author={Todor{\v{c}}evi{\'c}, Stevo},
	journal={Institute for Mathematical Sciences. National University of Singapore. World Scientific Publishing Co., Pte. Ltd., Hackensack, NJ: Edited and with a foreword by Chong, C., Feng, Q., Yang, Y., Slaman, TA, Woodin, WH},
	year={2014}
}

@article{Velikovic2014iterationofsemiproper,
	title={Iteration of Semiproper Forcing Revisited},
	author={Veli{\v{c}}kovi{\'c}, Boban},
	journal={arXiv preprint arXiv:1410.5095 [math.LO]},
	url={https://arxiv.org/abs/1410.5095},
	year={2014}
}

@incollection {Velickovic:Venturi2011properforcingremastered,
	AUTHOR = {Veli\v{c}kovi\'{c}, Boban and Venturi, Giorgio},
	TITLE = {Proper forcing remastered},
	BOOKTITLE = {Appalachian set theory 2006--2012},
	SERIES = {London Math. Soc. Lecture Note Ser.},
	VOLUME = {406},
	PAGES = {331--362},
	PUBLISHER = {Cambridge Univ. Press, Cambridge},
	YEAR = {2013},
	ISBN = {978-1-107-60850-4; 978-1-107-60850-1},
	MRCLASS = {03E57 (03E55)},
	MRNUMBER = {3821634},
}

@article{Velleman1982morassesdiamondandforcing,
	title={Morasses, diamond, and forcing},
	author={Velleman, Daniel J.},
	journal={Annals of Mathematical Logic},
	volume={23},
	number={2-3},
	pages={199--281},
	year={1982},
	publisher={North-Holland}
}

@article{Velleman1984simplifiedmorasses,
	title={Simplified morasses},
	author={Velleman, Daniel J.},
	journal={The Journal of Symbolic Logic},
	volume={49},
	number={1},
	pages={257--271},
	year={1984},
	publisher={Cambridge University Press}
}

@article {Zapletal1997Stronglyalmostdisjointfunctions,
	AUTHOR = {Zapletal, Jind\v{r}ich},
	TITLE = {Strongly almost disjoint functions},
	JOURNAL = {Israel J. Math.},
	FJOURNAL = {Israel Journal of Mathematics},
	VOLUME = {97},
	YEAR = {1997},
	PAGES = {101--111},
	ISSN = {0021-2172,1565-8511},
	MRCLASS = {03E35 (03E05)},
	MRNUMBER = {1441241},
	MRREVIEWER = {Douglas\ R.\ Burke},
	DOI = {10.1007/BF02774029},
	URL = {https://doi.org/10.1007/BF02774029},
}

\end{document}